  \documentclass[11pt,a4paper,twoside]{article}

  \usepackage[margin=3.15cm]{geometry}
  \usepackage{fourier}
  \usepackage{amsmath,amssymb,amsthm,amsfonts}
  \usepackage{bm}
  \usepackage{mathabx}
  \usepackage{xspace,xcolor}
  \usepackage{parskip,fancyhdr,url}
  
  \usepackage{enumitem}
  \usepackage{todonotes}
  \usepackage[
     breaklinks,colorlinks,
     citecolor=blue,linkcolor=blue,urlcolor=teal
  ]{hyperref}
  \usepackage{tikz, tikz-cd}
  \usetikzlibrary{
    arrows, calc, matrix, decorations.pathmorphing, shapes, arrows
    }
  \usepackage{mathtools}

  \newcommand\Rsquigarrow[1]{%
\mathrel{%
\begin{tikzpicture}[baseline= {( $ (current bounding box.south) + (0,-0.5ex) $ )}]
  \node[inner sep=.5ex] (a) {$\scriptstyle #1$};
  \path[draw,implies-,double distance between line centers=1.5pt,decorate,
    decoration={snake,amplitude=0.7pt,segment length=1.2mm,pre=lineto,
    pre   length=4pt}]
    (a.south east) -- (a.south west);
\end{tikzpicture}}%
}

  \begingroup
    \makeatletter
    \@for\theoremstyle:=definition,remark,plain\do{%
      \expandafter\g@addto@macro\csname th@\theoremstyle\endcsname{%
        \addtolength\thm@preskip\parskip
        }%
      }
  \endgroup

  \makeatletter
    \def\tagform@#1{\maketag@@@{%
     \textbf{(\ignorespaces#1\unskip\@@italiccorr)}}}%
     \renewcommand{\eqref}[1]{\textup{\maketag@@@{(\ignorespaces%
          {\ref{#1}}\unskip\@@italiccorr)}}}
  \makeatother

  \newcommand\address[1]{}
  \newcommand\email[1]{}
  \newcommand\dedicatory[1]{}

  \theoremstyle{plain}
  \newtheorem{theorem}[equation]{Theorem}
  \newtheorem{proposition}[equation]{Proposition}
  \newtheorem{corollary}[equation]{Corollary}
  \newtheorem{lemma}[equation]{Lemma}

  \theoremstyle{definition}
  \newtheorem{definition}[equation]{Definition}
  \newtheorem{remark}[equation]{Remark}
  \newtheorem{example}[equation]{Example}

  \theoremstyle{axiom}

  \numberwithin{equation}{section}

  \theoremstyle{problem}

  \newcommand{\bsf}[1]{\ensuremath{\bm{\mathnormal{\mathsf{#1}}}}\xspace}

  \newcommand{\cat}{\ensuremath{\bsf{C}\xspace}}

  \newcommand{\CLip}{\ensuremath{\bsf{CLip}\xspace}}
  \newcommand{\QGeod}{\ensuremath{\bsf{QGeod}\xspace}}
  \newcommand{\CGeod}{\ensuremath{\bsf{CGeod}\xspace}}

  \newcommand{\Set}{\ensuremath{\bsf{Set}\xspace}}

  \newcommand{\Born}{\ensuremath{\bsf{Coarse}\xspace}}

  \newcommand{\barcat}{\overline{\cat}}
  \newcommand{\barCL}{\overline{\CLip}}
  \newcommand{\barBorn}{\overline{\Born}}
  \newcommand{\barQG}{\overline{\QGeod}}
  \newcommand{\barCG}{\overline{\CGeod}}

  \newcommand{\Gr}{{\ensuremath{\mathsf{G}\xspace}}} 
  \newcommand{\Grhat}{\ensuremath{\widehat{\Gr}\xspace}} 
  \newcommand{\Rc}{{\ensuremath{\mathsf{R}\xspace}}} 

   \newcommand{\F}{\ensuremath{\mathcal{F}}\xspace}

  \DeclareMathOperator{\Obj}{Obj} 
  \DeclareMathOperator{\Mor}{Mor} 
  \DeclareMathOperator{\dom}{dom} 
  \DeclareMathOperator{\cod}{cod} 
  \DeclareMathOperator{\Cone}{Cone} 

  \DeclareMathOperator{\Tuple}{Tuple} 

  \newcommand{\limit}{\ensuremath{\varprojlim\xspace}}
  
  \newcommand{\colimit}{\ensuremath{\varinjlim\xspace}}

  \newcommand{\uc}{u.c.\xspace}

  \newcommand{\uclim}{\mathop{\ooalign{$\lim$\cr
  \hidewidth\raise-1.10ex\hbox{$\leftsquigarrow\mkern-0.5mu$}\cr}}}

  \DeclareMathOperator{\Rips}{Rips}

  \DeclareMathOperator{\Eq}{Eq} 

  \newcommand{\st}{\ensuremath{\,\, \colon \,\,}\xspace}
  \newcommand{\from}{\ensuremath{\colon \thinspace}\xspace}
  \newcommand{\tto}{\ensuremath{\Rightarrow}\xspace}
  \newcommand{\sqto}{\ensuremath{\rightsquigarrow}\xspace}
  \newcommand{\sqtto}{\ensuremath{\Rsquigarrow{\;\;\;\;}}\xspace}

  \newcommand{\floor}[1]{\ensuremath{\lfloor {#1} \rfloor}\xspace}
  \newcommand{\floorC}[1]{\floor{{#1}}_{\cat}}

  \DeclareMathOperator{\diam}{diam}

  \newcommand{\N}{\ensuremath{\mathbb{N}}\xspace}

        \newcommand{\C}{\ensuremath{\mathcal{C}}\xspace}
        
	\newcommand{\X}{\ensuremath{\mathcal{X}}\xspace}
	\newcommand{\Y}{\ensuremath{\mathcal{Y}}\xspace}

	\newcommand{\calS}{\ensuremath{\mathcal{S}}\xspace}
       \newcommand{\calR}{\ensuremath{\mathcal{R}}\xspace}

  \newcommand{\sJ}{\mathsf{J}}

  \newcommand{\proofof}[1]{\hfill\newline\noindent\emph{Proof of {#1}.} }

  \newcommand{\param}{{\mathchoice{\mkern1mu\mbox{\raise2.2pt\hbox{$
  \centerdot$}}
  \mkern1mu}{\mkern1mu\mbox{\raise2.2pt\hbox{$\centerdot$}}\mkern1mu}{
  \mkern1.5mu\centerdot\mkern1.5mu}{\mkern1.5mu\centerdot\mkern1.5mu}}}

  \fancypagestyle{plain}

  \AtEndDocument{\bigskip{\footnotesize%
  Department of Pure Mathematics, Xi'an Jiaotong--Liverpool University, 111 Ren'ai Road, Suzhou Industrial Park, Suzhou, Jiangsu, 215123, China
  \par \texttt{robert.tang@xjtlu.edu.cn}}}

\begin{document}

          \title    {The metric Rips filtration, universal quasigeodesic cones, and hierarchically hyperbolic spaces}
    \author   {Robert Tang} \date{20 November 2025}

  \maketitle \thispagestyle{empty}

  \begin{abstract}

  We introduce a flexible, categorical framework for large-scale geometry that clarifies basic behaviour of the metric Rips filtration and streamlines some constructions in geometric group theory. The paper has two main parts. First, we develop the theory of the metric Rips filtration and its colimit in natural coarse categories: informally, we characterise when the Rips colimit produces a canonical large‑scale model of a metric space and use this to prove that the quasigeodesic subcategory is closed under colimits in the coarsely Lipschitz category. We also establish adjointness properties of the Rips colimit and use them to characterise extremal metrics and universal morphisms from quasigeodesic sources.

  Second, we apply this machinery to characterise universal quasigeodesic cones via an explicit Rips--Tuple recipe. In the HHS setting this yields a concrete, canonical model of the total space: an HHS is quasi‑isometric to a Rips graph of the space of coarsely consistent tuples in the product of its factor spaces. Moreover, we give a local‑to‑global criterion that promotes uniformly controlled, factorwise retractions to a canonical global hierarchical retraction. Because the approach is based on universal properties and uniformly controlled coarse data rather than inductive constructions, distance formulae, or hierarchy paths, it applies equally well to arbitrary families of metric spaces equipped with pairwise constraints.
  \end{abstract}

  \tableofcontents

  \section{Introduction}

  \subsection*{Metric Rips filtration and its colimit}

  Coarse or large‑scale geometry studies properties of metric spaces that persist after ignoring small‑scale detail. A standard approach in geometric group theory and coarse geometry is to approximate complicated spaces by simpler objects -- typically graphs or hyperbolic spaces -- at varying scales. The metric Rips filtration formalises this: for a metric space $X$, the \emph{Rips graph} at scale $\sigma \geq 0$ is the graph $\Rips_\sigma X$ with vertex set $X$, with $x,y\in X$ declared adjacent whenever $d(x,y) \leq \sigma$. These graphs assemble via the 1‑Lipschitz inclusions $\Rips_\sigma X \hookrightarrow \Rips_\tau X$ for $\sigma \leq \tau$ to form the \emph{metric Rips filtration} $\Rips_* X$.

  The simplicial Rips complex traces back to Vietoris’s work in algebraic topology \cite{Vie27} and Rips' treatment of hyperbolic groups \cite{Gro87, BH99}. Subsequent developments connect the simplicial Rips filtration to manifold approximations up to homotopy \cite{Hau95, Lat01}, the coarse Baum–Connes conjecture \cite{HR95, Yu95}, topological data analysis \cite{ELZ02, EH08, Car09, CdSO14}, and other applications in geometry and topology \cite{CDVV12, Zar22}. Motivated by the natural question of ``what happens as the scale parameter tends to infinity’’ arising in the contexts of persistent homology and the coarse Baum–Connes conjecture, we study the analogous behaviour for the purely metric Rips filtration.

  To formulate a precise notion of ``limit'' for the metric Rips filtration, we work in two categories: the metric coarse category $\barBorn$ (objects are extended metric spaces; morphisms are controlled maps up to closeness) and the coarsely Lipschitz category $\barCL $ (the wide subcategory with coarsely Lipschitz morphisms). Isomorphisms in these categories are, respectively, (closeness classes of) coarse equivalences and quasi‑isometries. The \emph{Rips colimit} of $X$ is the categorical colimit $\colimit \Rips_* X$ taken in either $\barBorn$ or $\barCL$.

  Our first theorem shows that, up to coarse equivalence, the metric Rips filtration recovers every metric space.

   \begin{theorem}[$\barBorn$--Rips colimit]\label{thm:CE-Rips}
    Every metric space is canonically isomorphic to its Rips colimit in $\barBorn$.
   \end{theorem}

    Hence, up to coarse equivalence, any metric space is a canonical colimit of graphs, and the full subcategory of $\barBorn$ comprising the graphs (equivalently coarsely geodesic spaces) is colimit‑dense.

    Quasi‑isometry presents a different picture. We say that $\Rips_* X$ \emph{stabilises} if the inclusions ${\Rips_\sigma X \hookrightarrow \Rips_\tau X}$ are isomorphisms for all sufficiently large $\sigma \leq \tau$ (this notion is independent of working in $\barBorn$ or $\barCL$). It is standard result that $\Rips_* X$ stabilises precisely when $X$ is coarsely geodesic, in which case $\Rips_\sigma X$ realises the Rips colimit for large $\sigma$. We show that stability is in fact necessary for existence of the Rips colimit in the coarsely Lipschitz category.

    \begin{theorem}[$\barCL$--Rips colimit]\label{thm:QI-Rips}
    Let $X$ be a metric space. Then $\colimit \Rips_* X$ exists in $\barCL$ if and only if $\Rips_* X$ stabilises. In particular, any Rips colimit in $\barCL$ is quasigeodesic.
    \end{theorem}

    Thus, unless $X$ is already quasigeodesic (equivalently quasi‑isometric to a graph), its Rips colimit cannot recover $X$ up to quasi‑isometry. One might then try taking colimits of arbitrary small diagrams of graphs; our main structural result shows this yields no new spaces up to quasi-isometry.

    \begin{theorem}[Coclosed subcategory]\label{thm:coclosed}
    The full subcategory $\barQG \subset \barCL$ of quasigeodesic spaces is closed under all colimits.
    \end{theorem}

    Since $\barCL$ admits all finite colimits (see \cite[Propositon 3.4]{Tang-mono}), we obtain:

    \begin{corollary}[Finite cocompleteness]
     The category $\barQG$ admits all finite colimits. \qed
    \end{corollary}

     Theorem \ref{thm:coclosed} makes no assertion about the existence of colimits in $\barQG$; rather, it imposes geometric regularity on colimits when they do exist. Notably, we make no uniformity assumptions on quasigeodesicity constants or on the coarsely Lipschitz maps appearing in the diagrams.

     \begin{remark}
      Our results concerning the coarsely Lipschitz category generalise to wide subcategories of $\barBorn$ where the morphisms admit \emph{dominated controls} (see Definition \ref{def:dominated}). This yields a notion of $\cat$--geodesic spaces, generalising the quasigeodesic spaces.
     \end{remark}

     \subsection*{Adjointness and universal properties}

     A key technical reduction used repeatedly (notably in the proof of Theorem \ref{thm:QI-Rips}) is to pass first to the metric coarse setting, apply Theorem \ref{thm:CE-Rips} to obtain the canonical cocone, and then verify that the induced morphism is coarsely Lipschitz. We now make this passage between the metric coarse and coarsely Lipschitz categories systematic.

     Write $\barCG \subset \barBorn$ for the full subcategory whose objects are the coarsely geodesic spaces; and let {$\barCG \cap \barCL \subset \barCG$} denote the wide subcategory with the same objects but only coarsely Lipschitz morphisms.

    \begin{theorem}[Adjointness]\label{thm:adjoint0}
    The inclusion functor $\barQG \hookrightarrow \barCG$ and the (stable) Rips colimit functor $\colimit \Rips_* \from \barCG \to \barQG$ form an adjoint equivalence. The inclusion $\barQG \hookrightarrow \barCG \cap \barCL$ is left-adjoint to $\colimit \Rips_* \from \barCG \cap \barCL \to \barQG$ but this adjunction is not an equivalence.
    \end{theorem}

    This yields a canonical functorial way to convert coarsely geodesic spaces to quasigeodesic spaces -- take the Rips graph at a stable scale -- so universal constructions transfer naturally. Conversely, the result shows that enlarging the class of spaces from quasigeodesic to coarsely geodesic forces one to allow the more flexible controlled maps (rather than only coarsely Lipschitz maps) if one wishes to retain good categorical behaviour.

    Adjointness yields concrete universal properties for the Rips colimit. Consider the problem of finding a ``best approximation'' to a given metric space $X$ using maps from arbitrary quasigeodesic spaces (or graphs). Call a coarsely Lipschitz map $f \from W \to X$ from a quasigeodesic space $W$ a \emph{universal morphism from} $\barQG$ if for every coarsely Lipschitz $f' \from W' \to X$ with $W'$ quasigeodesic, there exists a unique (up to closeness) coarsely Lipschitz $h \from W' \to W$ with $f'$ close to $fh$. We show that the (stable) Rips graph of $X$ solves this problem, serving as a ``universal receiver'' for all coarsely Lipschitz maps from quasigeodesic spaces to $X$.

    \begin{theorem}[Universal morphisms]\label{thm:univ-morph} A metric space $X$ admits a universal morphism from $\barQG$ if and only if $X$ is coarsely geodesic. When this holds, the universal morphism is realised by the underlying identity $\Rips_\sigma X \to X$ for $\sigma$ sufficiently large.
    \end{theorem}

    We apply this to obtain an intrinsic extremal characterisation of quasigeodesic spaces inspired by Rosendal’s notion of a maximal compatible left-invariant metric for topological groups \cite{Ros22, Ros23}. Rosendal established several equivalent conditions for the existence of such a maximal metric on Polish groups; in particular, if one exists then it must be quasigeodesic. This endows such groups with a canonical metric (up to quasi-isometry) which generalises the word metric on finitely generated groups.

    Given metrics $d_1, d_2$ on a set $X$, write $d_1 \prec d_2$ if the underlying identity $(X,d_2) \to (X,d_1)$ is coarsely Lipschitz; this yields a preorder on the set of metrics on $X$.

    \begin{theorem}[Extremal metrics]\label{thm:extremal} For a metric space $(X,d)$, the following are equivalent:
        \begin{enumerate}
            \item $(X,d)$ is coarsely geodesic,
            \item there exists a $\prec$--maximal metric within the coarse equivalence class of $d$, and
            \item there exists a $\prec$--minimal quasigeodesic metric $d'$ on $X$ with $d \prec d'$.
        \end{enumerate}
        Moreover, when such extremal metrics exist any quasigeodesic metric $d'$ coarsely equivalent to d realises both the maximal and minimal properties.
    \end{theorem}

    Thus, the stable Rips graph of $X$ (if it exists) simultaneously attains the maximal and minimal metrics. The maximality in (2) mirrors Rosendal’s concept; while the minimality in (3) is a large‑scale notion which differs from Rosendal’s small‑scale minimal metric \cite{Ros18}.

    \subsection*{Universal quasigeodesic cones and the Rips-Tuple recipe}

    Our main application of the Rips filtration and its universal properties is a concrete procedure for constructing universal quasigeodesic cones, with an eye toward hierarchically hyperbolic spaces.

    Informally, a universal quasigeodesic cone is a quasigeodesic space that universally receives coarsely Lipschitz maps from all quasigeodesic spaces constrained to factor through a given family of metric spaces. We treat two regimes: the non‑uniform setting (morphisms considered up to closeness, with no uniform control) and a uniformly controlled setting (where upper controls and coarse commutativity satisfy uniform bounds).

    Let $D \from \sJ \to \barBorn$ be a diagram indexed by a small category $\sJ$ (the objects of $D$ need not be quasigeodesic). Call a cone over $D$ \emph{quasigeodesic} when its apex is quasigeodesic, and call such a cone \emph{universal} if every quasigeodesic cone over $D$ factors through it uniquely (up to closeness). We show that, in the non-uniform regime, one constructs a universal quasigeodesic cone by first taking the limit of $D$ in $\barBorn$ and then passing to the stable Rips colimit of that limit; if either step fails, no universal quasigeodesic cone exists.

    \begin{theorem}[Universal quasigeodesic cones]\label{thm:univ-qgc}
    A diagram $D \from \sJ \to \barBorn$ admits a universal quasigeodesic cone if and only if $\limit D$ exists in $\barBorn$ and is coarsely geodesic. In that case, $\Rips_\sigma (\limit D)$ for $\sigma$ large realises the universal quasigeodesic cone.
    \end{theorem}

    For finite diagrams this gives an especially concrete recipe. Choose a representative map $\beta_\phi \in D_\phi \from D_i \to D_j$ for each arrow $\phi: i \to j$ in $\sJ$. For $\kappa \geq 0$, let the \emph{$\kappa$--consistent tuple space} be
    \[\Tuple_\kappa = \left\{ x \in \prod_j D_j \st d_{D_j}(x_j, \beta_\phi x_i) \leq \kappa \textrm{ for every arrow } \phi \from i \to j \textrm{ of } D\right\}\]
    equipped with the $\ell^\infty$--metric.
    The $\Tuple_\kappa$ for all $\kappa \geq 0$ form an increasing filtration $\Tuple_*$. This filtration stabilises precisely when the inclusions are eventually coarsely surjective; in other words, by some finite threshold, allowing larger error produces only perturbations of previously existing solutions.

    \begin{theorem}[Rips--Tuple recipe]\label{thm:recipe}
    For a finite diagram $D \from \sJ \to \barBorn$, if $\Tuple_*$ stabilises to a coarsely geodesic space then for some $\sigma, \kappa \geq 0$, the graph $\Rips_\sigma \Tuple_\kappa$ realises the universal quasigeodesic cone for $D$; otherwise, no universal quasigeodesic cone exists.
    \end{theorem}

    This recipe extends to arbitrary diagrams provided one imposes mild uniformity conditions on controls and coarse commutativity (see Theorem \ref{thm:univ-geod-cone}). This uniformly controlled regime shall be treated formally using the notions of \emph{uniformly controlled diagrams} and \emph{cones} in Section \ref{sec:uc-cone}. In practice, these uniformity hypotheses capture the usual controlled data that arise in coarse constructions. Therefore, the tradeoff -- greater generality at the cost of uniform bounds -- is practical for applications arising in geometric group theory.

    \subsection*{Hierarchically hyperbolic spaces}

    Finally, we apply the uniformly controlled quasigeodesic cone formalism to hierarchically hyperbolic spaces (HHS). HHS were introduced by Behrstock–Hagen–Sisto \cite{HHS1, HHS2, HHSsurvey} to capture large‑scale geometry governed by uniformly coarsely Lipschitz projections to a family of uniformly hyperbolic factor spaces subject to pairwise consistency constraints. This framework builds upon the Masur–Minsky machinery for mapping class groups and Teichmüller space via subsurface projections to curve graphs \cite{MM99, MM00, Bro03, Beh06, Raf07, BM08, RS09, BKMM12, MS13, Raf14, BBF15, Bow16, EMR17} and encompasses many examples: CAT(0)--cube complexes, right-angled Artin groups, 3--manifold groups (without \emph{Nil} or \emph{Sol} components) \cite{HHS2}; finitely generated Veech groups \cite{DHS17, Tang-affine} and their surface extensions \cite{DDLS24}; and many others \cite{HS20, BR22, Vok22}. The HHS toolkit has produced a range of structural results, including distance formula and hierarchical paths \cite{HHS2}; bounds on asymptotic dimension \cite{HHS-asdim}; the strong Tits alternative and boundary structure \cite{DHS17}; characterisation of quasiflats \cite{HHS-quasiflat}; and uniform exponential growth \cite{ANS+24}.

    In Section \ref{sec:hhs}, we encode HHS data -- the quasigeodesic \emph{total space} $X$, the uniformly coarsely Lipschitz projection family $\{\lambda_U \from \X \to \C U\}_{U \in \calS}$ to uniformly hyperbolic factor spaces, and the pairwise constraints (Bounded Geodesic Image, Behrstock inequalities) -- as a uniformly controlled diagram and construct the corresponding uniformly controlled quasigeodesic cone. Using only the Uniqueness axiom and the Realisation theorem from HHS theory, we show that the total space $\X$ is universal among all uniformly controlled quasigeodesic cones over this data. The precise statement appears as Theorem \ref{thm:hhs-uccone}; informally:

    \begin{theorem}[HHS as universal quasigeodesic cone]\label{thm:hhs-universal}
    Any hierarchically hyperbolic space is the universal uniformly controlled quasigeodesic cone over its associated uniformly controlled diagram as defined in Section \ref{sec:hhs}.
    \end{theorem}

    A raison d'\^etre of the HHS axioms is that they allow the large‑scale geometry of $\X$ to be recovered from its projection data -- effectively to replicate the Masur--Minsky paradigm in a far broader setting. The prevailing approach to this recovery has required substantial technical machinery (inductively constructing quasigeodesics in each factor to obtain hierarchy paths, estimations of projection distances, consistency checks, keeping track of error bounds, etc.). Taking the universal  property as our starting point and combining it with the Rips--Tuple recipe, we obtain a considerably more streamlined route to the same conclusions.

    \begin{theorem}[HHS via Rips--Tuple]\label{thm:hhs-image}
     Let $\Y \subseteq \prod_U \C U$ be the space of coarsely consistent tuples equipped with the $\ell^\infty$--metric. Then $\prod_U \lambda_U$ induces a quasi-isometry $\X \to \Rips_\sigma \Y$ for $\sigma \geq 0$ large.
    \end{theorem}

    This statement does not itself reprove the distance formula or construct hierarchy paths, but it does show that the distance formula is really estimating distances in $\Rips_\sigma \Y$ (for large $\sigma$). In practice, many HHS arguments are already carried out in the auxiliary space of coarsely consistent tuples in the product of factor spaces rather than in the total space; the theorem says that, after an appropriate quasigeodesic remetrisation, that tuple space is canonically identified up to quasi‑isometry with the HHS.

    Viewed categorically, this decouples the total space from the relative HHS data (factor spaces, projection maps, consistency constraints) and admits a wider range of quasigeodesic cones built over the same family of factor spaces (for example, curve graphs) that need not satisfy every HHS axiom. This perspective complements forthcoming work of Hagen--Petyt--Zalloum on reconstructing total spaces from relative data; we will address the reconstruction problem via universal constructions and the Rips–Tuple recipe in future work.

    Another application of the universal property is to the construction of globally defined maps or subspaces compatible with factorwise information. Given a family of retractions on each factor space (e.g. nearest-point projections to chosen geodesics or approximating trees in each hyperbolic factor), we wish to assemble these to a retraction of the ambient HHS onto a subspace (eg. a hierarchical hull) compatible with the prescibed maps on each factor. The following criterion gives a clean solution under mild uniformity hypotheses -- see Corollary \ref{cor:hier-ret-subspace} for a precise formulation. The power lies in replacing inductive arguments up the hierarchy (where one must build objects level-by-level and perform careful bookkeeping of the accumulated errors) by a local-to-global principle: check a uniform compatibility condition on each pair of factors, and the global object follows.

    \begin{theorem}[Pairwise compatible retractions]\label{thm:hier-retractions}
    Let $(\X, \{\C U \}_{U \in \calS})$ be an HHS. Assume
    \begin{itemize}
        \item $\alpha_U \from \C U \to \C' U$ for $U \in \calS$ is a family of uniformly coarsely Lipschitz retractions to subspaces $\C' U \subseteq \C U$, and
        \item for all distinct $U,V\in\calS$, whenever $(x_U, x_V)\in \C U \times \C V$ is consistent then $(\alpha_U x_U, \alpha_V x_V)$ is consistent up to uniform error.
    \end{itemize}
    Then there exists a canonical (up to closeness) coarsely Lipschitz retraction $\tilde \alpha \from \X \to \X'$ onto a subspace $\X' \subseteq \X$ induced by the $\alpha_U$. Moreover, $\X'$ is a total space over $\{\C' U \}_{U \in \calS}$ with the pairwise consistency constraints inherited from the original HHS structure.
    \end{theorem}

    This criterion is really an instance of a more general statement about retractions between uniformly controlled diagrams and their universal quasigeodesic cones (Theorem \ref{thm:uc-retraction}); again, the proof requires nothing from the HHS structure beyond the universal property. Consequently our arguments apply equally to any pairwise constrained family of metric spaces equipped with uniformly controlled, pairwise compatible families of maps (see Section \ref{sec:hier-retract}).

    In closing, our perspective goes beyond the HHS setting: it applies to any diagram of metric spaces equipped with uniformly controlled maps and coarse commutativity relations, characterises when canonical global objects exist, and furnishes explicit constructions via the Rips--Tuple recipe. This level of abstraction produces a compact toolkit of general lemmas and equivalences that can replace many bespoke inductive arguments and support more diagrammatic reasoning in coarse geometry.
    Although it is a common rule of thumb that ``things ought to work once one has uniform control,'' a concise, systematic account isolating the precise uniformity hypotheses has not appeared to the best of our knowledge. This paper provides such an account: it pinpoints the required hypotheses, streamlines standard HHS reconstructions, increases flexibility (e.g.~allows more general retractions), and reduces the ad hoc bookkeeping of constants and error estimates in routine applications.

  \subsection*{Organisation}

  Section \ref{sec:background} collects background on coarse geometry and the categorical language used throughout. In Section \ref{sec:Rips}, we develop  foundational results about the metric Rips filtration and its behaviour in our categories and characterise existence of Rips colimits (Theorems \ref{thm:CE-Rips}, \ref{thm:QI-Rips}). In particular we establish that the quasigeodesic subcategory is colimit‑closed in the coarse Lipschitz setting (Theorem \ref{thm:coclosed}) and introduce the notion of $\cat$--geodesic spaces, which generalise coarsely geodesic and quasigeodesic notions. Section \ref{sec:adjoint} then focusses on adjointness properties of the Rips colimit and applications to universal mapping problems from quasigeodesic sources (Theorems \ref{thm:adjoint0}, \ref{thm:univ-morph}). The section contains a short detour on extremal metrics (Theorem \ref{thm:extremal}), in the spirit of Rosendal’s work.

  Section \ref{sec:universal} introduces uniformly controlled diagrams in coarse geometry and the associated ($\cat$‑geodesic) universal cones. We give equivalent criteria for existence of these universal objects, present the explicit Rips--Tuple construction, and define a uniformly controlled analogue of natural transformations; we also show how universal cones behave under controlled retractions of the input data.

  Section \ref{sec:hhs} applies the formalism to hierarchically hyperbolic spaces: we encode HHS data as a uniformly controlled diagram and prove that the total space realises the universal quasigeodesic cone over that data (Theorem \ref{thm:hhs-universal}), deducing further consequences such as the Rips--Tuple description (Theorem \ref{thm:hhs-image}). We conclude by treating pairwise constrained families of metric spaces and proving that uniformly compatible families of retractions induce canonical retractions of the corresponding total spaces, yielding Theorem \ref{thm:hier-retractions}.

  \subsection*{Acknowledgements}

  The author is grateful to Mark Hagen for numerous enlightening discussions regarding hierarchically hyperbolic spaces, and to Jing Tao for helpful feedback. Thanks to Valentina Disarlo, Harry Petyt, Alessandro Sisto, Davide Spriano,  Federico Vigolo, Richard Webb, and Abdul Zalloum for interesting conversations, and to Sanghyun Kim for pointing out connections to work of Rosendal.
  The author acknowledges support from the National Natural Science Foundation of China (NSFC 12101503); the Suzhou Science and Technology Development Planning Programme (ZXL2022473); and the XJTLU Research Development Fund (RDF-23-01-121).

  \section{Coarse geometry}\label{sec:background}

  In this section, we establish terminology and notation for working with coarse geometry. Background on this subject can found in \cite{Roe03, NY12, CH16, DK18, LV23}. We work with extended metric spaces, that is, where the metric is permitted to take the value infinity.

  Given a metric space $X$ and $\kappa \geq 0$, write $x \approx_\kappa x'$ to mean $d_X(x,x') \leq \kappa$ for $x,x' \in X$.

  \begin{definition}[Closeness]
  Let $X$ and $Y$ be metric spaces. Given $\kappa \geq 0$, say that two maps $f,f' \from X \to Y$ are \emph{$\kappa$--close}, denoted as $f \approx_\kappa f'$, if $fx \approx_\kappa fx'$ for all $x \in X$. Say that $f$ and $f'$ are \emph{close}, denoted as $f \approx f'$, if they are $\kappa$--close for some $\kappa \geq 0$.
  \end{definition}

  The notion of closeness is an equivalence relation on the set of all maps between a given pair of metric spaces. Let $\bar f$ denote the closeness class of $f$. In fact, closeness depends only on the metric on the codomain.

  \begin{definition}[Controlled map]
  Let $X$ and $Y$ be metric spaces and $f \from X \to Y$ be a map. A proper increasing function $\rho \from [0,\infty) \to [0,\infty)$ is an \emph{upper (resp. lower) control} for $f$ if
  \[d_Y(fx, fx') \leq \rho d_X(x,x') \qquad\left( \textrm{resp.}\quad d_Y(fx, fx') \geq \rho d_X(x,x') \;\right)\]
  for all $x,x'\in X$. We say that $f$ is \emph{controlled} if it admits an upper control.
  \end{definition}

  To handle extended metrics, we adopt a standard convention: if $f$ admits an upper control, then $d_Y(fx, fx') = \infty$ implies that $d_X(x,x') = \infty$; similarly, if $f$ admits a lower control then $d_X(x,x') = \infty$ implies that $d_Y(fx, fx') = \infty$.

   The following is immediate from the definitions and the triangle inequality.

   \begin{lemma}[Controls, closeness, and composition]\label{lem:control_close}
   Let $f \from X \to Y$ and $g \from Y \to Z$ be maps admitting upper controls $\rho$ and $\rho'$ respectively. Then $gf$ admits an upper control $\rho' \rho$. Moreover, if a map $f' \from X \to Y$ is $\kappa$--close to $f$ then $t \mapsto \rho(t) + 2\kappa$ is an upper control for $f'$. \qed
   \end{lemma}

    Motivated by this observation, we introduce the notion of control classes.  These shall be used to define some familiar categories appearing in coarse geometry.

  \begin{definition}[Control class]
   A \emph{control class} is a family $\F\subseteq [0,\infty)^{[0,\infty)}$ of functions satisfying the following:
   \begin{enumerate}
    \item $\F$ is closed under composition,
       \item if $\rho \in \F$ and $\rho' \in [0,\infty)^{[0,\infty)}$ satisfy $\rho'(t) \leq \rho(t)$ for all $t \geq 0$ then $\rho' \in \F$,
    \item any affine function $t \mapsto at + b$, where $a,b \geq 0$, belongs to $\F$.
   \end{enumerate}
  \end{definition}

  \begin{example}
   The following are examples of control classes.
   \begin{itemize}
    \item $\mathsf{Aff}$: functions bounded above by some affine function,
    \item $\mathsf{Poly}$: functions bounded above by some polynomial function, \item $\mathsf{All} = [0,\infty)^{[0,\infty)}$.
   \end{itemize}
   \end{example}

   Consequently, it is meaningful to define a category of metric spaces where the morphisms are those admiting an upper control from a given control class.

   \begin{definition}[Category associated to control class]
  Given a control class $\F$, let $\cat = \cat(\F)$ be the category whose objects are metric spaces, and whose morphisms are the maps which admit some upper control from $\F$. Let $\barcat = \barcat(\F)$ denote the quotient category of $\cat$ where the morphisms are considered up to closeness (and with the same objects as $\cat$).
  \end{definition}

  Write $\Born := \cat(\mathsf{All})$ and $\CLip := \cat(\mathsf{Aff})$. Their respective quotient categories are the \emph{metric coarse category} $\barBorn = \barcat(\mathsf{All})$ and the \emph{coarsely Lipschitz category} $\barCL = \barcat(\mathsf{Aff})$. Note that $\cat(\F)$ (resp. $\barcat(\F)$) is a wide subcategory of $\Born$ (resp. $\barBorn$) for any control class $\F$.

  Given a proper increasing function $\rho \from [0,\infty) \to [0,\infty)$, define $\rho^T \from [0,\infty) \to [0,\infty)$ by
  \[\rho^T(t) := \inf\{s\geq 0 ~|~ t \leq \rho(s)\}.\]
  This is also increasing and proper. Note that $\rho^T(t) \leq s \iff t \leq \rho(s)$ for all $s,t \geq 0$. Moreover, $\rho^T(\rho(t)) \leq t \leq \rho(\rho^T(t))$ for all $t \geq 0$. Define $\F^T := \{\rho^T ~|~ \rho \in \F\}$.

  The fact that $\CLip$ has restrictions on its controls, while $\Born$ does not, makes a crucual difference throughout this paper. Let us formalise this distinction as follows.

  \begin{definition}[Dominated control class]\label{def:dominated}
   A control class $\F$ is \emph{dominated} if there exists a function $\Theta \from [0,\infty) \to [1,\infty)$ that \emph{eventually always exceeds} every function in $\F$:    for all $\rho \in \F$, there exists some $t_0 \geq 0$ such that $\rho(t) < \Theta(t)$ for all $t \geq t_0$. We call $\Theta$ a \emph{dominating function} for $\F$. We say that $\cat(\F)$ or $\barcat(\F)$ has \emph{dominated controls} if $\F$ is a dominated control class.
  \end{definition}

   For example, $\Theta(t) = 2^t$ serves as a dominating function for both  $\mathsf{Aff}$ and $\mathsf{Poly}$. Any dominating function $\Theta$ must be proper since it eventually always exceeds $t \mapsto t$.

  Next, we present some useful characterisations of morphisms in $\barcat$.

  \begin{proposition}[Monomorphisms and epimorphisms in $\barcat$]\label{prop:mono-epi}
   Let $\cat$ be a category associated to some control class.
   Suppose that $f \from X \to Y$ is a morphism in $\cat$. Then
   \begin{itemize}
    \item $\bar f$ is a monomorphism in $\barcat$ if and only if $f$ admits a lower control, and
    \item $\bar f$ is an epimorphism in $\barcat$ if and only if $f$ is coarsely surjective: there exists $r \geq 0$ such that for all $y \in Y$, there exists some $x \in X$ satisfying $d_Y(y, fx) \leq r$.
   \end{itemize}
  \end{proposition}
  
  \proof
  The result holds for $\barBorn$ due to \cite[Proposition 3.A.16]{CH16} (see also \cite[Theorem 2.5]{Tang-mono} for a modified proof to handle extended metric spaces). Since left or right cancellativity persist when passing to any subcategory of $\barBorn$, the desired result also holds for $\barcat$.
  \endproof

  A map $f \from X \to Y$ is called a \emph{coarse equivalence} if it is coarsely surjective and admits both upper and lower controls. By \cite[Proposition 3.A.16]{CH16}, (closeness classes of) coarse equivalences are precisely the isomorphisms in $\barBorn$. Let us now characterise the isomorphisms in $\barcat$.

  \begin{proposition}[Isomorphisms in $\barcat$]\label{prop:iso-cat}
   Let $\cat = \cat(\F)$ be a category associated to a control class $\F$. Suppose that $f \from X \to Y$ is a morphism in $\cat$. Then $\bar f$ is an isomorphism in $\barcat$ if and only if $f$ is coarsely surjective and admits a lower control $\rho^T \in \F^T$.
  \end{proposition}

  \proof
  Assume that $\bar f \from X \to Y$ is an isomorphism in $\barcat$. In particular, $\bar f$ is an epimorphism in $\barcat$, hence $f$ is coarsely surjective by Proposition \ref{prop:mono-epi}. Choose any representative $g \from Y \to X$ for the inverse of $\bar f$ in $\barcat$. Then $g$ admits an upper control $\theta \in \F$ and there exists some $\kappa \geq 0$ such that $gf \approx_\kappa 1_X$.   By the triangle inequality, we deduce that
  \[d_X(x,x')  \leq d_X(gfx, gfx') + 2 \kappa\leq \theta d_Y(fx, fx') + 2\kappa = \rho(d_Y(fx, fx'))\]
  where $\rho(t) := \theta(t) + 2\kappa$. Note that $\rho \in \F$. Moreover,
  \[\rho^T\left(d_X(x,x')\right) \leq d_Y(fx, fx'),\]
  hence $\rho^T$ is a lower control for $f$.

  For the converse, assume that $f$ is coarsely surjective and admits a lower control $\rho^T\in\F^T$. Then $f$ is a coarse equivalence, hence $\bar f$ is an isomorphism in $\barBorn$ by \cite[Proposition 3.A.16]{CH16}. Let $\bar g \from Y \to X$ be the inverse of $\bar f$ in $\barBorn$. We claim that $\bar g$ is a morphism in $\barcat$, which would imply that $\bar f$ has an inverse in $\barcat$. Let $g$ be any representative of $\bar g$.
  Then $fg \approx_\kappa 1_Y$ for some $\kappa \geq 0$. Then
  \[\rho^T d_X(gy, gy') \leq d_Y(fgy, fgy') \leq d_Y(y,y') + 2\kappa,\]
  hence
  \[d_X(gy, gy') \leq \rho\rho^T d_X(gy, gy') \leq \rho \left(d_Y(y,y') + 2\kappa\right)\]
  for any $y,y'\in Y$ satisfying $d_Y(y,y') < \infty$. Consequently, $t \mapsto \rho(t + 2\kappa)$ is an upper control for $g$. Since $\rho$ and $t \mapsto t + 2\kappa$ belong to $\F$, and $\F$ is closed under composition, it follows that $g$ admits an upper control from $\F$. Thus, $\bar g$ is a morphism in $\barcat$.
  \endproof

   In particular, we recover the fact that isomorphisms in $\barCL$ are precisely the closesness classes of \emph{quasi-isometries}, that is, the coarse equivalences admitting both affine upper and lower controls. Furthermore, any coarse equivalence can be converted to a quasi-isometry upon replacing the metric on the domain with a coarsely equivalent one.

  \begin{lemma}[Coarse pullback metric]\label{lem:wlog_metric}
  Let $f \from (X,d_X) \to (Y,d_Y)$ be a coarse equivalence. Then there exists a metric $d'$ on $X$ such that the map $f' \from (X,d') \to (Y,d_Y)$ coinciding with $f$ on underlying sets is a quasi-isometry. Moreover, the underying identity on $X$ is a coarse equivalence between $(X,d_X)$ and $(X,d')$.
  \end{lemma}

  \proof
  Define a metric $d'$ on $X$ by setting by $d'(x,x') = \max\{1, \diam_Y(fx \cup fx')\}$ for $x,x' \in X$ distinct; and $d'(x,x') = 0$ when $x = x'$. Then the map $f' \from (X,d') \to (Y,d_Y)$ coinciding with $f$ on underlying sets is coarsely surjective since its image agrees with that of $f$. Moreover,
  \[d'(x,x') - 1 \leq d_Y(fx, fx') \leq d'(x,x')\]
  for all $x,x'\in X$. Hence $f'$ is a quasi-isometry. Furthermore, if $\rho$ is an upper (resp. lower) control for $f$, then $t \mapsto \rho(t) + 1$ (resp. $\rho)$ is an upper (resp. lower) control for the underlying identity $(X,d_X) \to (X,d')$.
  \endproof

  A subset $U \subseteq X$ is a \emph{$\cat$--retract} if there exists a morphism $r \from X \to U$ in $\cat$ such that $ri \approx 1_U$ where $i \from U \hookrightarrow X$ is the inclusion. Left- and right-invertibility in $\barcat$ characterise coarse analogues of sections and retractions with appropriate control.

  \begin{lemma}[Left-invertibility in $\barcat$]\label{lem:left-inverse}
   Let $\bar f \from X \to Y$ be a morphism in $\barcat$. Then $\bar f$ is left-invertible in $\barcat$ if and only if for any representative $f$, there exists a lower control $\rho^T \in \F^T$ for $f$ and $f(X) \subseteq Y$ is a $\cat$--retract.   \end{lemma}

  \proof
  Any representative $f$ admits a factorisation $f = if'$ via the corestriction $f' \from X \twoheadrightarrow f(X)$ and the inclusion $i \from f(X) \hookrightarrow Y$, where $f(X)$ is equipped with the induced metric from $Y$. Since $f'$ is surjective, $\bar f'$ is an isomorphism in $\barcat$ if and only if $f'$ (or equivalently $f$) admits a lower control $\rho^T \in \F^T$ by Proposition \ref{prop:iso-cat}.

  Assume there exists $\bar g \from Y \to X$ satisfying $\bar g \bar f = \bar 1_X$. Then, using the proof of Proposition \ref{prop:iso-cat}, $f$ admits a lower control $\rho^T \in \F^T$. Thus, we may write $\bar i = \bar f(\bar f')^{-1}$. Define a morphism $\bar r = \bar f'\bar g\from Y \to f(X)$ in $\barcat$. Then $\bar r \bar i = \bar f' \bar g \bar f(\bar f')^{-1} = \bar 1_{f(X)}$.

  For the converse, let $f \from X \to Y$ be a representative of $\bar f$ with lower control $\rho^T \in \F^T$.  Then $\bar f'$ is an isomorphism in $\barcat$. By assumption, there exists a morphism $r \from Y \to f(X)$ in $\cat$ satisfying $ri \approx 1_{f(X)}$. Let $\bar g = (\bar f')^{-1} \bar r \from Y \to X$. Then $\bar g \bar f = (\bar f')^{-1} \bar r \bar i \bar{f}' = \bar 1_X$ .
  \endproof

  \section{The metric Rips filtration and its colimit}\label{sec:Rips}

   \subsection{The Rips filtration and controlled maps}

   In this section, we investigate properties of the metric Rips filtration in the category $\barcat$. We will not consider any simplicial structure, rather, we focus primarily on its (large-scale) metric properties.
  \begin{definition}[Rips graph]
    The \emph{Rips graph} of a metric space $(X,d)$ at \emph{scale} $\sigma\geq 0$ is the graph $\Rips_\sigma (X,d)$ with $X$ as its vertex set, where distinct points $x,x' \in X$ are declared adjacent if and only if $d(x,x') \leq \sigma$. We shall write $\Rips_\sigma X$ for $\Rips_\sigma (X,d)$ if the metric is understood.
  \end{definition}

  It will be convenient to regard $\Rips_\sigma (X,d) = (X,d_\sigma)$ as a metric space with underlying set $X$, equipped with the standard graph metric $d_\sigma$. That is, $d_\sigma(x,x')$ is the infimum over integers $n\geq 0$ for which there exists a edge-path $x = x_0, \ldots, x_n = x'$ in $\Rips_\sigma (X,d)$.

  The graph $\Rips_\sigma X$ can be viewed as an approximation to $X$ which uses only knowledge of the metric $d$ up to scale $\sigma$. It may seem reasonable to expect that the approximations improve as the scale increases.  Indeed, given $\sigma \leq \tau$, whenever distinct points $x,x'\in X$ are adjacent in $\Rips_\sigma X$ then they are also adjacent in $\Rips_\tau X$,  hence the underlying identity map yields an inclusion  $\Rips_\sigma X \hookrightarrow \Rips_\tau X$. Rather than fixing a specific scale, we can simultaneously consider these approximations at all scales $\sigma \geq 0$ via the Rips filtration.

    \begin{definition}[Rips filtration]
    The \emph{(metric) Rips filtration of $X$} is the directed system
    \[\Rips_*X \from [0,\infty) \to \cat,\]
    where the morphisms $\Rips_\sigma X \to \Rips_\tau X$ are the (1--Lipschitz) identity maps between the underlying sets, defined whenever $\sigma\leq\tau$. We shall also write $\Rips_*X \from [0,\infty) \to \barcat$ for the associated directed system in the quotient category, where the morphisms are considered up to closeness.
    \end{definition}

    A standard fact is that the Rips filtration behaves nicely with respect to controlled maps, so long as we adjust the scale parameter for the target Rips filtration appropriately. Given a proper increasing function $\rho$, let $\Rips_{\rho *} X$ be the directed system whose morphisms are given by underlying identities $\Rips_{\rho \sigma} X \to \Rips_{\rho \tau} X$ for all $\sigma\leq\tau$. Note that when $\rho$ is the identity map on $[0,\infty)$ then $\Rips_{\rho *} X = \Rips_{*} X$.

    \begin{lemma}[Induced natural transformation]\label{lem:Rip_control}
    Suppose that $f \from X \to Y$ is a map with upper control $\rho$. Then $f$ induces a natural transformation $\Rips_* (f,\rho) \from \Rips_*X \tto \Rips_{\rho*}Y$ in $\cat$ whose components $\Rips_\sigma (f,\rho) \from \Rips_{\sigma} X \to \Rips_{\rho\sigma} Y$ coincide with $f$ on the underlying sets.
    \end{lemma}

    \proof
    Suppose $x,x' \in X$ are adjacent in $\Rips_{\sigma} X$. Then
    \[d_Y(fx, fx') \leq \rho(d_X(x,x')) \leq \rho(\sigma),\]
    hence $y$ and $y'$ are adjacent in $\Rips_{\rho\sigma} Y$ for all $y\in fx$, $y'\in fx'$. Naturality is immediate.
    \endproof

    Next, we consider the colimit of the Rips filtration in the category $\barcat$. Recall that the \emph{colimit} of a diagram (in any category) is the universal cocone under the given diagram. Given a metric space $Y$, a \emph{cocone under} $\bar\mu\from \Rips_{\rho*} X \tto Y$ with \emph{nadir} $Y$ in $\barcat$ is specified by a family of morphisms $\bar\mu_\sigma \from \Rips_{\rho \sigma} X \to Y$ in $\barcat$ (known as the colegs) for each $\sigma \geq 0$ such that whenever $\sigma \leq \tau$, we have that $\bar\mu_\sigma = \bar\mu_\tau\bar\iota$ where $\iota \from \Rips_{\rho \sigma} X \to \Rips_{\rho\tau} X$ is (the closeness class of) the inclusion. Equivalently, we require that each closeness class $\bar\mu_\sigma$ agrees with the closeness class of a common underlying function $m \from X \to Y$.

    \begin{definition}[Rips colimit of a metric space]
    Let $X$ be a metric space. The \emph{Rips colimit} of $X$ in $\barcat$ is defined as the colimit $\floorC{X} := \colimit \Rips_* X$ in the category $\barcat$.
    \end{definition}

    In other words, the limit cocone $\bar\lambda \from \Rips_* X \tto \floorC{X}$ is characterised by the property: for any cocone $\bar\mu \from \Rips_* X \tto Y$, there exists a unique morphism $\bar f \from \floorC{X} \to Y$ in $\barcat$ such that $\bar\mu = \bar f\bar\lambda$.    Note that while the colimit $\floorC{X}$ (if it exists) is not necessarily unique, by the universal property, it is unique up to unique isomorphism in $\barcat$. For example, if we work in the coarsely Lipschitz category, the Rips colimit $\floor{X}_{\CLip}$ is unique up to quasi-isometry. We shall henceforth write $\floorC{X}$ to denote a choice of Rips colimit for $X$.

    If we instead choose to work in the category $\cat$, then a cocone $\mu\from \Rips_{\rho*} X \tto Y$ must have all colegs $\mu_\sigma$ agreeing \emph{exactly} with a common underlying function $m \from X \to Y$. It turns out that this is a fairly mild constraint, as we may always arrange for this to hold for some choice of coleg representatives of a given cocone $\bar\mu\from \Rips_{\rho*} X \tto Y$.

    \begin{lemma}[Common underlying colegs]\label{lem:colegs}
     Let $\rho \from [0,\infty) \to [0,\infty)$ be an increasing function and $\bar\mu \from \Rips_{\rho*} X \tto Y$ be a cocone in $\barcat$. Then there is a function $m \from X \to Y$ such that each coleg $\bar\mu_{\sigma} \from \Rips_{\rho\sigma} X \to Y$ admits a representative coinciding with $m$ on underlying sets. In particular, $\bar\mu$ is realised by a cocone $\mu\from\Rips_{\rho *} X \tto Y$ in $\cat$.
    \end{lemma}

    \proof
    Let $m \from X \to Y$ be the underlying function of some representative of ${\bar\mu_0 \from \Rips_{\rho(0)} X \to Y}$. For $\sigma \geq 0$, write $\iota \from \Rips_{\rho(0)} X \to \Rips_{\rho\sigma} X$ for the 1--Lipschitz map coinciding with the underlying identity. Now $\bar\mu_0 = \bar\mu_\sigma \bar\iota$ since $\bar\mu$ is a cocone, hence $m$ is close to $\mu_\sigma$ as functions from $X$ to $Y$. Therefore, the map $\mu'_\sigma \from \Rips_\sigma X \to Y$ coinciding with $m$ on underlying sets is also close to $\mu_\sigma$. By assumption, ${\mu_\sigma \from \Rips_\sigma X \to Y}$ is a morphism in $\cat$ and so, by Lemma \ref{lem:control_close}, the same is true for $\mu'_\sigma$. Since this holds for all $\sigma \geq 0$, the result follows.
    \endproof

    Suppose that $X$ and $Y$ are metric spaces whose Rips colimits in $\barcat$ exist.  Given a (closeness class of a) controlled map $\bar f \from X \to Y$, we wish to define a morphism $\floorC{\bar f} \from \floorC{X} \to \floorC{Y}$ in $\barcat$. As a first step, choose a representative $f \in \bar f$ and corresponding upper control $\rho$. Taking the colimit of the natural transformation $\Rips_* (f,\rho) \from \Rips_*X \tto \Rips_{\rho*}Y$ in $\barcat$ yields a morphism
     \[\floorC{(f,\rho)} := \colimit \Rips_* (f,\rho) \from \floorC{X} \to \colimit \Rips_{\rho*}Y\]
    in $\barcat$. Our strategy is to verify that the right-hand colimit is canonically isomorphic to $\floorC{Y}$, and that the morphism does not depend on the choice of representative map or upper control.

    Let us examine how adjusting the scale parameter affects cocones under the Rips filtration, working in $\barcat$. Given a proper increasing function $\rho$ and cocone $\bar\mu \from \Rips_* X \tto Y$, define a cocone $F_Y\bar\mu \from \Rips_{\rho *} X \tto Y$ whose colegs are $(F_Y\bar\mu)_\sigma :=\bar\mu_{\rho\sigma}$.

    \begin{lemma}[Reparameterising scale]\label{lem:reparam}
     Let $X$ be a metric space and $\rho$ be a proper increasing function. Then there is a natural isomorphism
     \[F \from \Cone(\Rips_* X, -) \tto \Cone(\Rips_{\rho *} X, -)\]
     whose components $F_Y \from \Cone(\Rips_* X, Y) \to \Cone(\Rips_{\rho *} X, Y)$ are given by $\bar\mu \mapsto F_Y\bar\mu$.
    \end{lemma}
    \proof
    First, let us verify naturality of $F$. Let $\bar f \from Y \to Z$ be a morphism in $\barcat$. Then for any cocone $\bar\mu \from \Rips_* X \tto Y$, we have that
    \[(F_Z\bar f \bar\mu)_\sigma = (\bar f \bar \mu)_{\rho \sigma} = \bar f (\bar\mu_{\rho\sigma}) = \bar f  (F_Y \bar\mu)_\sigma = (\bar f F_Y \bar\mu)_\sigma\]
    for all $\sigma \geq 0$, hence $F_Z\bar f \bar\mu = \bar f F_Y \bar\mu$.

    Next, we verify injectivity of $F_Y$. Let $\bar\mu, \bar\mu' \from \Rips_* X \tto Y$ be cocones satisfying $F_Y\bar\mu = F_Y\bar\mu'$. Then
    \[\bar\mu_\sigma = (F_Y\bar\mu)_\sigma = (F_Y\bar\mu')_\sigma = \bar\mu'_\sigma \from \Rips_{\rho\sigma}X \to Y\] for all $\sigma \geq 0$, hence $\bar\mu = \bar\mu'$ as desired.

    For surjectivity, let $\bar\nu\from \Rips_{\rho*} X \tto Y$ be a cocone. By Lemma \ref{lem:colegs}, we may assume all representative colegs $\nu_\tau \from \Rips_{\rho(\tau)} X \to Y$ have the same underlying function $n \from X \to Y$. Let $\mu_\sigma \from \Rips_\sigma X \to Y$ be the map coinciding with $n$ on underlying sets. Since $\rho$ is proper, for each $\sigma \geq 0$ there exists $\tau\geq 0$ such that $\sigma \leq \rho(\tau)$. Therefore, $\mu_\sigma = \nu_\tau \iota$ where ${\iota \from \Rips_\sigma X \to \Rips_{\rho(\tau)} X}$ is the 1--Lipschitz inclusion map. Consequently, each $\bar\mu_\sigma$ is a morphism in $\barcat$, and so they assemble to form a cocone $\bar\mu \from \Rips_* X \tto Y$ in $\barcat$. By construction, this satisfies $F_Y\bar\mu = \bar\nu$.
    \endproof

    Thus, $\bar\mu\from\Rips_* X \tto Y$ is a limit cocone if and only if the same is true for $F_Y\bar\mu \from\Rips_{\rho *} X \tto Y$. In particular, a reparameterision of the Rips filtration does not affect its colimit.

    \begin{corollary}[Reparameterised Rips colimit]\label{cor:adjust-colimit}
     Let $X$ be a metric space and $\rho$ be a proper increasing function. Then there is a canonical isomorphism $\floorC{X} := \colimit \Rips_{*} X \cong \colimit \Rips_{\rho *} X$ in $\barcat$, where each colimit exists if and only if the other does. \qed
    \end{corollary}

    Consequently, we may choose a common colimit $\floorC{X}$  to realise $\colimit \Rips_{\rho *} X$ for any reparameterisation $\rho$ (so long as the colimit exists). Thus, we may define a morphism
    \[\floorC{f} := \floorC{(f,\rho)} \from \floorC{X} \to \floorC{Y}\]
    in $\barcat$ where $\rho$ is any upper control for $f\from X \to Y$. To verify that this depends only on the closeness class of $f$, we consider the behaviour of $\Rips_* (f,\rho)$ under change of representative map.

    \begin{lemma}[Adjusting representative]\label{lem:adjust-rep}
    Let $f,f' \from X \to Y$ be $\kappa$--close maps, and suppose $f$ has an upper control $\rho \in \F$. Let $\rho'(t) = \rho(t) + 2\kappa$. Then the natural transformations $\Rips_* (f,\rho')$ and $\Rips_* (f',\rho')$ coincide in $\barcat$.
    \end{lemma}

    \proof
    By Lemma \ref{lem:control_close}, $f'$ admits $\rho'$ as an upper control, and so the natural transformation $\Rips_* (f',\rho')$ is defined. Let $\sigma \geq 0$. Since $f \approx_\kappa f'$, we deduce that $d_Y(fx, f'x) \leq \kappa \leq \rho'(\sigma)$ for all $x \in X$. Therefore, $fx$ and $fx'$ are adjacent in $\Rips_{\rho'\sigma} Y$ for all $x\in X$, hence $\Rips_\sigma (f,\rho')$ and $\Rips_\sigma (f',\rho')$ are $1$--close.    Since this holds for every $\sigma \geq 0$, it follows that $\Rips_* (f,\rho')$ and $\Rips_* (f',\rho')$ coincide in $\barcat$.
    \endproof

    \begin{lemma}[Well-definedness]\label{lem:well-defined}
     Let $X$ and $Y$ be metric spaces such that $\floorC{X}$ and $\floorC{Y}$ exist. Let $\bar f \from X \to Y$ be (the closeness class of) a controlled map. Then for any representatives $f,f'$ for $\bar f$, we have that $\floorC{f} = \floorC{f'} \from \floorC{X} \to \floorC{Y}$ in $\barcat$.
    \end{lemma}

    \proof
    Let $\rho, \rho'$ be upper controls for $f,f'$ respectively. By Lemma \ref{lem:adjust-rep}, there exists an upper control $\rho''$ for $f,f'$ such that $\Rips_* (f,\rho'') = \Rips_* (f',\rho'')$ in $\barcat$. Therefore
     \[\floorC{f} = \floorC{(f,\rho'')} = \floorC{(f', \rho'')} = \floorC{f'} \]
    in $\barcat$ as desired.
    \endproof

     We are now ready to define the action of the Rips colimit on morphisms in $\barBorn$.

    \begin{definition}[Rips colimit of a closeness class]\label{def:Rips_f}
     Let $\bar f \from X \to Y$ be the closeness class of a controlled map. Then its Rips colimit in $\barcat$ is defined as
     \[\floorC{\bar f} := \floorC{f} \from \floorC{X} \to \floorC{Y} \]
     for any choice of representative $f \in \bar f$.
    \end{definition}

    To verify functoriality, we consider how the Rips colimit behaves with respect to composition.

    \begin{lemma}[Rips colimit and composition]
    Let $\bar f \from X \to Y$ and $\bar g \from Y \to Z$ be morphisms in $\barBorn$. Then $\floorC{\bar g} \floorC{\bar f} = \floorC{\bar g \bar f}$, so long as the involved Rips colimits in $\barcat$ exist.
    \end{lemma}

    \proof
    Choose representatives $f\in \bar f$ and $g\in \bar g$, with respective (proper) upper controls $\rho, \rho'$. Then
    \[\Rips_* (g,\rho') \circ \Rips_* (f,\rho) = \Rips_* (gf,\rho'\rho) \from \Rips_* X \tto \Rips_{\rho'\rho *} Z\]
    since all components coincide with $gf$ on the underlying sets. The result follows by taking colimits in $\barcat$ and using the fact that colimits behave functorially whenever they exist (see \cite[Proposition 3.6.1]{Rie16}).
    \endproof

    The preceding results are summarised as follows. Let $\Rc(\barcat)$ denote the full subcategory of $\barBorn$ whose objects are metric spaces which admit a Rips colimit in $\barcat$.

    \begin{proposition}[Functoriality of Rips colimit]\label{prop:functorial}
     A choice of Rips colimit $\floorC{X}$ in $\barcat$ for each object $X \in \Rc(\cat)$ defines the action on objects of a functor $\floorC{-} \from \Rc(\barcat) \to \barcat$. The action on morphisms is given by Definition \ref{def:Rips_f}.\qed
    \end{proposition}

     \begin{corollary}[Rips colimit of coarse equivalence]
     Let $X$ and $Y$ be metric spaces whose Rips colimits in $\barcat$ exist.     If $f \from X \to Y$ is a coarse equivalence then $\floorC{\overline{f}} \from \floorC{X} \to \floorC{Y}$ is an isomorphism in $\barcat$. \qed
    \end{corollary}

  \subsection{Rips colimit in the metric coarse category}\label{coarse-colimit}

  A natural question is whether the Rips colimit $\floorC{X}$ of a given metric space $X$ exists. As we shall see, the answer is dependent on the choice of ambient category $\barcat$. We focus here on the metric coarse category, and defer the general case to Section \ref{sec:colimit-geodesic}.

  For each $\sigma \geq 0$, let $\xi_\sigma \from \Rips_\sigma X \to X$ be the map coinciding with the underlying identity. Observe that $\xi_\sigma$ is $\sigma$--Lipschitz for $\sigma > 0$; while $\xi_0$ is $1$--Lipschitz since $\Rips_0 X$ is the empty graph with vertex set $X$. Therefore, each $\xi_\sigma$ is a morphism in $\cat$.

    \begin{definition}[Canonical cocone]
     Let $X$ be a metric space. The \emph{canonical cocone} under $\Rips_* X$ is the cocone $\xi \from \Rips_* X \tto X$ in $\cat$ whose colegs $\xi_\sigma \from \Rips_\sigma X \to X$ coincide with the underlying identity on $X$. Let $\bar\xi \from \Rips_* X \tto X$ denote the corresponding cocone in $\barcat$ with colegs $\bar\xi_\sigma$.
    \end{definition}

    Appealing to the universal property, the limit cocone under $\Rips_* X$ (if it exists) and the canonical cocone are related by a canonical morphism in $\barcat$.

    \begin{definition}[Rips counit at $X$]
    Let $X$ be a metric space such that $\floorC{X}$ exists. The \emph{Rips counit at $X$ in $\barcat$} is
    \[\bar\epsilon_X := \colimit \bar\xi \from \floorC{X} \to X.\]
    That is, $\bar\epsilon_X$ is the unique morphism in $\barcat$ satisfying $\bar\xi = \bar\epsilon_X \bar\lambda$, where $\bar\lambda \from \Rips_* X \tto \floorC{X}$ is the limit cocone. We shall write $\epsilon_X \from \floorC{X} \to X$ for any representative of $\bar\epsilon_X$.
     \end{definition}

   We now give a proof of Theorem \ref{thm:CE-Rips}.

    \begin{proposition}[Rips colimit in $\barBorn$]\label{prop:colimit_coarse}
      Let $X$ be a metric space. Then, the canonical cocone $\bar\xi\from \Rips_* X \tto X$ is the limit cocone in $\barBorn$. In particular, $\bar\epsilon_X \from \floor{X}_{\Born} \to X$ is an isomorphism.
    \end{proposition}

    \proof
    Let $\overline{\mu} \from \Rips_* X \tto Y$ be a cocone. By Lemma \ref{lem:colegs}, we may choose a representative map $\mu_\sigma$ for each coleg $\bar\mu_\sigma$ to agree with a common underlying function $f \from X \to Y$.    We claim that $f$ is controlled and that $\bar f \from X \to Y$ is the unique morphism in $\barBorn$ satisfying $\bar\mu = \bar f \bar\xi$. Let $\rho_\sigma$ be an upper control for $\mu_\sigma$ for each $\sigma\geq 0$. Suppose that $x, x' \in X$ satisfy $d_X(x,x') \leq \sigma$. Then $x,x'$ are adjacent in $\Rips_\sigma X$, hence
    \[d_Y(fx, fx') = d_Y(\mu_\sigma x , \mu_\sigma x') \leq \rho_\sigma(1).\]
    Consequently, $\rho(t) := \inf_{\sigma \geq t} \rho_\sigma(1)$ serves as an upper control for $f$. By considering the functions on underlying sets, we deduce that $f\xi_\sigma = \mu_\sigma$ for all $\sigma \geq 0$, hence $\bar\mu = \bar f \bar \xi$. Finally, suppose $\bar h \from X \to Y$ is a morphism in $\barBorn$ satisfying $\bar\mu = \bar h \bar\xi$. Then $\bar h \bar\xi_0 = \bar\mu_0 = \bar f \bar\xi_0$. Since $\xi_0$ is surjective, $\bar\xi_0$ is an epimorphism in $\barBorn$ by Proposition \ref{prop:mono-epi}, and thus $\bar h = \bar f$.
    \endproof

    \begin{corollary}[Sequential cocompleteness] Let $X$ be a metric space. Then $X$ is a (sequential) colimit of graphs in $\barBorn$. \qed
    \end{corollary}

  \subsection{\texorpdfstring{$\cat$}{$\mathsf{C}$}--geodesic spaces}\label{sec:geodesic}

  A standard fact is that any quasigeodesic space $X$ is quasi‑isometric to a geodesic space -- concretely, to its Rips graph at sufficiently large scale.     Replacing quasi‑isometry by isomorphism in $\barcat$ yields analogous classes of spaces. For example, the coarsely geodesic metric spaces are exactly those that are coarsely equivalent to a genuine geodesic space.

  \begin{definition}[$\cat$--geodesic space]
   A metric space $X$ is \emph{$\cat$--geodesic} if it is isomorphic in $\barcat$ to a geodesic space or, equivalently, to a graph with the combinatorial metric. Let $\Gr(\cat)$ (resp. $\Gr(\barcat)$) denote the full subcategory of $\cat$ (resp. $\barcat$) whose objects are $\cat$--geodesic spaces.
  \end{definition}

  Thus the $\CLip$--geodesic (resp. $\Born$--geodesic) spaces are precisely the quasigeodesic (resp. coarsely geodesic) spaces. Let us also generalise the notion of a quasigeodesic path.

   \begin{definition}[$\cat$--geodesic path]
   Let $\sigma\geq 0$ and $\rho$ be a proper increasing function. Let $X$ be a metric space and $x,x'\in X$. A \emph{$(\sigma,\rho)$--path} from $x$ to $x'$ is a sequence $x = x_0, x_1, \ldots, x_n = x'$ in $X$ such that
   \[\rho^T(|k-j|) \leq d(x_j, x_k) \leq \sigma|k-j| \]
   for all integers $j,k\in [0,n]$. If $\rho\in\F$ then we also call such a sequence a \emph{$\cat$--geodesic path}.
   \end{definition}

  If $\rho$ is affine, this recover the usual notion of a quasigeodesic path.  Using a standard argument, we show that the $\cat$--geodesic spaces are precisely those which admit uniform $\cat$--geodesic paths between evey pair of points (which are at finite distance).

  \begin{proposition}[Characterisations of $\cat$--geodesic spaces]\label{prop:C-geodesic}
   Let $X = (X,d)$ be a metric space. Then the following are equivalent.
   \begin{enumerate}
   \item $X$ is $\cat$--geodesic,
   \item There exists $\sigma \geq 0$ such that $\bar\xi_\sigma \from \Rips_\sigma X \to X$ is an isomorphism in $\barcat$,
   \item There exists $\sigma\geq 0$ such that $\xi_\sigma \from \Rips_\sigma X \to X$ admits a lower control $\rho^T\in\F^T$, and    \item There exists $\sigma\geq 0$ and a proper increasing function $\rho\in\F$    such that for all $x,x'\in X$ where $d(x,x') < \infty$, there exists a $(\sigma,\rho)$--path in $X$ from $x$ to $x'$.
   \end{enumerate}
   Moreover, any choice of $\sigma$ and $\rho$ satisfying either item (3) or (4) will also satisfy the other.
  \end{proposition}

  \proof
  
  (1) $\implies$ (2). It suffices to show that the inverse function $\xi_\sigma^{-1} \from X \to \Rips_\sigma X$ admits an upper control from $\F$ for $\sigma \geq 0$ sufficiently large. Let $\Gamma$ be a graph and $\bar f \from \Gamma \to X$ be an isomorphism in $\barcat$. Choose a representative $f$. Since $f$ is controlled, it must be $\sigma$--Lipschitz for $\sigma \geq 0$ sufficiently large.  Therefore, $\xi_\sigma \from \Rips_\sigma X \to X$ factors through the 1--Lipschitz map $h \from \Gamma \to \Rips_\sigma X$ which coincides with $f$ on underlying sets. Let $g \from X \to \Gamma$ be a representative of the inverse of $\bar f$ in $\barcat$. Then there exists some $\kappa \geq 0$ such that $fg \approx_\kappa 1_X$. Let us also assume that $\sigma \geq \kappa$. Then for any $x\in X$, we have that $d(x, fgx) \leq \kappa \leq \sigma$, hence $d_\sigma(\xi_\sigma^{-1} x, hgx) \leq 1$. Therefore, the function $\xi_\sigma^{-1}$ is 1--close to $hg$. Since $g, h$ admit upper controls from $\F$, so does $\xi_\sigma^{-1}$. Therefore, $\bar\xi_\sigma$ is an isomorphism in $\barcat$ for $\sigma\geq0$ sufficiently large.

  (2) $\implies$ (1). Trivial.

  (2) $\iff$ (3). This follows from Proposition \ref{prop:iso-cat}.

  (3) $\implies$ (4).
  Suppose that $x,x' \in X$ satisfy $d(x,x') < \infty$.   Then $d_\sigma(x,x') < \infty$ since $\xi_\sigma$ is a coarse equivalence, and so there exists a path $x = x_0, x_1, \ldots, x_n = x'$ from $x$ to $x'$ in $\Rips_\sigma X$, where $n = d_\sigma(x,x')$. Therefore,
  \[\rho^T(|k-j|) = \rho^T(d_\sigma(x_j, x_k)) \leq d(x_j, x_k) \leq \sigma d_\sigma(x_j, x_k) = \sigma |k-j| \]
  for all integers $j,k \in [0,n]$.

  (4) $\implies$ (3). We claim that $\rho^T$ is a lower control for $\xi_\sigma \from \Rips_\sigma X \to X$. Suppose that $x,x'\in X$ satisfy $d_\sigma(x,x') < \infty$. Then $d(x,x') \leq \sigma d_\sigma(x,x') < \infty$ since $\xi_\sigma$ is $\sigma$--Lipschitz . Therefore, by assumption, there is a sequence $x = x_0, x_1, \ldots, x_n = x'$ in $X$ such that $\rho^T(n) \leq d(x,x')$ and $d(x_i, x_{i+1}) \leq \sigma$ for all $i$. This yields a path in $\Rips_\sigma X$ from $x$ to $x'$, hence $d_\sigma(x,x') \leq n$. Consequently, $\rho^T(d_\sigma(x,x')) \leq \rho^T(n) \leq d(x,x')$ as required.
  \endproof

  \begin{remark}
  If $X$ is a geodesic space then we may take $\sigma = 1$ and $\rho(t) = t+1$.
  \end{remark}

  We conclude this section with two useful facts. The first is a characterisation of coarsely geodesic spaces in terms stability of the Rips filtration. This is a standard result; a proof in the more general setting of coarse structures is given by Roe \cite{Roe03}. We say that the Rips filtration $\Rips_* X$ \emph{stabilises} if the underlying identities $\Rips_\sigma X \to \Rips_\tau X$ represent isomorphisms in $\barcat$ for all $\sigma \leq \tau$ large. Since each $\Rips_\sigma X$ is a graph, the underlying identity $\Rips_\sigma X \to \Rips_\tau X$ is a quasi-isometry if and only if it is a coarse equivalence. Therefore, the stability of $\Rips_* X$ does not depend on the choice of ambient category $\barcat$.

  \begin{lemma}[Stable Rips colimit]\label{lem:stable}
   Let $X$ be a metric space. Then $\Rips_* X$ stabilises if and only if $X$ is coarsely geodesic.
  \end{lemma}

  \proof
  Assume $X$ is coarsely geodesic. By Proposition \ref{prop:C-geodesic}, there exists some $\sigma \geq 0$ such that $\bar\xi_\sigma \from \Rips_\sigma X \to X$ is an isomorphism in $\barBorn$. Moreover, this also holds for $\bar\xi_\tau$ whenever $\sigma \leq \tau$. Therefore, the (closeness class of the) underlying identity $\bar\xi_\tau^{-1}\bar\xi_\sigma \from \Rips_\sigma X \to \Rips_\tau X$ is an isomorphism for $\sigma \leq \tau$.

  For the converse, suppose that $\Rips_* X$ stabilises. Appealing to Proposition \ref{prop:colimit_coarse}, we deduce that $X \cong \floor{X}_{\Born} \cong \Rips_\sigma X$ in $\barBorn$ for $\sigma$ large. Therefore $X$ is coarsely equivalent to a graph and is hence coarsely geodesic.
  \endproof

  The second is that the property of being $\cat$--geodesic is preserved under retractions in $\barcat$.

   \begin{lemma}[Retractions and $\cat$--geodesicity]\label{lem:cgeod-ret}
   Let $X$ be a $\cat$--geodesic space. Let $Y$ be a metric space and $\bar r \from X \to Y$ be a right-invertible morphism in $\barcat$. Then $Y$ is $\cat$--geodesic.
  \end{lemma}

  \proof
  Let $\bar s \from Y \to X$ be a right inverse of $\bar r$ in $\barcat$.
  Choose representatives $r,s$. Then there exists some $\kappa \geq 0$ such that $rs \approx_\kappa 1_Y$. Let $\rho \in \F$ be an upper control for $r$.  Let $\xi \from \Rips_*X \tto X$ and $\zeta \from \Rips_*Y \tto Y$ denote the respective canonical cocones. Since $X$ is $\cat$--geodesic, there exists some $\sigma \geq 0$ such that $\bar\xi_\sigma \from \Rips_\sigma X \to X$ is an isomorphism in $\barcat$.   Choose $\tau \geq \max(\rho\sigma, \kappa)$. We claim that $\bar\zeta_\tau \from \Rips_\tau Y \to Y$ is an isomorphism in $\barcat$. Consider the maps
  \[Y \xrightarrow{s} X \xrightarrow{\xi_{\sigma}^{-1}} \Rips_{\sigma} X \xrightarrow{r'} \Rips_\tau Y \xrightarrow{\zeta_\tau} Y, \]
  where $r'$ is the 1--Lipschitz map coinciding with $r$ on underlying sets. Each arrow is a morphism in $\cat$, and the composite coincides with $rs \approx 1_Y$ on underlying sets. Therefore, $\bar\zeta_\tau$ admits a right inverse $\bar r' \bar\xi_\sigma^{-1} \bar s$ in $\barcat$. For the other direction, consider
  \[\Rips_\tau Y \xrightarrow{\zeta_\tau} Y \xrightarrow{s} X \xrightarrow{\xi_{\sigma}^{-1}} \Rips_{\sigma} X \xrightarrow{r'} \Rips_\tau Y.\]
  The composite again coincides with $rs$ on underlying sets; the assumption that $\tau \geq \kappa$ ensures that it is 1--close to the identity on $\Rips_\tau Y$. Consequently, $\bar\zeta_\tau$ admits a left inverse $\bar r' \bar\xi_\sigma^{-1} \bar s$ in $\barcat$, and is hence an isomorphism in $\barcat$. Therefore, $Y$ is $\cat$--geodesic.
  \endproof

  \subsection{Colimits of \texorpdfstring{$\cat$}{$\mathsf{C}$}--geodesic spaces}\label{sec:colimit-geodesic}

  We now characterise metric spaces arising as colimits of $\cat$--geodesic spaces in $\barcat$.

  Let $\sJ$ be a small category, called the \emph{indexing category}.
  A \emph{diagram} (of shape $\sJ$) in $\barcat$ is a functor $D \from \sJ \to \barcat$. Let us write $D_j$ (resp.~$D_\phi$) to denote the object (resp.~morphism) of $D$ at $j \in \Obj(\sJ)$ (resp.~$\phi \in \Mor(\sJ)$).  We wish to determine when a given metric space $X$ realises the colimit $\colimit_{\sJ} D$ for some diagram $D$ whose objects are $\cat$--geodesic.
  The main technical step is to reduce such colimits to some Rips colimit.

  \begin{theorem}[Colimits of $\cat$--geodesic spaces]\label{thm:colimits}
   Let $X$ be a metric space. Then the following are equivalent:
   \begin{enumerate}
    \item There exists a limit cocone $\bar\lambda \from D \tto X$ in $\barcat$ for some diagram $D \from \sJ \to \barcat$ whose objects $D_j$ are $\cat$--geodesic spaces,
    \item There exists a limit cocone $\bar\lambda \from \Gamma \tto X$ in $\barcat$ for some diagram $\Gamma \from \sJ \to \barcat$ whose objects $\Gamma_j$ are graphs (equipped with the combinatorial metric),
    \item Every controlled map $f \from X \to Y$ to a metric space $Y$ is a morphism in $\cat$,
    \item $\bar\xi \from \Rips_*X \tto X$ is the limit cocone in $\barcat$, and
    \item $\floorC{X}$ exists and $\bar\epsilon_X \from \floorC{X} \to X$ is an isomorphism in $\barcat$.
   \end{enumerate}
\end{theorem}

  \proof

  (1) $\implies$ (2). Since $D_j$ is a $\cat$--geodesic space for each $j \in\Obj(\sJ)$, there exists an isomorphism $\bar\alpha_j \from D_j \to \Gamma_j$ in $\barcat$ to some graph $\Gamma_j$. For each arrow $\phi \from i \to j$ in $\sJ$, the morphism $\Gamma_\phi := \bar\alpha_j D_\phi\bar\alpha_i^{-1} \from \Gamma_i \to \Gamma_j$ in $\barcat$ satisfies $\Gamma_\phi \bar\alpha_i = \bar\alpha_j D_\phi $. Moreover, if $\phi,\phi' \in \Mor(\sJ)$ are arrows for which $\phi\phi'$ is defined, then $\Gamma_{\phi\phi'} = \Gamma_\phi\Gamma_{\phi'}$. Therefore, we obtain a diagram $\Gamma \from \sJ \to \barcat$ in $\barcat$ whose objects and morphisms are respectively $\Gamma_j$ for $j \in \Obj(\sJ)$ and $\Gamma_\phi$ for $\phi\in\Mor(\sJ)$. Furthermore, the $\bar\alpha_j$ form the components of a natural isomorphism $\bar\alpha \from D \tto \Gamma$ in $\barcat$. Consequently, if $\bar\lambda \from D \tto X$ is a limit cocone in $\barcat$ then the same is also true of $\bar\lambda\bar\alpha^{-1} \from \Gamma \tto X$.

  (2) $\implies$ (3).
  Assume $f \from X \to Y$ has upper control $\rho$. Consider the metric space $X' = (X,d')$ equipped with the metric $d'(x,x') := d_X(x,x') + d_Y(fx, fx')$ for all $x,x' \in X$. Then the map $g \from X' \to X$ coinciding with the underlying identity, and the map $f':= fg \from X' \to Y$ coinciding with $f$ on underlying sets are both $1$--Lipschitz. Therefore, $\bar g, \bar f'$ are morphisms in $\barcat$. Moreover, the inverse function $h = g^{-1} \from X \to X'$ has upper control $t \mapsto t + \rho(t)$, hence $\bar h = \bar g^{-1}$ is a morphism in $\barBorn$. We claim that $\bar h$ is a morphism in $\barcat$. It would then follow that $\bar f = \bar f' \bar g^{-1} = \bar f'\bar h$ is a morphism in $\barcat$.

  \begin{center}

    \tikzcdset{arrow style=math font}
    \begin{tikzcd}[column sep=huge, row sep=huge]
    \Gamma \arrow[Rightarrow, r, "\bar\lambda'"] \arrow[Rightarrow, d, "\bar\lambda"'] & X' \arrow[rightarrow, d, "\bar f'"] \arrow[leftarrow, ld, "\bar h"', bend right=10] \arrow[rightarrow, ld, "\bar g", bend left=10] \\
    X \arrow[r, "\bar f"'] & Y
    \end{tikzcd}
    \end{center}

  By assumption, there is a limit cocone $\bar\lambda \from \Gamma \tto X$ in $\barcat$ where $\Gamma \from \sJ \to \barcat$ is a diagram whose objects are graphs. Consider the cocone $\bar\lambda' := \bar h\bar\lambda \from \Gamma \tto X'$ in $\barBorn$. For each $j \in \sJ$, any coleg representative $h\lambda_j \from \Gamma_j \to X'$ is a controlled map from a graph and is hence coarsely Lipschitz. Therefore, $\bar\lambda'$ is also a cocone in $\barcat$. Since $\bar\lambda \from \Gamma \tto X$ is the limit cocone in $\barcat$, there is a unique morphism $\bar h' \from X \to X'$ in $\barcat$ such that $\bar\lambda' = \bar h'\bar\lambda$. Thus
  \[\bar g \bar h ' \bar\lambda = \bar g \bar\lambda' = \bar g \bar h \bar\lambda = \bar g \bar g^{-1} \bar\lambda = \bar\lambda \]
  holds in $\barBorn$. Since $\bar\lambda$ is a limit cocone in $\barcat$ and $\bar g \bar h'$ is a morphism in $\barcat$, it follows that $\bar g \bar h' = \bar 1_X$. Therefore, $\bar h = \bar h \bar g \bar h ' = \bar h '$ in $\barBorn$, hence $\bar h = \bar h'$ is a morphism $\barcat$ as claimed.

  (3) $\implies$ (4). Let $\bar\mu \from \Rips_*X \tto Y$ be a cocone in $\barcat$. By regarding $\bar\mu$ as a cocone in $\barBorn$ and applying Proposition \ref{prop:colimit_coarse}, there exists a unique morphism $\bar f \from X \to Y$ in $\barBorn$ such that $\bar\mu = \bar f \bar\xi$. By assumption, $\bar f$ is a morphism in $\barcat$. Therefore $\bar\xi$ is the limit cocone in $\barcat$.

  (4) $\iff$ (5). This follows from the definition $\bar\epsilon_X = \colimit \bar\xi$.

  (4) $\implies$ (1). This holds since each $\Rips_\sigma X$ is a graph, and hence $\cat$--geodesic.
  \endproof

  We obtain several immediate consequences of this result. Since $\floorC{X}$ is itself a colimit of graphs in $\barcat$, we may apply Theorem \ref{thm:colimits} to it to deduce idempotence of the Rips colimit.

  \begin{corollary}[Idempotence]\label{cor:idempotence}
   Assume that $X$ is a metric space such that $\floorC{X}$ exists. Then $\floorC{\floorC{X}}$ exists and $\bar\epsilon_{\floorC{X}} \from \floorC{\floorC{X}} \to \floorC{X}$ is an isomorphism in $\barcat$. \qed
  \end{corollary}

   Another consequence is that colimits of graphs persists when passing to a larger control class.

  \begin{lemma}[Enlargening control class]\label{lem:enlarge}
   Assume that $\barcat$ and $\barcat'$ are categories respectively associated to control classes $\F$ and $\F'$ satisfying $\F \subseteq \F'$.  Suppose that $\bar \lambda \from \Gamma \tto X$ is a limit cocone in $\barcat$, where $\Gamma \from \sJ \to \barcat$ is a diagram whose objects are graphs.   Then $\bar\lambda$ is also a limit cocone in $\barcat'$.
  \end{lemma}

  \proof
  Let $\bar\mu \from \Gamma \tto Y$ be a cocone in $\barcat'$. Since each object $\Gamma_j$ is a graph, all colegs of $\bar\mu$ are coarsely Lipschitz, hence $\bar\mu$ is a cocone in $\barcat$. Thus, $\bar\mu$ factors through a unique morphism $\bar f \from \Gamma \to Y$ in $\barcat$. Since $X$ is a colimit of graphs in $\barcat$, we may apply Theorem \ref{thm:colimits}(3) to deduce that any morphism in $\barcat'$ with domain $X$ must be a morphism in $\barcat$. Therefore, $\bar f$ is the unique morphism in $\barcat'$ through which $\bar\mu$ factors, and so $\bar\lambda$ is also the limit cocone in $\barcat'$.
  \endproof

  In the case where $\barcat' = \barBorn$, we may use the fact the canonical cocone is the limit cocone in $\barBorn$ to obtain desirable properties regarding the Rips counit $\bar\epsilon_X$ in $\barcat$.

  \begin{lemma}[Rips counit is a coarse equivalence]\label{lem:counit_equiv}
     Let $X$ be a metric space whose Rips colimit in $\barcat$ exists. Then the Rips counit $\bar\epsilon_X \from \floorC{X} \to X$ in $\barcat$ is an isomorphism in $\barBorn$. In particular, $\bar\epsilon_X$ is monic and epic in $\barcat$, and any representative $\epsilon_X \from \floorC{X} \to X$ is a coarse equivalence.
  \end{lemma}

    \proof
     By Lemma \ref{lem:enlarge}, the limit cocone $\bar\lambda \from \Rips_*X \tto \floorC{X}$ in $\barcat$ is also the limit cocone in $\barBorn$.    By Proposition \ref{prop:colimit_coarse}, the canonical cocone $\bar\xi \from \Rips_* X \tto X$ is also the limit cocone in $\barBorn$. Therefore, by the universal property, there is a unique isomorphism $\bar f \from \floorC{X} \to X$ in $\barBorn$ such that $\bar\xi = \bar f \bar\lambda$. Now, $\bar\epsilon_X \from \floorC{X} \to X$ satisfies $\bar\xi = \bar\epsilon_X \bar\lambda$ by definition. Since $\bar\lambda$ is the limit cocone in $\barBorn$, it follows that $\bar\epsilon_X = \bar f$ using the universal property.
    \endproof

    Consequently, by applying Lemma \ref{lem:wlog_metric} to (any representative of) $\bar\epsilon_X^{-1} \from X \to \floorC{X}$, we may assume that $\floorC{X}$ is realised by a metric on $X$ coarsely equivalent to the original.

    \begin{corollary}[Realising Rips colimit]\label{cor:colimit_realise}
     Let $(X,d)$ be a metric space whose Rips colimit in $\barcat$ exists. Then there exists a metric $\partial$ on $X$ such that $\floorC{(X,d)} \cong (X,\partial)$ in $\barcat$, where the Rips counit is a coarse equivalence realised by the underlying identity $\epsilon_X \from (X,\partial) \to (X,d)$. \qed
    \end{corollary}

    \subsection{Rips colimit under dominated controls}\label{sec:dominated_colimit}

    We now characterise the existence of the Rips colimit in a category with dominated controls. The following result, together with Lemma \ref{lem:stable}, yields Theorem \ref{thm:QI-Rips}.

    \begin{theorem}[Rips colimit under dominated controls]\label{thm:dominated-colimit}
     Assume that $\cat$ has dominated controls. Let $X$ be a metric space. If $\floorC{X}$ exists then $X$ is coarsely geodesic.
    \end{theorem}

    Before going into the details of the proof, let us state some immediate corollaries.

  \begin{theorem}\label{thm:colimits2}
   Assume that $\cat$ has dominated controls. Then for a given metric space $X$, items (1) -- (5) in Theorem \ref{thm:colimits} are equivalent to each of the following statements:
   \begin{enumerate}
    \item[(6)] $\Rips_* X$ stabilises and $\bar\epsilon_X \from \floorC{X} \to X$ is an isomorphism in $\barcat$, and
    \item[(7)] $X$ is $\cat$--geodesic.
   \end{enumerate}
  \end{theorem}
  \proof

  (5) $\implies$ (6). This follows from Theorem \ref{thm:dominated-colimit}.

  (6) $\implies$ (7). The assumptions imply that $X \cong \floorC{X} \cong \Rips_\sigma X$ in $\barcat$ for some $\sigma \geq 0$. Since $\Rips_\sigma X$ is a graph, $X$ is $\cat$--geodesic.

  (7) $\implies$ (1). The space $X$ is the colimit of the diagram $X \from \bullet \tto \barcat$, indexed by the trivial category, whose unique object is $X$.
  \endproof

  The implication (1) $\implies$ (7) deserves special mention.

  \begin{corollary}[Colimit-closed]
   Assume that $\barcat$ has dominated controls. Then the subcategory $\Gr(\barcat)$ is closed under taking colimits in $\barcat$. \qed
  \end{corollary}

  As a special case, we deduce Thereom \ref{thm:coclosed} which states that $\barQG$ is colimit-closed in $\barCL$.

  Since the Rips colimit is itself a colimit of graphs, we also obtain the following.

  \begin{corollary}
   Assume that $\barcat$ has dominated controls. If $X$ is a metric space whose Rips colimit in $\barcat$ exists, then $\floorC{X}$ is $\cat$--geodesic. \qed
  \end{corollary}

    We now return our attention to the proof of Theorem \ref{thm:dominated-colimit}. For the remainder of this section, let $(X,d)$ be a metric space such that $\floorC{(X,d)}$ exists. By Lemma \ref{lem:counit_equiv} and Corollary \ref{cor:colimit_realise}, we may assume that $\floorC{(X,d)}$ is realised by $(X,\partial)$ for some metric $\partial$ on $X$ such that the underlying identity $\epsilon_X \from (X,\partial) \to (X,d)$ is a coarse equivalence.

    The standard Rips graph construction assigns each edge unit weight. For this proof, we consider a weighted version of the Rips filtration for $(X,d)$. This shall be used to produce a modified metric $d'$ on $X$ which still admits a cocone $\Rips_*X \tto (X,d')$ via underlying identity maps.

  \begin{definition}[Weighted Rips graph]\label{def:weighted_Rips}
  Let $(X,\partial)$ be a metric space and $\Theta\from [0,\infty) \to [1,\infty)$ be an increasing function. Define $\Rips^\Theta_\infty (X,\partial)$ to be the weighted graph with $X$ as its vertex set, with an edge of weight $\Theta(\partial(x,x'))$ between every distinct pair $x,x' \in X$ satisfying $\partial(x,x') < \infty$. Given $\sigma \geq 0$, let $\Rips^\Theta_\sigma (X,\partial)$ be the subgraph of $\Rips^\Theta_\infty (X,\partial)$ with vertex set $X$, which has an edge between a pair $x,x' \in X$ if and only if $\partial(x,x') \leq \sigma$.   \end{definition}

  The increasing function $\Theta$ shall be called a \emph{weight function}.  We shall regard $\Rips^\Theta_\sigma (X,\partial)$ (resp. $\Rips^\Theta_\infty (X,\partial)$) as a metric space $(X,\partial^\Theta_\sigma)$ (resp. $(X,\partial^\Theta_\infty)$) with underlying set $X$, where the metric $\partial^\Theta_\sigma$ (resp. $\partial^\Theta_\infty$) is given by restriction of the induced path metric to $X$. Note that if $\Theta(t) = 1$ for all $t$ then $\Rips^\Theta_\sigma (X,\partial) = \Rips_\sigma(X,\partial)$.

   \begin{remark}
   If $\Theta$ is assumed to take only positive integer values, then any finite distance $\partial^\Theta_\sigma(x,x')$ is attained by some weighted edge-path in $\Rips^\Theta_\sigma(X,\partial)$ from $x$ to $x'$.
  \end{remark}

    Let us now compare the metric $\partial^\Theta_\infty$ to $\partial$. The following is immediate from construction.

    \begin{lemma}[Weight control]\label{lem:weight_control}
     Let $\Theta$ be a weight function. Then the underlying identity map $(X,\partial) \to (X,\partial^\Theta_\infty)$ has upper control $\Theta$. \qed
    \end{lemma}

    If $\Theta$ is proper, then we also obtain a lower control.
    Define $\Theta^\perp(t) := \sup\{s\geq 0~|~ \Theta(s) \leq t \}$.
    Note that $\Theta(s) \leq t \iff s \leq \Theta^\perp(t)$.

   \begin{lemma}[Proper weight function]\label{lem:proper_ce}
        If $\Theta$ is a proper function then the underlying identity $(X,\partial) \to (X,\partial^\Theta_\infty)$ admits a lower control.     In addition, if $\partial$ is unbounded and does not take the value infinity, then the converse also holds.
    \end{lemma}

    \proof
    Assume $\Theta$ is proper.
    We will show that the underlying identity $(X,\partial^\Theta_\infty) \to (X,\partial)$ admits an upper control. Suppose $x,x' \in X$ are points satisfying $\partial^\Theta_\infty(x,x') = t < \infty$. Then there exists a weighted edge-path $x = x_0, \ldots, x_n = x'$ in $\Rips^\Theta_\infty (X,\partial)$ of length $L < t + 1$. We may assume that all $x_i$ are distinct. Note that $n \leq L$ since all edge weights are at least 1. Moreover, each edge along this path has weight at most $L$. In other words, $\Theta(\partial(x_{i-1}, x_i)) \leq L$ for all $i$. Since $\Theta$ is proper, it follows that $\partial(x_{i-1}, x_i) \leq \Theta^\perp(L) <\infty$ for all $i$. Therefore, by the triangle inequality,
    \[\partial(x,x') \leq \sum_{i=0}^n \partial(x_{i-1}, x_i) \leq n\Theta^\perp(L) \leq (t+1)\Theta^\perp(t+1). \]
    Hence $t\mapsto (t+1)\Theta^\perp(t+1)$ serves as an upper control.

    Now, assume that $\partial$ is unbounded and does not take the value infinity. This implies that $\Rips^\Theta_\infty(X,\partial)$ is (combinatorially) a complete graph. If $\Theta$ is not proper, then all edge weights in $\Rips^\Theta_\infty(X,\partial)$ are uniformly bounded above. Thus, $\partial^\Theta_\infty$ is bounded, and so the underlying identity $(X,\partial) \to (X,\partial^\Theta_\infty)$ cannot admit a lower control.
   \endproof

    Consequently, for any proper weight function $\Theta$, the underlying identity ${(X,\partial) \to (X,\partial^\Theta_\infty)}$ is a coarse equivalence. If this map admits an upper control $\rho$ which eventually remains below $\Theta$ by a definite amount, then we obtain stability for the weighted Rips filtration of $X$.

   \begin{lemma}[Surplus weight]\label{lem:surplus}
    Suppose that the underlying identity ${(X,\partial) \to (X,\partial^\Theta_\infty)}$ has upper control $\rho$. Assume there exists $\sigma \geq 0$ such that $\rho(t) + 1 \leq \Theta(t)$ for all $t \geq \sigma$. Then the underlying identity $\Rips^\Theta_\sigma  (X,\partial) \to \Rips^\Theta_\infty (X,\partial)$ is an isometry.
    \end{lemma}

   \proof
    Let $x,x' \in X$. If $\partial^\Theta_\infty(x,x') = \infty$ then there is no weighted path from $x$ to $x'$ in $\Rips^\Theta_\infty (X,\partial)$. Since  $\Rips^\Theta_\sigma(X,\partial)$ is a subgraph of $\Rips^\Theta_\infty (X,\partial)$, it follows that $\partial^\Theta_\sigma(x,x') = \infty$.

    Now suppose that $\partial^\Theta_\infty(x,x') < \infty$. Let $x = x_0, \ldots, x_n = x'$ be a weighted path in $\Rips^\Theta_\infty (X,\partial)$ of length $L$. We claim that if $L \leq \partial^\Theta_\infty(x,x') + \frac{1}{2}$ then $\partial(x_{i-1}, x_i) \leq \sigma$ for all $i$. This would imply that the given path lies in the subgraph $\Rips^\Theta_\sigma(X,\partial)$. Consequently, the distance $\partial^\Theta_\infty(x,x')$ would equal the infimum of the lengths of all weighted paths from $x$ to $x'$ in $\Rips^\Theta_\sigma(X,\partial)$, thence $\partial^\Theta_\infty(x,x') = \partial^\Theta_\sigma(x,x')$.

    Suppose for a contradiction that $\partial(x_{i-1}, x_i) > \sigma$ for some $i$. Then the edge from $x_{i-1}$ and $x_i$ in $\Rips^\Theta_\infty (X,\partial)$ has weight $\Theta(\partial(x_{i-1}, x_i)) \geq \rho(\partial(x_{i-1}, x_i)) + 1$. Since $\rho$ is an upper control for the underlying identity ${(X,\partial) \to (X,\partial^\Theta_\infty)}$, it follows that $\partial^\Theta_\infty(x_{i-1},x_i) \leq \rho(\partial(x_{i-1}, x_i))$. Therefore, there exists a path from $x_{i-1}$ to $x_i$ in $\Rips^\Theta_\infty (X,\partial)$ of length at most $\rho(\partial(x_{i-1}, x_i)) + \frac{1}{3}$. By replacing the edge from $x_{i-1}$ to $x_i$ in the original path by this detour, we obtain a path from $x$ to $x'$ in $\Rips^\Theta_\infty (X,\partial)$ of length at most $L - \frac{2}{3} \leq \partial^\Theta_\infty(x,x') - \frac{1}{6} < \partial^\Theta_\infty(x,x')$, a contradiction.
   \endproof

   We are now ready to prove the main result.

   \proofof{Theorem \ref{thm:dominated-colimit}}
   By assumption, $\cat = \cat(\F)$ has dominated controls.
   Let $\Theta$ be a proper dominating function for $\F$ and consider the metric space $(X,\partial^\Theta_\infty)$. By Lemma \ref{lem:weight_control}, the underlying identity $f \from (X,\partial) \to (X,\partial^\Theta_\infty)$ admits $\Theta$ as an upper control. Since $(X,\partial)$ is a colimit of graphs in $\barcat$, by Theorem \ref{thm:colimits}, $f$ must admit an upper control $\rho \in \F$. Then $t \mapsto \rho(t) + 1$ also belongs to $\F$, hence there exists some $\sigma \geq 0$ such that $\rho(t) + 1 \leq \Theta(t)$ for all $t \geq \sigma$.
   Appealing to Lemma \ref{lem:surplus}, we deduce that $(X,\partial^\Theta_\infty) = (X,\partial^\Theta_\sigma)$. Since all edge weights of $\Rips^\Theta_\sigma(X,\partial)$ lie in a bounded interval $[1, \Theta(\sigma)]$, it follows that $(X,\partial^\Theta_\sigma)$ is biLipschitz equivalent to a (combinatorial) graph. Thus, $(X,\partial^\Theta_\infty)$ is coarsely geodesic. Since $\Theta$ is proper, we may invoke Lemmas \ref{lem:weight_control} and \ref{lem:proper_ce} to deduce that $(X,\partial)$ is coarsely equivalent to $(X,\partial^\Theta_\infty)$, and hence coarsely geodesic. Consequently, $(X,d)$ is coarsely geodesic.\qed

   Let us conclude this section with a useful observation using the weighted Rips construction. Here, instead of using the (fast growing) dominating function $\Theta$ as the weight function, we use (the slowly growing) $\Theta^\perp$ to introduce many shortcuts.

  \begin{example}[Coarsely but not $\cat$--geodesic]\label{exa:notcg}
   Let $\cat = \cat(\F)$ be a category with dominated controls. Choose a dominating function $\Theta$ for $\F$. Let $(\Gamma, d)$ be an unbounded connected graph equipped with the standard graph metric. Let $d' = d^{\Theta^\perp}_\infty$ be the metric on $\Gamma$ obtained via the weighted Rips graph construction using $\Theta^\perp$ as a (proper) weight function. Then the underlying identity $f \from (\Gamma, d') \to (\Gamma, d)$ is a coarse equivalence by Lemmas \ref{lem:weight_control} and \ref{lem:proper_ce}. Therefore, $(\Gamma, d')$ is coarsely geodesic.    Let $\rho$ be an upper control for $f$. By construction, $d'(x,y) \leq \Theta^\perp(d(x,y))$, and so \[\Theta(d'(x,y)) \leq d(x,y) \leq \rho(d'(x,y))\] for all vertices $x,y \in \Gamma$. Therefore $\rho \not\in \F$ since $d'$ takes unbounded values. Consequently, $f$ is not a morphism in $\cat$, and so $(\Gamma, d')$ is not $\cat$--geodesic by Theorem \ref{thm:colimits}.
  \end{example}

  \section{Adjointness and consequences}\label{sec:adjoint}

  We now establish adjointness properties of the Rips colimit. To set up some notation, write
  \begin{itemize}
   \item $\Gr(\cat)$ (resp. $\Gr(\barcat)$) for the full subcategory of $\cat$ (resp. $\barcat$) whose objects are $\cat$--geodesic spaces,
   \item $\Grhat(\cat)$ (resp. $\Grhat(\barcat)$) for the full subcategory of $\cat$ (resp. $\barcat$) whose objects are those metric spaces which arise as a colimit of graphs in $\barcat$, and
   \item $\Rc(\cat)$ (resp. $\Rc(\barcat)$) for the full subcategory of $\cat$ (resp. $\barcat$) whose objects are the metric spaces $X$ such that $\floorC{X}$ exists.
  \end{itemize}
  Note that $\Gr(\cat)$ is a subcategory of $\Grhat(\cat)$. By Theorem \ref{thm:colimits} and Corollary \ref{cor:idempotence}, $\Grhat(\cat)$ is a subcategory of $\Rc(\cat)$. Some of these subcategories coincide, depending on the ambient category:
  \begin{itemize}
   \item If $\cat$ has dominated controls then, by Theorems \ref{thm:dominated-colimit} and \ref{thm:colimits2},
   \[\Grhat(\cat) = \Gr(\cat) \qquad \textrm{and} \qquad \Rc(\cat) = \CGeod \cap \cat.\]
   In particular, $\Grhat(\CLip) = {\QGeod}$ and $\Rc({\CLip}) = {\CGeod} \cap {\CLip}$.
   \item Proposition \ref{prop:colimit_coarse} implies that $\Grhat({\Born}) = \Rc({\Born}) = \Born$.
  \end{itemize}

  Appealing to Proposition \ref{prop:functorial}, the Rips colimit yields a functor $\floorC{-} \from \Rc(\barcat) \to \Grhat(\barcat)$ (where we have restricted the codomain to the essential image). Write $I \from \Grhat(\barcat) \to \Rc(\barcat)$ for the inclusion functor.

   \begin{proposition}[Naturality of Rips counit]\label{prop:naturality}
    Let $X$ and $Y$ be metric spaces and $f \from X \to Y$ be a controlled map. Then the diagram
    \begin{center}
    \begin{tikzcd}[column sep=large, row sep=large]
    \floorC{X} \arrow[r, "\floorC{\overline{f}}"] \arrow[d, "\bar\epsilon_X"'] & \floorC{Y} \arrow[d, "\bar\epsilon_Y"] \\
    X \arrow[r, "\bar f"'] & Y
    \end{tikzcd}
    \end{center}
    commutes in $\barBorn$, so long as the involved Rips colimits exist.    Consequently, there exists a natural transformation $\bar\epsilon \from I\floorC{-} \tto 1_{\Rc(\barcat)}$     whose component at $X \in \Rc(\barcat)$ is $\bar\epsilon_X \from \floorC{X} \to X$.
    \end{proposition}

    \proof
    Assume $\floorC{X}$ and $\floorC{Y}$ exist. Let $\xi^X \from \Rips_* X \tto X$ and $\xi^Y \from \Rips_* Y \tto Y$ denote the respective canonical cocones. Let $f \from X \to Y$ be a representative for $\overline{f}$, and $\rho$ be an upper control for $f$. We may assume that $\rho(t) \geq t$ for all $t \geq 0$. By Lemma \ref{lem:Rip_control}, $f$ induces a natural transformation $\Rips_* (f,\rho) \from \Rips_*X \tto \Rips_{\rho*}Y$ in $\cat$. Consider the commutative diagram
     \begin{center}
        \begin{tikzcd}[column sep=huge, row sep=large]
    \Rips_* X \arrow[Rightarrow, r, "\Rips_*(f{,}\rho)"] \arrow[Rightarrow, d, "\xi^X"'] & \Rips_{\rho*}Y \arrow[Rightarrow, d, "\xi^Y \rho"] \arrow[Leftarrow, r, "\Rips_*(1_Y{,}\rho)"] & \Rips_*Y \arrow[Rightarrow, d, "\xi^Y"] \\
    X \arrow[r, "f"'] & Y \arrow[leftarrow, r, "1_Y"', "\sim"] & Y
    \end{tikzcd}
    \end{center}
    in $\cat$, where each object is a diagram indexed by $[0,\infty)$. By working in $\barcat$ and taking colimits, we obtain a commutative diagram
    \begin{center}
    \begin{tikzcd}[column sep=large, row sep=large]
    \floorC{X} \arrow[r, "\floorC{\overline{f}}"] \arrow[d, "\bar\epsilon_X"'] & \floorC{Y} \arrow[d, "\colimit \xi^Y \rho"] \arrow[leftrightarrow, r, "\bar 1_{\floorC{Y}}", "\sim"'] & \floorC{Y} \arrow[d, "\bar\epsilon_Y"]\\
    X \arrow[r, "\bar f"'] & Y \arrow[leftrightarrow, r, "\bar 1_Y"', "\sim"] & Y
    \end{tikzcd}
    \end{center}
    as desired, where the identity $\bar 1_{\floorC{Y}}$ realises the canonical isomorphism between the common colimit $\floorC{Y}$ for $\Rips_*Y$ and $\Rips_{\rho*}Y$ as given by Corollary \ref{cor:adjust-colimit}.
    \endproof

    \begin{theorem}[Rips colimit is a right adjoint]\label{thm:adjoint}
     The Rips colimit $\floorC{-} \from \Rc(\barcat) \to \Grhat(\barcat)$ is right adjoint to the inclusion $I \from \Grhat(\barcat) \to \Rc(\barcat)$, where the counit is given by $\bar\epsilon$.
    \end{theorem}

     \proof
     We need to show that for all $X\in\Grhat({\barcat})$ and $Y \in \Rc(\barcat)$, the function on hom-sets
     \[(\bar\epsilon_Y)_\# \from \Grhat({\barcat})(X,\floorC{Y}) \to \Rc(\barcat)(X, Y) \]
     given by $(\bar\epsilon_Y)_\# \bar f := \bar\epsilon_Y \bar f$ is a bijection. Naturality in $X$ and $Y$ follows from naturality of $\bar\epsilon$.

     By Lemma \ref{lem:counit_equiv}, $\bar\epsilon_Y$ is a monomorphism in $\barcat$ and so $(\bar \epsilon_Y)_\#$ is injective. Let $\bar f \in \Rc(\barcat)(X, Y)$, that is, a morphism $\bar f \from X \to Y$ in $\barcat$. Since $X\in\Grhat({\barcat})$, it follows from Theorem \ref{thm:colimits} that $\bar\epsilon_X$ is an isomorphism in $\barcat$. Therefore, we may define a morphism $\bar f':= \floorC{\bar f}\bar\epsilon_X^{-1} \from X \to \floorC{Y}$ in $\barcat$, yielding an element of $\Grhat({\barcat})(X,\floorC{Y})$. By Proposition \ref{prop:naturality}, we deduce that
     \[(\bar\epsilon_Y)_\# \bar f' = (\bar\epsilon_Y)_\# \floorC{\bar f}\bar\epsilon_X^{-1} = \bar\epsilon_Y \floorC{\bar f} \bar\epsilon_X^{-1} = \bar f \bar\epsilon_X \bar\epsilon_X^{-1} = \bar f, \]
     hence $(\bar\epsilon_Y)_\#$ is surjective.
    \endproof

     We shall refer to $\bar\epsilon$ as the \emph{Rips counit} in $\barcat$.       The unit $\bar\eta \from 1_{\Grhat(\cat)} \tto \floorC{I-}$ of this adjunction is given componentwise by $\bar\eta_X := \bar\epsilon_X^{-1}$ for $X \in \Grhat(\cat)$.

     Using the standard fact that right adjoints preserve limits \cite[Theorem 4.5.2]{Rie16}, we immediately obtain the following.

    \begin{corollary}
     The Rips colimit $\floorC{-} \from \Rc(\barcat) \to \Grhat(\barcat)$ is faithful and preserves limits. \qed
    \end{corollary}

     If $\cat$ has dominated controls, then Example \ref{exa:notcg} shows that there exist coarsely geodesic spaces which are not $\cat$--geodesic. This implies that
     \[I \from \Grhat(\barcat) = \Gr(\barcat) \to \barCG \cap \barcat = \Rc(\barcat)\]
     is not an equivalence of categories. To remedy this, we should allow for coarse equivalences when dealing with coarsely geodesic spaces. In other words, we consider instead the inclusion
     \[I \from \Grhat(\barcat) = \Gr(\barcat) \to \barCG = \Rc(\barBorn).\]
     Note that this is essentially surjective since every coarsely geodesic space is coarsely equivalent to a graph. For the other direction, observe that, in the proof of Theorem \ref{thm:adjoint}, we never needed the fact that $\bar f \from X \to Y$ was a morphism in $\barcat$, only that it was controlled. Indeed, we only need that $f$ is controlled to define the morphism $\floorC{\bar f}$ in $\barcat$ (Definition \ref{def:Rips_f}) and for Proposition \ref{prop:naturality} to hold. Consequently, we may replace $\Rc(\barcat)$ with $\barCG$ to obtain an alternative Rips colimit functor
     \[\floorC{-} \from \barCG \to \Grhat(\barcat) = \Gr(\barcat)\]
     and Rips counit $\bar\epsilon \from I\floorC{-} \tto 1_{\barCG}$. We show that this remedy indeed promotes the adjunction to an adjoint equivalence when $\barcat$ has dominated controls.

  \begin{theorem}[Adjoint equivalence]\label{thm:adjeq}
  Assume that $\barcat$ has dominated controls. Then the pair $I \from \Gr(\barcat) \to \barCG$ and $\floorC{-} \from \barCG \to \Gr(\barcat)$ form an adjoint equivalence.
  \end{theorem}

  \proof
  Since $\barcat$ has dominated controls, $\Grhat(\barcat) = \Gr(\barcat)$ and the objects of $\Rc(\barcat)$ are precisely the coarsely geodesic spaces by Theorems \ref{thm:dominated-colimit} and \ref{thm:colimits2}. Using the same proof as for Theorem \ref{thm:adjoint}, with $\barCG$ in place of $\Rc(\barcat)$, we deduce that $\floorC{-} \from \barCG \to \Gr(\barcat)$ is right-adjoint to $I \from \Gr(\barcat) \to \barCG$. Let us verify that the associated unit and counits are natural isomorphisms. For any $Y \in \barCG$, the map $\epsilon_Y \from \floorC{Y} \to Y$ is a coarse equivalence by Lemma \ref{lem:counit_equiv}. Therefore, $\bar\epsilon \from I\floorC{-} \tto 1_{\barCG}$ is a natural isomorphism (since we consider the components in $\barBorn$). The unit $\bar\eta \from 1_{\Gr(\barcat)} \tto \floorC{I-}$ is a natural isomorphism since, by definition, the component $\bar\eta_X$ at $X \in \Grhat(\barcat)$ is given by the inverse of $\bar\epsilon_X$ in $\barcat$ (which exists by Theorem \ref{thm:colimits}).
  \endproof

   Theorem \ref{thm:adjoint0} follows from Theorem \ref{thm:adjoint}, Example \ref{exa:notcg}, and Theorem \ref{thm:adjeq}.

   \subsection{Universal morphisms}

   We now apply adjointness to give a proof of Theorem \ref{thm:univ-morph}.

   Consider the problem of finding a ``best approximation'' to a metric space $X$ using controlled maps from a graph or, equivalently, from a $\cat$--geodesic space. Say $\bar f \from W \to X$ is a \emph{universal morphism from $\Gr(\cat)$} if $W \in \Gr(\cat)$ and for any $\bar f' \from W' \to X$ where $W' \in \Gr(\cat)$, there exists a unique morphism $\bar g \from W' \to W$ such that $\bar f' = \bar f \bar g$;  moreover, such $\bar f'$ is universal if and only if $\bar g$ is an isomorphism.\footnote{Equivalently, $\bar f$ is the terminal object in the full subcategory of the slice category $\barcat / X$ whose objects are (closeness classes of) controlled maps to $X$ with domain in $\Gr(\cat)$.}

  \begin{proposition}[Universal morphisms from $\Gr(\cat)$]\label{prop:universal}
   Let $X$ be a metric space. Then $X$ admits a universal morphism from $\Gr(\cat)$ if and only if $X$ is coarsely geodesic. In particular, a morphism $\bar f \from W \to X$, where $W \in \Gr(\cat)$, is universal if and only if $f$ is a coarse equivalence.
  \end{proposition}

  \proof
  Suppose that $\bar f \from W \to X$ is a universal morphism from $\Gr(\cat)$.   Then for all $\sigma \geq 0$, there exists a unique morphism $\bar \mu_\sigma \from \Rips_\sigma X \to W$ such that $\bar\xi_\sigma = \bar f \bar \mu_\sigma$. Given $\sigma \leq \tau$, write $\iota \from \Rips_\sigma X \to \Rips_\tau X$ for the map coinciding with the underlying identity. Then
  \[\bar f \bar \mu_\sigma = \bar\xi_\sigma = \bar\xi_\tau\bar\iota = \bar f \bar \mu_\tau\bar\iota,\]
  hence $\bar \mu_\sigma = \mu_\tau\bar\iota$ by the universal property of $\bar f$. Consequently, the $\bar\mu_\sigma$ assemble to form the colegs of a cocone $\bar\mu \from \Rips_* X \tto W$ in $\barcat$ satisfying $\bar\xi = \bar f \bar\mu$.

    \begin{center}
    \tikzcdset{arrow style=math font}
    \begin{tikzcd}[column sep=huge, row sep=large]
    \Rips_* X \arrow[Rightarrow, r, "\bar\mu"] \arrow[Rightarrow, rd, "\bar\xi"'] & W \arrow[rightarrow, d, "\bar f", bend left = 10] \arrow[leftarrow, d, "\bar h"', bend right=10] \\
    & X
    \end{tikzcd}
    \end{center}

  By Proposition \ref{prop:colimit_coarse}, $\bar\xi \from \Rips_* X \tto X$ is the limit cocone in $\barBorn$, and so there exists a unique morphism $\bar h \from X \to W$ in $\barBorn$ such that $\bar\mu = \bar h\bar\xi$. Then $\bar f\bar h \bar\xi = \bar f \bar \mu = \bar\xi$, hence $\bar f\bar h = \bar 1_X$ by the universal property of $\bar\xi$. Now, $\bar h \bar f \from W \to W$ is a morphism in $\barcat$ since $W$ is $\cat$--geodesic. Therefore, $\bar h \bar f = \bar 1_W$ by the universal property of $\bar f$. It follows that $\bar f \from W \to X$ is a coarse equivalence, and so $X$ is coarsely geodesic.

  For the converse, assume that $X$ is coarsely geodesic. Then $\floorC{X}$ exists and is an object in $\Gr(\cat)$. Let $\bar f \from W \to X$ be a morphism from some $W \in \Gr(\cat)$. Note that $W$ is also an object in $\Grhat(\cat)$. By Theorem \ref{thm:adjoint}, there exists a unique morphism $\bar f' \from W \to \floorC{X}$ satisfying $\bar f = \bar\epsilon_X \bar f'$, hence $\bar\epsilon_X$ is a universal morphism from $\Gr(\cat)$.

   \begin{center}
   \tikzcdset{arrow style=math font}
    \begin{tikzcd}[column sep=huge, row sep=large]
    W \arrow[rightarrow, r, "\bar f'"] \arrow[rightarrow, rd, "\bar f"'] & \floorC{X} \arrow[rightarrow, d, "\bar \epsilon_X'"] \\
    & X
    \end{tikzcd}
  \end{center}

  Note that $\bar f$ is itself a universal morphism from $\Gr(\cat)$ if and only if $\bar f'$ is an isomorphism in $\barcat$. Now, $W$ and $\floorC{X}$ are both colimits of $\cat$--geodesic spaces in $\barcat$ and so, by Theorem \ref{thm:colimits}(3), $\bar f'$ being an isomorphism in $\barcat$ is equivalent to it being so in $\barBorn$. Since $\bar\epsilon_X$ is an isomorphism in $\barBorn$, by Lemma \ref{lem:counit_equiv}, this holds precisely when $\bar f$ is an isomorphism in $\barBorn$.
  \endproof

  \begin{corollary}[Realising universal morphisms]\label{cor:univ_mor}
   Let $X$ be a coarsely geodesic space. Then $\bar\xi_\sigma \from \Rips_\sigma X \to X$ is a universal morphism from $\Gr(\cat)$ for $\sigma \geq 0$ sufficiently large. \qed
  \end{corollary}

  \subsection{Extremality of the stable Rips colimit}\label{sec:extreme}

  In \cite{Ros22}, Rosendal defines a partial order on the set of compatible left-invariant metric on a topological group. Rosendal shows that a Polish group admits a maximal such metric, with respect to this partial order, if and only if it admits a left-invariant coarsely geodesic metric; moreover, the maximal metric is realised by a left-invariant quasigeodesic metric, unique up to quasi-isometry. Consequently, such groups are endowed with a canonical quasi-isometry class of quasigeodesic metrics.

  We draw inspiration from Rosendal's work to define a preorder on the set of metrics on a set $X$. Let $\cat$ be a category associated ot a control class $\F$. Given metrics $d,d'$ on $X$, declare $d \prec_{\cat} d'$ if and only if the identity map $(X,d') \to (X,d)$ is a morphism in $\cat$. Explicitly, this means that there exists some $\rho \in \F$ such that $d(x,y) \leq \rho d'(x,y)$ for all $x,y\in X$. Equivalence of two metrics, denoted by $\asymp_{\cat}$, under this preorder holds exactly when the (closeness class of the) underlying identity is an isomorphism in $\barcat$.

  \begin{remark}
  Rosendal's partial order is defined using quasimetrics (quasi-isometry classes of metrics) rather than metrics. Thus, we recover Rosendal's definition by taking the quotient partial order in the case where $\cat=\CLip$.
  \end{remark}

  For each $\sigma \geq 0$, write $\Rips_\sigma (X,d) = (X,d_\sigma)$. Each $\Rips_\sigma (X,d)$ is a graph, hence belongs to $\Gr(\cat)$. Moreover, these metrics satisfy $d \prec_{\cat} d_\sigma \prec_{\cat} d_\tau$ for all $0 \leq \tau \leq \sigma$.

  Given a metric $d$ on $X$, consider the set of $\cat$--geodesic metrics $d'$ which $\prec_{\cat}$--upper bound $d$. We show that there exists a $\prec_{\cat}$--least element in this set precisely when $(X,d)$ is coarsely geodesic. Note that any  such least element must be unique up to $\asymp_{\cat}$--equivalence.

  \begin{proposition}[Minimal $\cat$--geodesic upper bound]
   Let $(X,d)$ be a metric space. Then the following are equivalent:
   \begin{enumerate}
    \item There is a $\prec_{\cat}$--least element among all $\cat$--geodesic metrics $\partial$ satisfying  $d \prec_{\cat} \partial$,
    \item There is a $\prec_{\cat}$--minimal element among all $\cat$--geodesic metrics $\partial$ satisfying  $d \prec_{\cat} \partial$, and
    \item $(X,d)$ is coarsely geodesic.
   \end{enumerate}
  Furthermore, the minimal metric (if it exists) is attained by $\Rips_\sigma (X,d)$ for $\sigma$ large.
  \end{proposition}

  \proof

  (1) $\implies$ (2). Trivial.

  (2) $\implies$ (3).
  Assume that $\partial$ is minimal among $\cat$--geodesic metrics $d'$ on $X$ satisfying $d \prec_{\cat} d'$. Write $\Rips_\sigma (X,\partial) = (X,\partial_\sigma)$ for $\sigma \geq 0$. Appealing to Proposition \ref{prop:C-geodesic}, we deduce that $\partial \asymp_{\cat} \partial_\sigma$ for some $\sigma \geq 0$. The underlying identity $(X,\partial) \to (X,d)$ admits an upper control $\rho \in \F$, hence $d_{\rho\sigma} \prec_{\cat} \partial_\sigma $ by Lemma \ref{lem:Rip_control}. Therefore,
  \[d_\tau \prec_{\cat} d_{\rho \sigma} \prec_{\cat} \partial_\sigma \prec_{\cat} \partial \]
  for all $\tau \geq \rho\sigma$. Consequently, $\partial \asymp_{\cat} d_\tau$ by the minimality assumption on $\partial$. Since this holds for any $\tau \geq \rho\sigma$, it follows that $\Rips_* (X,d)$ stabilises, hence $(X,d)$ is coarsely geodesic by Lemma \ref{lem:stable}.

  (3) $\implies$ (1).
  Assume $(X,d)$ is coarsely geodesic. By Corollary \ref{cor:univ_mor}, $\bar\xi_\sigma \from (X, d_\sigma) \to (X,d)$ is a universal morphism from $\Gr(\cat)$ for $\sigma \geq 0$ large. Let $d'$ be a $\cat$--geodesic metric on $X$ such that $d \prec_{\cat} d'$, and write $f \from (X,d') \to (X,d)$ for the underlying identity. Then there exists a unique morphism $\bar f ' \from (X,d') \to (X,d_\sigma)$ in $\barcat$ such that $\bar f = \bar\xi_\sigma \bar f'$. Now, $\xi_\sigma$ is a coarse equivalence, hence $\xi^{-1}_\sigma$ is controlled. Therefore, the underlying identity $\xi_\sigma^{-1}f \from (X,d') \to (X,d_\sigma)$ is controlled. Since $d'$ is $\cat$--geodesic, it follows using Theorem \ref{thm:colimits2} that $d_\sigma \prec_{\cat} d'$. Thus, $d_\sigma$ is $\prec_{\cat}$--least among $\cat$--geodesic metrics on $X$ which $\prec_{\cat}$--upper bound $d$.
  \endproof

  Let us now present a counterpart result for maximal metrics analogous to Rosendal's characterisation.

  \begin{proposition}[Maximal coarsely equivalent metric]
   Assume $\cat$ has dominated controls. Let $(X,d)$ be a metric space. Then the following are equivalent:
   \begin{enumerate}
    \item there exists a $\prec_{\cat}$--greatest element among metrics $\partial$  on $X$ satisfying $\partial \asymp_{\Born} d$,
    \item there exists a $\prec_{\cat}$--maximal element among metrics $\partial$ on $X$ satisfying $\partial \asymp_{\Born} d$,
    \item $(X,d)$ is coarsely geodesic.
   \end{enumerate}

  Moreover, the $\prec_{\cat}$--maximal metric (if it exists) is attained by $\Rips_\sigma (X,d)$ for $\sigma$ large.
  \end{proposition}

  \proof

   (1) $\implies$ (2). Trivial.

   (2) $\implies$ (3). Suppose that $\partial$ is $\prec_{\cat}$--maximal among metrics on $X$ coarsely equivalent to $d$. Let $\Theta$ be a dominating function for $\F$. By replacing $\Theta$ with $t \mapsto \Theta(t) + t + 1$, we may assume that the following hold: $\Theta(t) \geq t$ for all $t\geq 0$ and for any $\rho \in \F$, there exists $\sigma \geq 0$ such that $\Theta(t) \geq \rho(t) + 1$ for all $t \geq \sigma$. Let $(X,\partial^\Theta_\infty) = \Rips^\Theta_\infty(X,d)$ be the space obtained by the weighted Rips graph contruction (Definition \ref{def:weighted_Rips}). By Lemmas \ref{lem:weight_control} and \ref{lem:proper_ce}, we deduce that $\partial^\Theta_\infty \asymp_{\Born} \partial \asymp_{\Born} d$.   Therefore $\partial^\Theta_\infty \prec_{\cat} \partial$ by the maximality assumption on $\partial$, and so the underlying identity $(X,\partial) \to (X, \partial^\Theta_\infty)$ admits an upper control $\rho \in \F$. Appealing to Lemma \ref{lem:surplus}, there exist some $\sigma \geq 0$ such that the underlying identity $ \Rips^\Theta_\sigma(X,\partial) \to \Rips^\Theta_\infty(X,\partial)$ is an isometry. Since $\Rips^\Theta_\sigma(X,\partial)$ is biLipschitz equivalent to a graph, it follows that $(X,\partial^\Theta_\infty)$, hence $(X,d)$, is coarsely geodesic.

   (3) $\implies$ (1). Assume $(X,d)$ is coarsely geodesic. Then $d_\sigma \asymp_{\Born} d$ for some $\sigma \geq 0$ by Lemma \ref{lem:stable}. Suppose that $\partial$ is a metric on $X$ such that $\partial \asymp_{\Born} d$. Then the underlying identity $(X,d_\sigma) \to (X, \partial)$ is controlled. Since $(X,d_\sigma)$ is a graph, it follows that $\partial \prec_{\cat} d_\sigma$. Therefore, $d_\sigma$ is $\prec_{\cat}$--greatest (hence $\prec_{\cat}$--maximal) among metrics on $X$ coarsely equivalent to $d$.
  \endproof

  In conclusion, if $\cat$ has dominated controls then, under the $\prec_{\cat}$--ordering, the least $\cat$--geodesic upper bound for $d$ and the greatest metric $\asymp_{\Born}$--equivalent to $d$ both coincide with the stable Rips colimit of $(X,d)$.

  This completes the proof of Theorem \ref{thm:extremal}.

  \section{Universal \texorpdfstring{$\cat$}{$\mathsf{C}$}--geodesic cones}\label{sec:universal}

  In this section, we apply results from Section \ref{sec:adjoint} to compute universal $\cat$--geodesic cones. We shall first give a general two-step procedure, and then specialise to the case of equalisers. In Section \ref{sec:uc-cone}, we shall treat the uniformly controlled version.

  Let us prove Theorem \ref{thm:univ-qgc}. Consider a diagram $D \from \sJ \to \barBorn$. Here, we do not assume the objects $D_j$ are $\cat$--geodesic. A \emph{$\cat$--geodesic cone} above $D$ is a cone $\bar \mu\from X \tto D$ in $\barBorn$ whose apex $X$ is $\cat$--geodesic. A $\cat$--geodesic cone $\bar \mu\from X \tto D$ is \emph{universal} if for any $\cat$--geodesic cone $\bar \mu'\from X' \tto D$, there exists a unique morphism $\bar f \from X' \to X$ in $\barBorn$ such that $\bar\mu' = \bar\mu \bar f$. Note that any $\cat$--geodesic cone must be a cone in $\barcat$.

  \begin{theorem}[Universal $\cat$--geodesic cones]\label{thm:cones}
   Let $D \from \sJ \to \barBorn$ be a diagram. Then $D$ admits a universal $\cat$--geodesic cone if and only if $\limit_\sJ D$ exists in $\barBorn$ and is coarsely geodesic. Moreover, a universal $\cat$--geodesic cone for $D$ (if it exists) is realised by $\floorC{\limit_\sJ D} \cong \Rips_\sigma \limit_\sJ D$ for $\sigma \geq 0$ large.
  \end{theorem}

  \proof
  Assume $\bar\lambda \from W \tto D$ is a universal $\cat$--geodesic cone. We claim that $\bar\lambda$ is a limit cone in $\barBorn$. Suppose that $\bar\nu \from X \tto D$ is an arbitrary cone in $\barBorn$. Let $\bar\xi \from \Rips_* X \tto X$ be the associated canonical cocone. Then $\bar\nu\bar\xi_\sigma \from \Rips_\sigma X \tto D$ is a $\cat$--geodesic cone for each $\sigma \geq 0$. By the universal property of $\bar\lambda$, there exists a unique morphism $\bar \mu_\sigma \from \Rips_\sigma X \to W$ in $\barBorn$ such that $\bar\nu \bar\xi_\sigma= \bar\lambda \bar \mu_\sigma$. For all $\sigma \leq \tau$, we have that
  \[\bar\lambda \bar \mu_\sigma = \bar\nu \bar\xi_\sigma = \bar\nu \bar\xi_\tau \bar\iota = \bar\lambda \bar \mu_\tau \bar\iota,\]
  where $\bar\iota \from \Rips_\sigma X \to \Rips_\tau X$ is the inclusion, hence $\bar \mu_\sigma = \bar \mu_\tau \bar\iota$. Therefore, there exists a cocone $\bar \mu \from \Rips_* X \tto W$ whose colegs are the $\bar \mu_\sigma$.

  \begin{center}
   \begin{tikzcd}[column sep=huge, row sep=huge] \Rips_\sigma X \arrow[r, "\bar\mu_\sigma", dashed] \arrow[d, "\bar\xi_\sigma" swap] \arrow[dr, "\bar\nu_j\bar\xi_\sigma"] & W \arrow[d, "\bar\lambda_j"] \\ X \arrow[r, "\bar \nu_j" swap]  & D_j  \end{tikzcd}
   \qquad\qquad
   \begin{tikzcd}[column sep=huge, row sep=huge] \Rips_\sigma X \arrow[r, "\bar\mu_\sigma"] \arrow[d, "\bar\xi_\sigma" swap]  & W \arrow[d, "\bar\lambda_j"] \\ X \arrow[ur, "\bar m", dashed] \arrow[r, "\bar \nu_j" swap] & D_j   \end{tikzcd}
  \end{center}

  By Proposition \ref{prop:colimit_coarse}, $\bar m := \colimit \bar\mu \from X \to W$ is the unique morphism in $\barBorn$ that satisfies $\bar\mu = \bar m \bar\xi$. Now, for each $\sigma \geq 0$ and $j \in \Obj(\sJ)$, we have that
  \[\bar\nu_j \bar\xi_\sigma = \bar\lambda_j\bar\mu_\sigma = \bar\lambda_j \bar m \bar\xi_\sigma.\]
  Since $\bar\xi_\sigma$ is an epimorphism in $\barBorn$ (as $\xi_\sigma$ is surjective), we deduce that $\bar\nu_j = \bar\lambda_j \bar m$ for all $j \in \Obj(\sJ)$, hence $\bar\nu = \bar\lambda \bar m$ factors through $\bar\lambda$.

  To verify uniqueness, suppose $\bar m'\from X \to W$ is a morphism in $\barBorn$ satisfying $\bar\nu = \bar\lambda \bar m'$. Then for all $\sigma \geq 0$, we obtain the equality
  \[\bar\lambda \bar m' \bar\xi_\sigma = \bar\nu\bar\xi_\sigma = \bar\lambda\bar\mu_\sigma = \bar\lambda \bar m \bar\xi_\sigma\]
  of $\cat$--geodesic cones. By the universal property of $\bar\lambda$, we deduce that $\bar m' \bar\xi_\sigma= \bar m \bar\xi_\sigma$, and hence $\bar m'= \bar m$ as desired.

  For the converse, assume that $\bar\lambda \from X \tto D$ is a limit cone in $\barBorn$ with coarsely geodesic apex. By Corollary \ref{cor:univ_mor}, $\bar\xi_\sigma \from \Rips_\sigma X\to X$ is a universal morphism from $\Gr(\barcat)$ for some fixed $\sigma \geq 0$ large. We claim that $\bar\lambda\bar\xi_\sigma \from \Rips_\sigma X \tto D$ is a universal $\cat$--geodesic cone.
  \begin{center}
   \begin{tikzcd}[column sep=huge, row sep=large] & W \arrow[Rightarrow, dr,"\bar\mu"] \arrow[d,"\bar f" swap, dashed] \arrow[dl,"\bar g" swap, dashed] &  \\ \Rips_\sigma X \arrow[r,"\bar\xi_\sigma"] \arrow[Rightarrow, rr, "\bar \lambda \bar \xi_\sigma", bend right=20, swap] & X \arrow[Rightarrow, r, "\bar \lambda"] & D\end{tikzcd}
  \end{center}
  Let $\bar\mu \from W \tto D$ be a $\cat$--geodesic cone. By the universal property of $\bar\lambda$, there exists a unique morphism $\bar f \from W \to X$ in $\barBorn$ such that $\bar\mu = \bar \lambda \bar f$. By the universal property of $\bar\xi_\sigma$, there exists a unique morphism $\bar g \from W \to \Rips_\sigma X$ such that $\bar f = \bar\xi_\sigma \bar g$.  Therefore $\bar\mu = \bar\lambda \bar f= \lambda \bar\xi_\sigma \bar g$. To verify uniqueness, suppose $\bar g' \from W \to \Rips_\sigma X$ is a morphism satisfying $\bar\mu = \lambda \bar\xi_\sigma \bar g'$. Since $\lambda$ is a limit cone, the equality $\lambda \bar\xi_\sigma \bar g = \lambda \bar\xi_\sigma \bar g'$ implies that $\bar\xi_\sigma \bar g = \bar\xi_\sigma \bar g'$. By Proposition \ref{prop:universal}, $\bar\xi_\sigma$ is realised by a coarse equivalence and is thus an isomorphism in $\barBorn$, hence $\bar g = \bar g'$.
  \endproof

  This turns the existence problem for universal $\cat$--geodesic cones into a two-step process: first, compute the limit in $\barBorn$ and then (if the limit exists) check whether it is coarsely geodesic. If we only need to show that a given candidate $\cat$--geodesic space $X$ serves as a desired $\cat$--geodesic cone, it suffices to check that $X$ realises $\limit_\sJ D$ in $\barBorn$.

  \begin{corollary}\label{cor:uc-cone-check}
   Let $D \from \sJ \to \barBorn$ be a diagram. Suppose that $X$ is a $\cat$--geodesic space such that $X \cong \limit_\sJ D$ in $\barBorn$. Then $X$ is a universal $\cat$--geodesic cone for $D$.
  \end{corollary}

  \proof
  Since $X$ is $\cat$--geodesic, we deduce that $X \cong \Rips_\sigma X$ in $\barcat$ for $\sigma \geq 0$ large using Proposition \ref{prop:C-geodesic}. The result follows from Theorem \ref{thm:cones}.
  \endproof

  \subsection{Equalisers}

  In any category, the \emph{equaliser} of a parallel pair $f,g \from X \rightrightarrows Y$ is a morphism $h \from W \to X$ satisfying $fh = gh$ with the universal property: whenever a morphism $k \from Z \to X$ satisfies $fk = gk$ there exists a unique morphism $k' \from Z \to W$ such that $k = hk'$.

  We now address the problem of computing the universal $\cat$--geodesic cone over a parallel pair $\bar f, \bar g \from X \rightrightarrows Y$ in $\barBorn$ (if it exists).   This is a morphism $\bar h \from W \to X$ from a $\cat$--geodesic space $W$ satisfying $\bar f \bar h = \bar g \bar h$ which is universal among all such morphisms.

  Let us recall some notions and results from \cite{Tang-mono} regarding equalisers in $\barcat$. Given a \emph{consistency parameter} $\kappa \geq 0$ and two functions $f,g \from X \to Y$ between metric spaces, define their \emph{$\kappa$--equaliser} to be
  \[\Eq^\kappa(f,g) := \{x\in X \st d_Y(fx, gx) \leq \kappa \}\]
  equipped with the induced metric from $X$. The inclusion maps $\iota^{\kappa'\kappa} \from \Eq^\kappa(f,g) \hookrightarrow \Eq^{\kappa'}(f,g)$ defined for $0 \leq \kappa \leq \kappa'$, are isometric embeddings. These assemble to form the \emph{equaliser filtreation}
  \[\Eq^*(f,g) \from [0,\infty) \to \barcat.\]
  Note that the inclusions $\iota^\kappa \from \Eq^\kappa(f,g) \hookrightarrow X$, for all $\kappa \geq 0$, are isometric embeddings.

  \begin{lemma}[Factoring through $\kappa$--equalisers {\cite[Lemma 4.2]{Tang-mono}}]
  Let ${\bar f, \bar g \from X \rightrightarrows Y}$ be a parallel pair in $\barBorn$, with representatives $f\in \bar f$, $g \in \bar g$. If $\bar h \from W \to X$ is a morphism in $\barBorn$ satisfying $\bar f \bar h = \bar g \bar h$ then, for all $\kappa \geq 0$ large, there is a unique morphism $\bar h' \from W \to \Eq^\kappa(f,g)$ such that $\bar h = \bar\iota^\kappa \bar h$. \qed
  \end{lemma}

  By considering the colimit of the equaliser filtration, we can attempt to compute the equaliser of a parallel pair $\bar f, \bar g$. The following result was originally proven for $\barcat = \barBorn$ or $\barCL$, however, the proof works for any category $\barcat$ associated to a control class.

  \begin{proposition}[Equalisers in $\barcat$ {\cite[Proposition 4.3]{Tang-mono}}]\label{prop:equaliser}
     Let ${\bar f, \bar g \from X \rightrightarrows Y}$ be a parallel pair in $\barcat$, and let $f\in \bar f$, $g \in \bar g$ be representatives.
    Then the following are equivalent:
      \begin{enumerate}
       \item There exists $\kappa \geq 0$ such that for all $\kappa' \geq \kappa$, there exists $r \geq 0$ such that $\Eq^{\kappa'}(f,g)$ lies in the (closed) metric $r$--neighbourhood of $\Eq^{\kappa}(f,g)$ in $X$,
       \item $\Eq^*(f,g)$ stabilises in $\barcat$,
       \item $\colimit \Eq^*(f,g)$ exists in $\overline{\cat}$,
       \item The equaliser of $\bar f, \bar g$ exists in $\barcat$,
       \item The equaliser for $\bar f, \bar g$ is realised by $\bar\iota^\kappa \from \Eq^\kappa(f,g)\to X$ for sufficiently large $\kappa \geq 0$. \qed
      \end{enumerate}
      \end{proposition}

     Informally, stability means that, up to closeness, no new ``solutions'' $x\in X$ to $fx \approx_\kappa gx$ appear beyond some consistency threshold. Note that stability of the equaliser filtration does not depend on choice of category $\barcat$. Moreover, $\Eq^\kappa(f,g)$ realises the equaliser of $\bar f, \bar g$ for $\kappa \geq 0$ large for any choice of $\barcat$. 

    Next, we consider the interaction between the Rips and equaliser filtrations. Specifically, we define the \emph{Rips--equaliser filtration} of $f, g$ to be the doubly-indexed system
    \[\Rips_*\Eq^*(f,g) \from [0,\infty) \times [0,\infty) \to \barcat,\]
    where the arrows are given by the $1$--Lipschitz maps $\Rips_\sigma \Eq^\kappa(f,g) \rightarrow \Rips_{\sigma'} \Eq^{\kappa'}(f,g)$ coinciding with the underlying inclusions, defined whenever $\sigma \leq \sigma'$ and $\kappa\leq\kappa'$. Here, our notational convention is to write the scale parameter $\sigma$ as the subscript, and the consistency parameter $\kappa$ as the superscript.

    We exhibit a natural relationship between the Rips and Rips--equaliser filtrations in the metric coarse category. For $\kappa \geq 0$, write $\zeta^\kappa \from \Rips_*\Eq^\kappa(f,g) \tto \Eq^\kappa(f,g)$ for the canonical cocone of $\Eq^\kappa(f,g)$.  Given a metric space $Z$ and cocone $\bar\mu \from \Eq^*(f,g) \tto Z$ in $\barBorn$, define $R_Z\bar\mu \from \Rips_*\Eq^*(f,g) \tto Z$ to be the cocone whose colegs are given by $(R_Z\bar\mu)_{\sigma}^{\kappa} := \bar\mu^\kappa\bar\zeta_\sigma^\kappa$. The map $\bar\mu \mapsto R_Z\bar\mu$ defines a function $R_Z \from \Cone(\Eq^*(f,g), Z) \to \Cone(\Rips_*\Eq^*(f,g), Z)$ between the respective sets of cocones in $\barBorn$ with nadir $Z$.

    \begin{center}
    \begin{tikzcd}[column sep=huge,row sep=huge]
    \Rips_\sigma \Eq^\kappa(f,g)    \arrow[r, "\bar\zeta^\kappa_\sigma"] \arrow[dr, "(R_Z\bar\mu)_{\sigma}^{\kappa} = \bar\mu^\kappa\bar\zeta_\sigma^\kappa" swap]  & \Eq^\kappa(f,g) \arrow[d, "\bar\mu^\kappa"] \\   & Z \end{tikzcd}
    \end{center}

     \begin{lemma}
     There exists a natural isomorphism
     \[R \from \Cone(\Eq^*(f,g), -) \tto \Cone(\Rips_*\Eq^*(f,g), -)\]
     in $\barBorn$ whose components are $R_Z$ for each metric space $Z$.
    \end{lemma}

    \proof
    Naturality is straightforward. We wish to show that $R_Z$ is a bijection.    For injectivity, suppose that cocones $\bar\mu, \bar\mu' \from \Eq^*(f,g) \tto Z$ satisfy $R_Z\bar\mu = R_Z\bar\mu'$. Then $\bar\mu^\kappa\bar\zeta_\sigma^\kappa = \bar\mu'^\kappa\bar\zeta_\sigma^\kappa$ for all $\sigma,\kappa \geq 0$. Since $\zeta^\kappa_\sigma$ is surjective, $\bar \zeta^\kappa_\sigma$ is an epimorphism, hence $\bar\mu = \bar\mu'$.
    \begin{center}
    \quad\quad
    \begin{tikzcd}[column sep=large,row sep=normal]
    \Rips_\sigma\Eq^\kappa(f,g) \arrow[r,"\bar h^{\kappa'\kappa}_\sigma"] \arrow[d,"\bar\zeta_\sigma^\kappa" swap] \arrow[dd,"\bar\nu_\sigma^\kappa" swap, bend right =60]  & \Rips_\sigma\Eq^{\kappa'}(f,g) \arrow[d,"\bar\zeta_\sigma^{\kappa'}"] \arrow[dd,"\bar\nu_\sigma^{\kappa'}", bend left =60]\\ \Eq^\kappa(f,g) \arrow[r,"\bar\iota^{\kappa'\kappa}" swap] \arrow[d,"\bar\mu^\kappa" swap]   & \Eq^{\kappa'}(f,g) \arrow[d,"\bar\mu^{\kappa'}"] \\ Z \arrow[r, equal] & Z \end{tikzcd}
    \end{center}
    To verify surjectivity of $R_Z$, consider a cocone $\bar\nu \from \Rips_* \Eq^*(f,g) \tto Z$. This induces a cocone $\bar\nu^\kappa \from \Rips_* \Eq^\kappa(f,g) \tto Z$ for each $\kappa \geq 0$. By Proposition \ref{prop:colimit_coarse},  $\bar\zeta^\kappa \from \Rips_*\Eq^\kappa(f,g) \tto \Eq^\kappa(f,g)$ is a limit cocone in $\barBorn$ and so there exists a unique morphism $\bar \mu^\kappa \from \Eq^\kappa(f,g) \to Z$ in $\barBorn$ such that $\bar\nu^\kappa = \bar\mu^\kappa\bar\zeta^\kappa$.    For $\sigma \geq 0$ and $\kappa \leq \kappa'$, let $h^{\kappa'\kappa}_\sigma \from \Rips_\sigma\Eq^\kappa(f,g) \rightarrow \Rips_\sigma \Eq^{\kappa'}(f,g)$ be the 1--Lipschitz map coinciding with the underlying inclusion $\iota^{\kappa'\kappa}$. Then
    \[\bar\mu^\kappa\bar\zeta^\kappa_\sigma = \bar\nu_\sigma^\kappa = \bar\nu_\sigma^{\kappa'}\bar h^{\kappa'\kappa} = \bar\mu^{\kappa'}\bar\zeta_\sigma^{\kappa'}\bar h^{\kappa'\kappa} = \bar\mu^{\kappa'}\bar\iota^{\kappa'\kappa}\bar\zeta_\sigma^\kappa\]
    for all $\kappa \leq \kappa'$, hence $\bar\mu^\kappa = \bar\mu^{\kappa'}\bar\iota^{\kappa'\kappa}$ since $\bar\zeta_\sigma^\kappa$ is an epimorphism. Therefore, we obtain a cocone $\bar\mu \from \Eq^*(f,g) \tto Z$ with colegs $\bar\mu^\kappa$. Since $\bar\nu^\kappa_\sigma = \bar\mu^\kappa\bar\zeta^\kappa_\sigma$ for all $\sigma, \kappa\geq 0$, it follows that $R_Z\bar\mu = \bar\nu$.
    \endproof

    \begin{corollary}\label{cor:RE}
     There exists a canonical isomorphism
     \[\colimit \Rips_* \Eq^*(f,g) \cong \colimit \Eq^*(f,g)\]
     in $\barBorn$, where each colimit exists if and only if $\Eq^*(f,g)$ stabilises. Moreover, either colimit, if it exists, is realised by $\Eq^\kappa(f,g)$ for $\kappa \geq 0$ sufficiently large. \qed
     \end{corollary}

    \begin{proposition}[Universal $\cat$--geodesic cones for parallel pairs]\label{prop:geod-eq}
       Let ${\bar f, \bar g \from X \rightrightarrows Y}$ be a parallel pair in $\barBorn$, with representatives $f\in \bar f$, $g \in \bar g$. Then the following are equivalent:
     \begin{enumerate}
       \item $\bar f, \bar g$ admits a universal $\cat$--geodesic cone;
       \item $\bar f, \bar g$ admits an equaliser in $\barBorn$ which is coarsely geodesic;
       \item $\Eq^*(f,g)$ stabilises and $\Eq^\kappa(f,g)$ is coarsely geodesic for $\kappa \geq 0$ large;
       \item The universal $\cat$--geodesic cone for $\bar f, \bar g$ is realised by $\bar\iota^{\kappa}\bar\zeta_\sigma \from \Rips_\sigma \Eq^\kappa(f,g) \to X$ for some $\sigma, \kappa \geq 0$;
       \item $\colimit \Rips_* \Eq^*(f,g)$ in $\barcat$ is realised by $\Rips_\sigma \Eq^\kappa(f,g)$ for some $\sigma,\kappa \geq 0$.
      \end{enumerate}
      In addition, if $\cat$ has dominated controls then each item above is equivalent to:
      \begin{enumerate}
       \item[6.] $\colimit \Rips_* \Eq^*(f,g)$ in $\barcat$ exists.
      \end{enumerate}
    \end{proposition}

    \proof
    The implications (1) $\iff$ (2) $\iff$ (3) $\implies$ (4) follow from Theorem \ref{thm:cones} and Proposition \ref{prop:equaliser}; while (4) $\implies$ (1) and (5) $\implies$ (6) are immediate.

    (3) $\implies$ (5).
    Choose $\kappa \geq 0$ sufficiently large so that $\Eq^\kappa(f,g)$ is coarsely geodesic and $\bar\iota^{\kappa'\kappa}$ is an isomorphism for all $\kappa' \geq \kappa$. By Corollary \ref{cor:RE}, there exist canonical isomorphisms
    \[\Eq^{\kappa}(f,g) \cong \colimit \Eq^*(f,g) \cong \colimit \Rips_* \Eq^*(f,g)\]
    in $\barBorn$. Let $\bar\lambda \from \Rips_* \Eq^*(f,g) \tto \Eq^\kappa(f,g)$ be the limit cocone in $\barBorn$.
    \begin{center}
    \begin{tikzcd}[column sep=huge,row sep=huge] \Rips_* \Eq^*(f,g) \arrow[Rightarrow, drr,"\bar\mu"] \arrow[Rightarrow, dr,"\bar\lambda'" swap] \arrow[Rightarrow, d,"\bar\lambda" swap]   &  &  \\ \Eq^\kappa(f,g) \arrow[r,"(\bar\zeta^\kappa_\sigma)^{-1}", swap]  & \Rips_\sigma\Eq^\kappa(f,g) \arrow[r,"\bar m", swap, dashed] & Y \end{tikzcd}
    \end{center}
    Since $\Eq^\kappa(f,g)$ is coarsely geodesic, by Proposition \ref{prop:C-geodesic}, $\bar\zeta^\kappa_\sigma \from \Rips_\sigma\Eq^\kappa(f,g) \to \Eq^\kappa(f,g)$ is an isomorphism in $\barBorn$ for $\sigma \geq 0$ large. Then
    \[\bar\lambda' := (\bar\zeta_\sigma^\kappa)^{-1} \bar\lambda \from \Rips_* \Eq^*(f,g) \tto \Rips_\sigma\Eq^\kappa(f,g)\]
    is also a limit cocone in $\barBorn$. Therefore, for any cocone $\bar\mu \from \Rips_* \Eq^*(f,g) \tto Y$, there is a unique morphism $\bar m \from \Rips_\sigma\Eq^\kappa(f,g) \to Y$ in $\barBorn$ satisfying $\bar\mu = \bar m \bar\lambda'$. Since $\Rips_\sigma\Eq^\kappa(f,g)$ is a graph, $\bar m$ is a morphism in $\barcat$. Thus, $\bar\lambda'$ is a limit cocone in $\barcat$.

    (5) $\implies$ (3). By Lemma \ref{lem:enlarge}, $\Rips_\sigma\Eq^\kappa(f,g)$ realises $\colimit \Rips_* \Eq^*(f,g)$ in $\barBorn$. Therefore, by Corollary \ref{cor:RE}, $\colimit \Eq^*(f,g) \cong \colimit \Rips_* \Eq^*(f,g)$ exists in $\barBorn$ and is coarsely geodesic. Consequently, by Proposition \ref{prop:equaliser}, $\Eq^*(f,g)$ stabilises and $\Eq^{\kappa'}(f,g) \cong \colimit \Eq^*(f,g)$ is coarsely geodesic for $\kappa' \geq 0$ large.

    \emph{Now assume that $\cat$ has dominated controls.}

    (6) $\implies$ (3).
    Suppose that $W \cong \colimit \Rips_* \Eq^*(f,g)$ in $\barcat$. Then $W$ is a colimit of graphs in $\barcat$ and so it must be $\cat$--geodesic by Theorem \ref{thm:colimits2}. The remainder of the proof is identical to that of {(5) $\implies$ (3)} with $W$ in place of of $\Rips_\sigma\Eq^\kappa(f,g)$.
    \endproof

    \begin{lemma}\label{lem:eq-threshold}
     Assume item (4) from Proposition \ref{prop:geod-eq} holds for constants $\sigma, \kappa$. Then $\Eq^*(f,g)$ stabilises at threshold $\kappa$.
    \end{lemma}

    \proof
    By assumption, $\iota^{\kappa}\zeta_\sigma \from \Rips_\sigma\Eq^{\kappa'}(f,g) \to X$ realises the universal $\cat$--geodesic cone over $\bar f, \bar g$. Let $\kappa' \geq \kappa$ and consider the map $\iota^{\kappa'}\zeta_\sigma \from \Rips_\sigma\Eq^{\kappa'}(f,g) \to X$. This is a $\cat$--geodesic cone over $\bar f, \bar g$, and so it uniquely factors through $\iota^{\kappa}\zeta_\sigma$ up to closeness. Therefore the image of $\iota^{\kappa'}\zeta_\sigma$ lies in a bounded neighbourhood of the image of $\iota^{\kappa}\zeta_\sigma$ in $X$. Consequently, the inclusion $\Eq^\kappa(f,g) \hookrightarrow \Eq^{\kappa'}(f,g)$ is coarsely surjective. The result follows using Proposition \ref{prop:equaliser}.
    \endproof

   \subsection{Uniformly controlled cones}\label{sec:uc-cone}

   We now wish to compute universal $\cat$--geodesic cones over diagrams $D \from \sJ \to \barBorn$ of arbitrary shape. By Theorem \ref{thm:cones}, the first step is to compute $\limit_\sJ D$, if it exist. For a finite family of spaces, their setwise product (with the $l^\infty$--metric) realises their categorical product in $\barBorn$. Therefore, for finite diagrams, we can follow a standard categorical argument to realise $\limit_\sJ D$ as an equaliser involving appropriate finite products, then apply Proposition \ref{prop:equaliser}. However, as $\barBorn$ does not admit arbitrary (categorical) products, it does not appear obvious how one can compute $\limit_\sJ D$ in general.

   In order to handle diagrams of arbitrary shape, we introduce the notion of \emph{uniformly controlled diagrams}. These are not (in general) true diagrams valued in either $\barBorn$ or $\Born$, but they do provide workable middle ground for situations requiring uniform control on maps and coarse commutivity.

   For this setup, we do not actually need the diagram to be indexed by a small category -- it suffices to assume that $\sJ$ is a directed multigraph. That is, we do not require compositions to be defined, nor the existence of identity arrows at each object. Nevertheless, we shall keep the terminology and notation mostly consistent with that of true diagrams in order to make the analogy more apparent. Thus, we shall refer to the vertices of $\sJ$ as objects, and its directed edges as arrows.

   \begin{definition}[Uniformly controlled diagram]
   Let $\cat = \cat(\F)$.
   A \emph{uniformly controlled (\uc) diagram} $D \from \sJ \sqto \cat$ in $\cat$ of shape $\sJ$ comprises the following data:
  \begin{itemize}
   \item a metric space $D_j$ for each $j\in\sJ$, and
   \item a map $\beta_\phi \from D_i \to D_j$ for each arrow $\phi \from i \to j$ in $\sJ$
  \end{itemize}
  such that all $\beta_\phi$ admit a common upper control from $\F$.
  \end{definition}

  \begin{definition}[U.C. cone]
   A \emph{\uc cone} $\mu \from X \sqtto D$ over a \uc diagram ${D \from \sJ \sqto \cat}$ with \emph{apex} $X$ comprises a family of maps ${\mu_j \from X \to D_j}$ in $\cat$ for $j\in\Obj(\sJ)$, called the \emph{legs}, such that
     \begin{itemize}
   \item all $\mu_j$ admit a common upper control from $\F$, and
   \item $\mu_j$ and $\beta_\phi \circ \mu_i$ are uniformly close for every arrow $\phi \from i \to j$ in $\sJ$.
  \end{itemize}
  \end{definition}

    Suppose that $D, D' \from \sJ \sqto \cat$ are \uc diagrams with respective spaces $D_j$, $D'_j$ for $j \in \sJ$ and maps $\beta_\phi$, $\beta'_\phi$ for $\phi \in \Mor(\sJ)$. We say that $D, D' \from \sJ \sqto \cat$ are $\kappa$--close if $D_j = D'_j$ for all $j \in \Obj(\sJ)$ and $\beta_\phi \approx_\kappa \beta'_\phi$ for all $\phi\in\Mor(\sJ)$. Similarly, two \uc cones $\mu,\mu' \from X \sqtto D$ are $\kappa$--close if $\mu_j \approx_\kappa \mu'_j$ for all respective legs. We write $D \approx D'$ or $\mu \approx \mu'$ if they are $\kappa$--close for some $\kappa \geq 0$.

    Given a \uc cone $\mu \from X \sqtto D$ and a controlled map $f \from W \to X$, we may define a \uc cone $\mu f \from W \sqtto D$ whose legs are $\mu_j f$. Moreover, if $f \approx f'$ and $\mu \approx \mu'$ then $\mu f \approx \mu' f'$.

    \begin{definition}[U.C. limit]
    Let ${D \from \sJ \sqto \cat}$ be a \uc diagram. A \uc cone $\lambda \from X \sqtto D$ over $D$ is a \emph{\uc limit cone} if for any \uc cone $\mu \from W \sqtto D$, there exists a unique morphism $\bar f \from W \to X$ in $\barcat$ such that $\mu \approx \lambda f$ for any $f \in \bar f$. The \emph{\uc limit} $\uclim_\sJ D$ of $D$ is the apex of any \uc limit cone over $D$.
    \end{definition}

    Any \uc limit cone $\lambda \from X \sqtto D$, if it exists, is unique in the following sense: if $\lambda' \from X' \sqtto D$ is also a \uc limit cone then there is a unique isomorphism $\bar f \from X' \to X$ in $\barcat$ such that $\lambda' \approx \lambda f$ for all $f \in \bar f$.

   Given a \uc diagram ${D \from \sJ \sqto \cat}$, consider the product space $\prod_{j} D_j$ with the $\ell^\infty$--metric. For a metric space $W$ equipped with a family of functions $h_j \from W \to D_j$ indexed by $j\in\Obj(\sJ)$,  the \emph{product map} $\prod_{j} h_j \from W \to \prod_{j} D_j$ is given by $w \mapsto (h_jw)_j$. Observe that $\prod_{j} h_j$ has upper control $\rho$ if and only if the family $h_j$ has uniform upper control $\rho$. Furthermore, another family $h'_j \from W \to D_j$ satisfies $h_j \approx_\kappa h'_j$ for all $j\in\Obj(\sJ)$ if and only if $\prod_{j} h'_j \approx_\kappa \prod_{j} h_j$. In other words, it is safe to work up to closeness on the level of products  once we have chosen a family of representatives $h_j$ (up to uniform closeness). If we instead take the closeness classes $\bar h_j$ first, then their product class may not be well-defined, as the following example shows.

    \begin{example}
    Consider the family of maps $h_k \from \N \to \N$, for $k\in \N$, given by ${h_k(n) := n + k}$. Then $\bar h_k = \bar 1_\N$ for all $k\in\N$, but the product maps $\prod_k 1_\N$ and $\prod_k h_k$ are not close.
   \end{example}

   Let us verify that $\ell^\infty$--products realise \uc analogues of categorical products. Call a \uc diagram $D \from \sJ \sqto \cat$ \emph{discrete} if all maps $\beta_\phi$ are identity maps. For a discrete \uc diagram $D$, the standard projection cone $\pi \from \prod_{j} D_j \tto D$ in $\cat$ whose legs $\pi_j$ are the 1--Lipschitz projections to each factor $D_j$ is itself a \uc cone.

   \begin{lemma}[\uc products exist in $\cat$]
    Let $\sJ$ be a discrete category and $D \from \sJ \sqto \cat$ be a \uc diagram. Then the standard projection cone $\pi \from \prod_{j} Dj \tto D$ is a \uc limit cone.
    \end{lemma}

    \proof
    Let $\mu \from X \sqtto D$ be a \uc cone. Then $\mu = \pi \circ \prod_j \mu_j$. Since the legs $\mu_j \from X \to D_j$ admit a uniform upper control $\rho \in \F$, the product map $\prod_j \mu_j$ also has upper control $\rho$ and is hence a morphism in $\cat$. To verify uniqueness, suppose $h \from X \to \prod_{j} D_j$ is a morphism in $\cat$ satisfying $\mu \approx_\kappa \pi h$ for some $\kappa$. Then for all $x \in X$, we have that $\pi_j h x \approx_\kappa \mu_j x$ for each $j \in \Obj(\sJ)$, hence $h x \approx_\kappa (\mu_jx)_j = (\prod_j \mu_j) x$. It follows that $h \approx_\kappa \prod_j \mu_j$.
    \endproof

   To handle \uc diagrams of arbitrary shape, we follow the standard construction of limits as equalisers in the category $\Set$ (see, for example, \cite[Theorem 3.2.13]{Rie16}), while keeping track of controls and closeness constants. Given a \uc diagram $D \from \sJ \sqto \cat$, consider a parallel pair
   \begin{center}
    \begin{tikzcd}
   \displaystyle\prod_{j \in \Obj(\sJ)} Dj \arrow[r, shift left, "\gamma"] \arrow[r, shift right, "\delta" below] & \displaystyle\prod_{\phi \in \Mor(\sJ)} D(\cod \phi)
  \end{tikzcd}
  \end{center}
  where $\gamma$ and $\delta$ are given by $(x_j)_j \mapsto (x_{\cod \phi})_\phi$ and $(x_j)_j \mapsto (\beta_\phi x_{\dom \phi})_\phi$ respectively; we shall refer to these as the \emph{codomain} and \emph{domain} maps associated to $D$ (and the choice of representative maps). Observe that $\gamma$ is $1$--Lipschitz, while $\delta$ admits an upper control $\rho$ if and only if $\rho$ serves as a uniform upper control for all $\beta_\phi$. If $D' \from \sJ \sqto \cat$ is a \uc diagram whose objects coincide with $D$ then the associated codomain map $\gamma'$ coincides with $\gamma$; in addition, the associated domain map $\delta'$ satisfies $\delta \approx_\kappa \delta'$ if and only if $D \approx_\kappa D'$.

   \begin{remark}
    If $\sJ$ is finite, then all uniform requirements hold automatically, so long as all involved maps are controlled. Moreover, any choice of representatives $\beta_\phi$ will be unique up to uniform closeness. Thus, $\bar\delta$ depends only on the closeness classes $\bar\beta_\phi$.
   \end{remark}

    There is a natural correspondence between closeness classes of \uc cones over $D$ and morphisms in $\barcat$ equalising $\bar\gamma$ and $\bar\delta$; this follows immediately from the definition of $\gamma$ and $\delta$.

   \begin{lemma}[\uc cones equalise]\label{lem:uc-factor}
     Let $\mu_j \from X \to D_j$ be a family of maps from a metric space $X$ indexed by $j \in \Obj(\sJ)$. Let $\kappa \geq 0$. Then $\gamma \circ \prod_{j} h_j \approx_\kappa \delta \circ \prod_{j} \mu_j$ if and only if $\mu_{\cod\phi} \approx_\kappa \beta_\phi \mu_{\dom\phi}$ for every arrow $\phi\in\Mor(\sJ)$. Consequently, the $\mu_j$ form the legs of a \uc cone $\mu \from X \sqtto D$ in $\cat$ if and only if they admit a common upper control $\rho \in \F$ and satisfy $\bar\gamma \circ \overline{\prod_{j} \mu_j} = \bar\delta \circ \overline{\prod_{j} \mu_j}$.
     \qed
   \end{lemma}

   \begin{corollary}[\uc limit cones as equalisers]\label{cor:uc-equaliser}
    A \uc cone $\lambda \from X \sqtto D$ is a \uc limit cone in $\cat$ if and only if $\prod_{j} \lambda_j$ realises the equaliser of $\bar\gamma, \bar\delta$ in $\barcat$. \qed
   \end{corollary}

   Thus, if $\uclim_\sJ D$ exists, we obtain an equaliser diagram in $\barcat$; moreover any \uc cone $\mu \from X \sqtto D$ factors uniquely through $\uclim_\sJ D$.

   \begin{center}
    \begin{tikzcd}
   X \arrow[d, dashed] \arrow[dr, "\overline{\prod_{j} \mu_j}"] & & \\
   \uclim_\sJ D \arrow[r] & \displaystyle\prod_{j \in \Obj(\sJ)} Dj \arrow[r, shift left, "\bar\gamma"] \arrow[r, shift right, "\bar\delta" below] & \displaystyle\prod_{\phi \in \Mor(\sJ)} D(\cod \phi)
  \end{tikzcd}
  \end{center}

   The problem of finding a universal uniformly controlled cone (up to uniform closeness) on $D$ reduces to computing the equaliser of $\bar\gamma, \bar\delta$ in $\barcat$. Explicitly, the $\kappa$--equaliser of $\gamma, \delta$ is
   \[\Eq^\kappa (\gamma, \delta) = \left\{(x_j)_j \in \prod_{j \in \Obj(\sJ)} D_j \; : \; x_{\cod\phi} \approx_\kappa \beta_\phi x_{\dom\phi} \textrm{ for all } \phi \right\}.\]
   This is the space of tuples in $\prod_{j \in \Obj(\sJ)} D_j$ which satisfy each constraint up to uniform error $\kappa$. Stability of this equaliser filtration means that, up to closeness, all coarse solutions must appear by some finite threshold $\kappa$. Write $\iota^\kappa \from \Eq^\kappa(\gamma, \delta) \hookrightarrow \prod_j D_j$ for the inclusion. Invoking Proposition \ref{prop:equaliser}, we immediately obtain the following.

   \begin{proposition}[Characterisation of \uc limits]\label{prop:tupstab}
     Let $D \from \sJ \sqto \cat$ be a \uc diagram. The following are equivalent:
     \begin{enumerate}
      \item $\uclim_\sJ D$ exists,
      \item $\bar\gamma, \bar\delta$ admits an equaliser,
      \item $\Eq^*(\gamma, \delta)$ stabilises, and
      \item for all $\kappa \geq 0$ sufficiently large, $\pi \iota^\kappa \from \Eq^\kappa(\gamma, \delta) \sqtto D$ is a \uc limit cone.
     \end{enumerate}
      In particular, $\uclim_\sJ D \cong \colimit\Eq^*(\gamma, \delta)$ in $\barcat$, where each (co)limit exists if and only if the other does. \qed
   \end{proposition}

   \begin{remark}
   Stability of $\Eq^*(\gamma, \delta)$ does not depend on choice of ambient category. So there is no loss of generality in assuming that $\cat = \Born$ when dealing with \uc limits. In other words, we may write $\uclim_\sJ D$ without referring to the ambient category.
   \end{remark}

   We now characterise universal $\cat$--geodesic \uc cones over $D$: a $\cat$--geodesic cone $\lambda \from X \sqtto D$ such that for any $\cat$--geodesic cone $\mu \from W \sqtto D$, there exists a unique morphism $\bar f \from W \to X$ in $\barcat$ such that $\mu \approx \lambda f$ for any representative $f$. In this situation, the choice of category $\cat$ does come into play. Write $\zeta^\kappa \from \Rips_*\Eq^\kappa(\gamma, \delta) \tto \Eq^\kappa(\gamma, \delta)$ for the canonical cocone of $\Eq^\kappa(\gamma, \delta)$.

   \begin{theorem}[Characterisations of universal $\cat$--geodesic \uc cones]\label{thm:univ-geod-cone}
   Let $D \from \sJ \sqto \cat$ be a \uc diagram. The following are equivalent:
    \begin{enumerate}
     \item $D$ admits a universal $\cat$--geodesic \uc cone,
     \item $\uclim_\sJ D$ exists and is coarsely geodesic,
     \item the stable Rips colimit $\floorC{\uclim_\sJ D}$ exists,
     \item the equaliser of $\bar\gamma, \bar\delta$ exists and is coarsely geodesic,
     \item $\Eq^*(\gamma, \delta)$ stabilises to a coarsely geodesic space, and
     \item $\pi \iota^\kappa \zeta_\sigma^\kappa \from \Rips_\sigma \Eq^\kappa(\gamma, \delta) \sqtto D$ is a universal $\cat$--geodesic cone for some $\sigma, \kappa \geq 0$.
    \end{enumerate}
    In addition, if $\cat$ has dominated controls then each item above is equivalent to:
    \begin{enumerate}
       \item[7.] $\colimit \Rips_* \Eq^*(\gamma,\delta)$ in $\barcat$ exists.
    \end{enumerate}
    \end{theorem}

   \proof
   This follows from Proposition \ref{prop:geod-eq}, Lemma \ref{lem:uc-factor}, and Proposition \ref{prop:tupstab}.
   \endproof

  This yields a proof of Theorem \ref{thm:recipe}.

  If $D$ admits a universal $\cat$--geodesic \uc cone then we obtain the following diagram, where $\mu \from X \sqtto D$ is any $\cat$--geodesic \uc cone.

    \begin{center}
    \begin{tikzcd}[row sep = large]
    X \arrow[d, dashed] \arrow[drr, "\overline{\prod_{j} \mu_j}"] & & & \\
   \left\lfloor\uclim_\sJ D\right\rfloor_{\cat} \arrow[r] & \uclim_\sJ D \arrow[r] & \displaystyle\prod_{j \in \Obj(\sJ)} Dj \arrow[r, shift left, "\bar\gamma"] \arrow[r, shift right, "\bar\delta" below] & \displaystyle\prod_{\phi \in \Mor(\sJ)} D(\cod \phi)
  \end{tikzcd}
  \end{center}

   The following is a useful criterion for checking that a given $\cat$--geodesic space realises a univeral $\cat$--geodesic cone.
  
  \begin{corollary}\label{cor:uc-cat-cone}
    Let $X$ be a $\cat$--geodesic space which realises the equaliser of $\bar\gamma, \bar\delta$ in $\barBorn$. Then $X$ realises the universal $\cat$--geodesic \uc cone above $D$.
    Consequently, for some $\sigma, \kappa \geq 0$, there exists a canonical isomorphism $X \cong \Rips_\sigma\Eq^\kappa(\gamma,\delta)$ in $\barcat$  coinciding (up to closeness) with the underlying product map $\prod_j \lambda_j \from X \to \prod_j D_j$.
   \end{corollary}

   \proof
   This follows from Corollary \ref{cor:uc-cone-check}, Theorem \ref{thm:univ-geod-cone}, and the universal property.
   \endproof

   \begin{corollary}
   Assume $\barcat$ has dominated controls. If $D$ admits a universal $\cat$--geodesic \uc cone then there are canonical isomorphisms
   \[\left\lfloor\uclim_\sJ D\right\rfloor_{\cat} \cong \left\lfloor\colimit\Eq^*(\gamma, \delta)\right\rfloor_{\cat} \cong \colimit\Rips_*\Eq^*(\gamma, \delta)\]
   in $\barcat$, and each displayed (co)limit exists and realises the universal $\cat$--geodesic \uc cone over $D$; otherwise, none of the displayed (co)limit exist.
   \end{corollary}

  \proof
  This follows from Proposition \ref{prop:geod-eq}, Proposition \ref{prop:tupstab}, and Theorem \ref{thm:univ-geod-cone}.
  \endproof

  \subsection{Interactions with retractions}

  The goal of this subsection is to prove that the existence of universal $\cat$--geodesic \uc cones persists under suitable retractions of the associated \uc diagrams. Our strategy is to show that existence of equalisers, existence of \uc limits, and coarse geodesicity are preserved under appropriate retractions of the data.

  \begin{lemma}[Retractions and factoring]\label{lem:nat-equal}
  Let $f, g \from X \rightrightarrows Y$ and $f', g' \from X \rightrightarrows Y$ be maps between metric spaces.   Let $\alpha_X \from X \to X'$ and $\alpha_Y \from Y \to Y'$ be maps such that $f'\alpha_X \approx_K \alpha_Y f$ and $g'\alpha_X \approx_K \alpha_Y g$ for some $K \geq 0$. Assume $\alpha_X, \alpha_Y$ have upper control $\rho\in\F$. Then for all $\kappa \geq 0$, there is a unique map $a \from \Eq^\kappa(f,g) \to \Eq^{\kappa'}(f',g')$ through which the composite $\Eq^\kappa(f,g) \hookrightarrow X \xrightarrow{\alpha_X} X'$ factors,   where $\kappa' = 2K + \rho\kappa$. Moreover, $a$   has upper control $\rho$.
  \end{lemma}

  \begin{center}
   \begin{tikzcd}[column sep=huge, row sep=large] \Eq^\kappa(f,g) \arrow[r, hookrightarrow] \arrow[d, dashed, "a" swap] & X \arrow[r, shift left, "f"] \arrow[r, shift right, "g"'] \arrow[d, "\alpha_X" swap]  & Y \arrow[d, "\alpha_Y"] \\ \Eq^{\kappa'}(f',g')  \arrow[r, hookrightarrow] & X'  \arrow[r, shift left, "f'"] \arrow[r, shift right, "g'"'] & Y'   \end{tikzcd}
   \end{center}

  \proof
  Let $x \in \Eq^\kappa(f,g) \subseteq X$. Then $fx \approx_\kappa gx$, hence
  \[f'\alpha_X x \approx_K \alpha_Y fx \approx_{\rho\kappa} \alpha_Y gx \approx_K g'\alpha_X x.\]
  Therefore, $\alpha_X x \in \Eq^{\kappa'}(f',g')$, where $\kappa' = 2K + \rho\kappa$. The induced map $a \from \Eq^\kappa(f,g) \to \Eq^{\kappa'}(f',g')$ is a (co)restriction of $\alpha_X$, and hence has upper control $\rho$. Uniqueness is immediate.
  \endproof

  \begin{lemma}[Retractions and stability]\label{lem:ret-stab}
   In addition to the assumptions given in Lemma \ref{lem:nat-equal}, let $\omega_X \from X' \to X$ and $\omega_Y \from Y' \to Y$ be controlled maps satisfying $f\omega_X \approx \omega_Y f'$, $g\omega_X \approx \omega_Y g'$, $\alpha_X\omega_X \approx 1_{X'}$, and $\alpha_Y\omega_Y \approx 1_{Y'}$.    Assume that $\Eq^*(f,g)$ stabilises at threshold $\kappa_0$. Then $\Eq^*(f',g')$ stabilises at threshold $\kappa'_0 = 2K + \rho\kappa_0$.
   \end{lemma}

  \proof
  Let $t' \geq 2K + \rho\kappa_0$ and $x' \in \Eq^{t'}(f',g')$. Since $\omega_Y$ is controlled, we may apply Lemma \ref{lem:nat-equal} to deduce $\omega_X x' \in \Eq^{t}(f,g)$ for some $t \geq 0$. We may also assume $t \geq \kappa_0$. By the stability assumption, there exists $r \geq 0$ such that $\omega_X x' \approx_r x$ for some $x \in \Eq^{\kappa_0}(f,g)$.

    \begin{center}
   \begin{tikzcd}[column sep=huge, row sep=large] \Eq^{\kappa_0}(f,g) \arrow[r, hookrightarrow] \arrow[d, dashed] & \Eq^{t}(f,g) \arrow[r, hookrightarrow]  & X \arrow[r, shift left, "f"] \arrow[r, shift right, "g"'] \arrow[d, "\alpha_X" swap, bend right=20]  & Y \arrow[d, "\alpha_Y" swap, bend right=20] \\ \Eq^{\kappa'_0}(f',g') \arrow[r, hookrightarrow] & \Eq^{t'}(f',g')  \arrow[r, hookrightarrow] \arrow[u, dashed] & X'  \arrow[r, shift left, "f'"] \arrow[r, shift right, "g'"'] \arrow[u, "\omega_X" swap, bend right=20]  & Y' \arrow[u, "\omega_Y" swap, bend right=20]  \end{tikzcd}
   \end{center}

  Therefore, by Lemma \ref{lem:nat-equal},
  \[\alpha_X\omega_X x' \approx_{\rho r} \alpha_X x \in \Eq^{\kappa'_0}(f',g')\]
  where $\kappa'_0 = 2K + \rho\kappa_0$. Since $\alpha_X\omega_X \approx 1_{X'}$, it follows that $x'$ lies in a bounded neighbourhood of $\Eq^{\kappa'_0}(f',g')$. Consequently, $\Eq^{\kappa'_0}(f',g') \hookrightarrow \Eq^{t'}(f',g')$ is coarsely surjective.
  \endproof

  \begin{corollary}[Retractions and equalisers]\label{cor:eq_ret}
   Let $\bar f, \bar g \from X \rightrightarrows Y$ and $\bar f', \bar g' \from X \rightrightarrows Y$ be parallel pairs in $\barcat$. Assume that $\bar\alpha_X \from X \to X'$, $\bar\alpha_Y \from Y \to Y'$, $\bar\omega_X \from X' \to X$, $\bar\omega_Y \from Y' \to Y$ are morphisms in $\barcat$ satisfying
   \[\bar f'\bar\alpha_X = \bar\alpha_Y \bar f, \quad \bar g'\bar\alpha_X = \bar\alpha_Y \bar g, \quad \bar f\bar\omega_X = \bar\omega_Y \bar f', \quad \bar g\bar\omega_X = \bar\omega_Y \bar g', \quad \bar\alpha_X\bar\omega_X = \bar 1_{X'}, \quad \bar\alpha_Y\bar\omega_Y = \bar 1_{Y'}.\]
   Then if $\bar f, \bar g$ admits an equaliser, so does $\bar f', \bar g'$. Moreover, $\bar\alpha_X$ induces a canonical right-invertible morphism in $\barcat$ between the respective limits. \qed
  \end{corollary}

   Next, we introduce a uniformly controlled analogue of a natural transformation.

   \begin{definition}[\uc natural transformation]
   Let $D, D' \from \sJ \sqto \Born$ be \uc diagrams, with morphisms denoted respectively by $\beta_\phi, \beta_\phi'$ maps for each arrow $\phi$. A \emph{\uc natural transformation} $\alpha \from D \sqtto D'$ in $\cat$ comprises a family of maps $\alpha_j \from D_j \to D'_j$ for $j \in \Obj(\sJ)$ admitting a common upper control from $\F$ such that the diagram       \begin{center}
   \begin{tikzcd}[column sep=large, row sep=normal] D_i \arrow[r, "\beta_\phi"] \arrow[d, "\alpha_i" swap]  & D_j \arrow[d, "\alpha_j"] \\ D'_i \arrow[r, "\beta'_\phi" swap] & D'_j   \end{tikzcd}
  \end{center}
  uniformly coarsely commutes for all arrows $\phi \from i \to j$ in $\sJ$. Composition of \uc natural transformations is defined componentwise.
  \end{definition}

  \begin{remark}
   A \uc cone $\mu \from X \sqtto D$ is a special case of a \uc natural transformation, where we interpret $X \from \sJ \tto \Born$ as a constant functor (all objects are $X$, all arrows $1_X$). Dualising this, we may define s \emph{\uc cocone} $\nu \from D \sqtto X$ as \uc natural transformation to a constant functor.
  \end{remark}

  By simultaneously considering all components $\alpha_j$ associated to a \uc natural transformation $\alpha \from D \sqtto D'$, we obtain a diagram of product spaces and product maps
  \begin{center}
   \begin{tikzcd}[column sep=huge, row sep=large] \prod_j D_j \arrow[r, shift left, "\gamma"] \arrow[r, shift right, "\delta"'] \arrow[d, "\prod_j \alpha_j" swap]  & \prod_\phi D_{\cod \phi} \arrow[d, "\prod_\phi \alpha_{\cod\phi}"] \\ \prod_j D'_j   \arrow[r, shift left, "\gamma'"] \arrow[r, shift right, "\delta'"'] & \prod_\phi D'_{\cod \phi}   \end{tikzcd}
   \end{center}
   where $\gamma, \delta$ (resp.~ $\gamma', \delta'$) are the codomain and domain maps for $D$ (resp.~$D'$). We make the following observations:
   \begin{itemize}
    \item any uniform upper control $\rho$ for the $\alpha_j$ serves as an upper control for $\prod_j \alpha_j$ and $\prod_\phi \alpha_{\cod \phi}$,
    \item $\gamma' \circ \prod_j \alpha_j = (\prod_\phi \alpha_{\cod \phi}) \circ \gamma$, and
    \item $\delta' \circ \prod_j \alpha_j \approx_K (\prod_\phi \alpha_{\cod \phi}) \circ \delta$ if and only if $\beta'_{\phi} \alpha_{i} \approx_K \alpha_{j} \beta_{\phi}$ for all arrows $\phi \from i \to j$.
   \end{itemize}

   Consequently, if we work up to closeness, any \uc natural transformation induces a true natural transformation in $\barBorn$.

   \begin{lemma}[\uc to Natural transformation]\label{lem:nat-trans}
    Let $\alpha \from D \sqtto D'$ be a \uc natural transformation. Then $\alpha$ induces a natural transformation
    \begin{center}
   \begin{tikzcd}[column sep=huge, row sep=large] \prod_j D_j \arrow[r, shift left, "\bar\gamma"] \arrow[r, shift right, "\bar\delta"'] \arrow[d, "\overline{\prod_j \alpha_j}" swap]  & \prod_\phi D_{\cod \phi} \arrow[d, "\overline{\prod_\phi \alpha_{\cod\phi}}"] \\ \prod_j D'_j   \arrow[r, shift left, "\bar\gamma'"] \arrow[r, shift right, "\bar\delta'"'] & \prod_\phi D'_{\cod \phi}   \end{tikzcd}
   \end{center}
   in $\barBorn$ whose components (vertical arrows) are morphisms in $\barcat$. \qed   \end{lemma}

   Appealing to Corollary \ref{cor:uc-equaliser}, any \uc natural transformation $\alpha \from D \sqtto D'$ in $\cat$ induces a canonical morphism $\uclim_{\sJ} \alpha \from \uclim_{\sJ} D \to \uclim_{\sJ} D'$ in $\barcat$ whenever the displayed \uc limits exist; moreover this behaves functorially where defined.

   Say that two \uc natural transformations $\alpha, \alpha' \from D \sqtto D'$ are \emph{$\kappa$--close}, denoted as $\alpha \approx_\kappa \alpha'$, if their components satisfy $\alpha_j \approx_\kappa \alpha'_j$ for all $j$. Say that $\alpha,\alpha'$ are \emph{close}, denoted as $\alpha \approx \alpha'$, if they are $\kappa$--close for some $\kappa \geq 0$.

   \begin{lemma}[Closeness of \uc natural transformations]
    Two \uc natural transformations $\alpha$, $\alpha'$ are close if and only if $\overline{\prod_j \alpha_j} = \overline{\prod_j \alpha'_j}$; in which case $\overline{\prod_\phi \alpha_{\cod \phi}} = \overline{\prod_\phi \alpha'_{\cod \phi}}$ also holds. \qed
   \end{lemma}

   It follows that $\uclim_{\sJ} \alpha$ depends only on the closeness class of $\alpha$.

   Given \uc natural transformations $\alpha \from D \sqtto D'$ and $\omega \from D' \sqtto D$ in $\cat$, say that $\omega$ is a \emph{right} (resp.~\emph{left}) \emph{inverse of $\alpha$ up to closeness} in $\cat$ if $\alpha\omega \approx 1_{D'}$ (resp.~$\omega\alpha \approx 1_{D})$; this is equivalent to the existence of some $\kappa \geq 0$ such that $\alpha_j \omega_j \approx_\kappa 1_{D'_j}$ (resp.~$\alpha_j\omega_j \approx_\kappa 1_{D_j})$ for all $j\in\sJ$. If $\alpha$ is left-invertible, right-invertible, or invertible (up to closeness) as a \uc natural transformation, then the induced morphism $\uclim_{\sJ} \alpha$ (if it exists) also satisfies the respective property in $\barcat$.
   \begin{lemma}[Right inverse]\label{lem:right-inv}
    Let $\alpha \from D \sqtto D'$ and $\omega \from D' \sqtto D$ be \uc natural transformations. Then $\overline{\prod_j \alpha_j} \circ \overline{\prod_j \omega_j} = \bar 1_{\prod_j D'_j}$ if and only if     $\omega$ is a right inverse of $\alpha$ up to closeness in $\Born$; in which case $\overline{\prod_\phi \alpha_{\cod\phi}} \circ \overline{\prod_\phi \omega_{\cod\phi}} = \bar 1_{\prod_\phi D'_{\cod\phi}}$ also holds. \qed
   \end{lemma}

    Next, we show that the existence of \uc limits persists under right-invertible \uc natural transformations.

  \begin{proposition}[Retractions and \uc limits]\label{prop:uclim_ret}
   Let $\alpha \from D \sqtto D'$ be a right-invertible \uc natural transformation (up to closeness) in $\cat$. Assume that $\uclim_{\sJ} D$ exists. Then $\uclim_{\sJ} D'$ exists and $\uclim_{\sJ} \alpha \from \uclim_{\sJ} D \to \uclim_{\sJ} D'$ is right-invertible in $\barcat$.
     \end{proposition}

  \proof
  This follows from Propositions \ref{prop:equaliser}, \ref{prop:tupstab}, Corollary \ref{cor:eq_ret}, and Lemma \ref{lem:right-inv}.
  \endproof

  \begin{theorem}[Retractions and universal $\cat$--geodesic cones]\label{thm:uc-retraction}
   Let $\alpha \from D \sqtto D'$ be a right-invertible morphism of \uc diagrams (up to closeness) in $\Born$. If $D$ admits a universal $\cat$--geodesic cone, then so does $D'$. Moreover, $\floorC{\uclim_{\sJ} \alpha} \from \floorC{\uclim_{\sJ} D} \to \floorC{\uclim_{\sJ} D'}$ is right-invertible in $\barcat$.
   \end{theorem}

  \proof
  By Theorem \ref{thm:univ-geod-cone} and Corollary \ref{cor:uc-cat-cone}, $\Eq^*(\gamma, \delta)$ stabilises to a coarsely geodesic space, hence $\uclim_{\sJ} D$ exists and is coarsely geodesic. By Proposition \ref{prop:uclim_ret}, $\uclim_{\sJ} D'$ exists and the induced morphism $\uclim_{\sJ} \alpha \from \uclim_{\sJ} D \to \uclim_{\sJ} D'$ is right-invertible in $\barBorn$. By Lemma \ref{lem:cgeod-ret}, $\uclim_{\sJ} D'$ is coarsely geodesic, hence $\floorC{\uclim_{\sJ} D'}$ exists and realises the universal $\cat$--geodesic cone for $D'$. The result follows by applying the Rips colimit functor $\floorC{-}$.
  \endproof

  The following is useful if we wish to keep track of constants.

  \begin{proposition}[Rips-tuple retraction]\label{prop:rips-tuple-ret}
   Let $\alpha \from D \sqtto D'$ and $\omega \from D \sqtto D'$ be \uc natural transformations in $\cat$ satisfying $\alpha\omega \approx_K 1_{D'}$ for some $K \geq 0$. Let $\rho$ be a uniform upper control for all $\alpha_j$. Suppose that $\floorC{\uclim_{\sJ} D}$ is realised by $\Rips_\sigma\Eq^\kappa(\gamma,\delta)$ for some $\sigma, \kappa \geq 0$. Then $\floorC{\uclim_{\sJ} D'}$ is realised by $\Rips_{\sigma'}\Eq^{\kappa'}(\gamma',\delta')$ where $\kappa' = 2K + \rho\kappa$ and $\sigma' = \max(\rho\sigma, K)$. Moreover, the canonical morphism $\Rips_\sigma\Eq^\kappa(\gamma,\delta) \to \Rips_{\sigma'}\Eq^{\kappa'}(\gamma',\delta')$ is induced by $\prod_j \alpha_j$; and this representative is 1--Lipschitz.
  \end{proposition}

  \proof
  The filtration $\Eq^*(\gamma,\delta)$ stabilises at threshold $\kappa$ by Lemma \ref{lem:eq-threshold}. Note that $\prod_j \alpha_j$ has upper control $\rho$ and $\prod_j \alpha_j \circ \prod_j \omega_j \approx_K 1_{\prod D_j}$. The constant $\kappa'$ is obtained by applying Lemmas \ref{lem:nat-equal} and \ref{lem:ret-stab} to the diagram of product maps and spaces associated to $D$ (appearing above Lemma \ref{lem:nat-trans}). By Lemma \ref{lem:nat-equal}, the  map $\Eq^\kappa(\gamma,\delta) \to \Eq^{\kappa'}(\gamma',\delta')$ induced by $\prod_j \alpha_j$ has upper control $\rho$. This constant $\sigma'$ is chosen as in the proof of Lemma \ref{lem:cgeod-ret}; in particular, the choice $\sigma' \geq \rho\sigma$ implies the desired 1--Lipschitz property.
  \endproof

 \section{Hierarchically hyperbolic spaces}\label{sec:hhs}

 \subsection{HHS as universal quasigeodesic cones}

    In this section, we explain how the data of a hierarchically hyperbolic space (HHS) together with its family of projection maps to hyperbolic spaces can be encoded using \uc diagrams. Let us briefly recall some essential ingredients for the HHS--setup; full details can be found in \cite{HHS1, HHS2, HHSsurvey}. Assume we are given the following data:

\begin{itemize}
  \item A family of uniformly hyperbolic spaces $\{\C U\}_U$, called \emph{factor spaces}, indexed by a set of \emph{domains} $U \in \calS$.
  \item A \emph{consistency constraint} for each unordered pair of distinct domains $U,V \in \calS$, specified by a prescribed subset $\calR_{UV} \subseteq \C U \times \C V$.   \item A quasigeodesic space $\X$, known as the \emph{total space}.
  \item A uniformly coarsely Lipschitz family of maps $\lambda_U \from \X \to \C U$ for $U \in \calS$ such that ${\lambda_U \times \lambda_V \from \X \to \C U \times \C V}$ factors through $\calR_{UV}$ for all $U,V$ distinct.
  \end{itemize}

    The exact nature of the constraints $\calR_{UV}$ depend on the type of relationship $U$ and $V$ satisfy. Specifically, there is a trichotomy of binary relations on $\calS$ known as \emph{nesting} $\sqsubseteq$, \emph{transversality} $\pitchfork$, or \emph{orthogonality} $\perp$. Exactly one of these relations hold for each $U,V \in \calS$ distinct. The sets $\calR_{UV}$ comprise pairs $(x_U, x_V)$ satisfying a certain inequality up to a fixed uniform error: when $U \sqsubseteq V$ (or $V \sqsubseteq U$), this is the \emph{Bounded Geodesic Image Inequality}; for $U \pitchfork V$, it is the \emph{Behrstock Inequality}; otherwise, for $U \perp V$ there is no constraint, in which case $\calR_{UV} = \C U \times \C V$.

    Let us convert encode the family of factor spaces and consistency constraints using a diagram. Define an indexing category $\sJ$ and a diagram $D \from \sJ \tto \Born$ as follows.
    \begin{itemize}
    \item Objects of $\sJ$: one \emph{primary} vertex $U$ for each $U \in \calS$; and one \emph{secondary} vertex $UV$ for each (unordered) pair of distinct $U , V \in \calS$   (choose either one of $UV$ or $VU$ as the label).
    \item (Non-identity) arrows of $\sJ$: two arrows $UV \to U$ and $UV \to V$ for  each secondary $UV$.
    \item Vertex spaces: \emph{primary spaces} $D_U:= \C U$ for each primary $U$; and \emph{secondary spaces} $D_{UV}:= \calR_{UV}$ with the $\ell^\infty$--metric for each secondary $UV$.
    \item Bonding maps: for each arrow $UV\to U$ define $\beta^{UV}_{U}\from  \calR_{UV}\to\mathcal C U$ as the projection $(x_U,x_V) \mapsto x_U$ to the $\C U$--factor; similarly, $\beta^{UV}_{V}$ is the projection to the $\C V$--factor.
    \end{itemize}
    
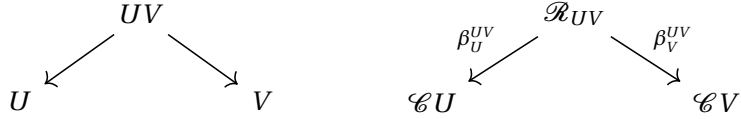
\begin{center}
    \begin{tikzcd}[column sep=normal,row sep=normal]
        & UV  \arrow[dl] \arrow[dr]  &   \\ U & & V \end{tikzcd}
    \qquad\qquad
        \begin{tikzcd}[column sep=normal,row sep=normal]
        & \calR_{UV}  \arrow[dl, "\beta^{UV}_U" swap] \arrow[dr, "\beta^{UV}_V"]  &   \\ \C U & & \C V \end{tikzcd}
        \end{center}

    As all bonding maps are $1$--Lipschitz, $D$ may be regarded as a \uc diagram. We now verify that the total space $\X$ together with its family of maps to factor spaces forms a \uc cone over $D$.

    \begin{lemma}[HHS--conversion]\label{lem:convert}
    Let $D \from \sJ \tto \Born$ be a diagram encoding an HHS--structure as defined above. Then there exists a quasigeodesic \uc cone $\lambda \from \X \sqtto D$ whose legs are
    \begin{itemize}
    \item $\lambda_U \from \X \to \C U$ as given in the HHS--data for each primary vertex $U$; and
    \item $ \lambda_{UV}:= \lambda_U \times \lambda_V \from \X\to \calR_{UV}$ for each secondary vertex $UV$.
    \end{itemize}
    \end{lemma}

    \proof
    All legs $\lambda_U$ (for primary vertices $U$) are uniformly coarsely Lipschitz by assumption, hence the same is true of $\lambda_{UV}$ for all secondary vertices $UV$. The consistency constraints ensure that $\lambda_{UV}$ indeed takes values in $\calR_{UV}$; moreover, compositions with bonding maps yield
    $\beta^{UV}_{U}\circ\lambda_{UV}=\lambda_U$ and $\beta^{UV}_{V}\circ\lambda_{UV}=\lambda_V$.
    \endproof

    \begin{remark}
    This conversion in fact yields a true cone $\lambda \from \X \tto D$ in $\Born$. Indeed, there are no non-trivial compositions involved.
    \end{remark}

    We wish to prove that $\lambda \from \X \sqtto D$ is a universal quasigeodesic \uc cone. Consider the codomain and domain maps $\gamma, \delta$ associated to $D$ as defined in Section \ref{sec:uc-cone}. These are both 1--Lipschitz since all bonding maps of $D$ are 1--Lipschitz.    Observe that a tuple
    \[(x_U)_U \times (x_{UV})_{UV} ~\in~ \prod_{U} \C U \;\times\; \prod_{UV} \calR_{UV} \;=\; \prod_{j \in \Obj(\sJ)} D_j \]
    belongs to $\Eq^\kappa(\gamma,\delta)$ if and only if
    \[d_{\C U}\left(x_U, \beta^{UV}_{U} x_{UV}\right) \leq \kappa \quad\textrm{and}\quad d_{\C V}\left(x_V, \beta^{UV}_{V} x_{UV}\right) \leq \kappa \]
    for every secondary vertex $UV$.

    To verify that the equaliser filtration $\Eq^*(\gamma,\delta)$ stabilises,we require one important property of HHS's: the Realisation Theorem. Roughly speaking, it says that the image of the product map $\prod_U \lambda_U \from \X \to \prod_U \C U$ is coarsely determined by all the pairwise constraints. This was originally proven by Behrstock--Kleiner--Minsky--Mosher in the case of the mapping class group equipped with its subsurface projections to curve graphs \cite{BKMM12}.

    \begin{theorem}[Realisation Theorem {\cite[Theorem 3.1]{HHS2}}]
    There exists an increasing function $r \from [0,\infty) \to [0,\infty)$ with the following property. Let $(x_U)_{U} \in \prod_{U\in\calS} \C U$ be a tuple satisfying the HHS consistency constraints up to error $\kappa\geq 0$, that is, for each distinct pair $U,V \in \calS$, the point $(x_U, x_V)$ lies in the closed $\kappa$--neighbourhood of $\calR_{UV}$ in $\C U \times \C V$. Then there exists $p\in \X$ such that $d_{\C U}(\lambda_U(p),x_U)\leq r(\kappa)$ for every $U \in \calS$. \qed
    \end{theorem}

    The Realisation Theorem concerns tuples in the product $\prod_{U\in\calS} \C U$ of the primary spaces. Realisation also holds for tuples in the product $\prod_{U} \C U \times \prod_{UV} \calR_{UV}$ of all primary and secondary spaces, at the cost of adjusting the error bound.

    \begin{lemma}\label{lem:dense-image}
     The image of ~$\prod_U \lambda_U \times \prod_{UV} \lambda_{UV}$ is a coarsely dense subset of $\Eq^\kappa(\gamma,\delta)$ for all $\kappa \geq 0$. Consequently, the inclusions $\Eq^0(\gamma,\delta) \hookrightarrow \Eq^\kappa(\gamma,\delta)$ are coarsely surjective.
    \end{lemma}

    \proof
    Suppose that $x:= (x_U)_U \times (x_{UV})_{UV} \in \Eq^\kappa(\gamma, \delta)$. This means that
    \[d_{\C U}\left(x_U, \beta^{UV}_{U} x_{UV}\right) \leq \kappa \quad\textrm{and}\quad d_{\C V}\left(x_V, \beta^{UV}_{V} x_{UV}\right) \leq \kappa \]
    for every secondary vertex $UV$. By the Realisation Theorem, there exists some $p \in \X$ such that $d_{\C U}(\lambda_U(p),x_U)\leq r(\kappa)$ for every $U \in \calS$. Now $y:=(\lambda_U p)_U \times (\lambda_{UV} p)_{UV} \in \Eq^0(f,g)$ by the Consistency Constraints. Then
    \[d_{\C U}\left(\beta^{UV}_{U}\lambda_{UV} p, \beta^{UV}_{U} x_{UV}\right) \leq d_{\C U}\left(\lambda_{U} p, x_U\right) + d_{\C U}\left(x_U, \beta^{UV}_{U} x_{UV}\right) \leq r(\kappa) + \kappa\]
    for every arrow $UV \to U$, hence $d_{\C U \times \C V}\left(\lambda_{UV} p,  x_{UV}\right) \leq r(\kappa) + \kappa$ for every secondary vertex $UV$.
    Therefore, $d_{\prod_{U} \C U \times \prod_{UV} \calR_{UV}} (x,y) \leq r(\kappa) + \kappa$.
    \endproof

    Next, we consider the Uniqueness Criterion for HHS. This is itself one of the HHS axioms, however, we shall state a version which allows for arbitrary \uc cones over $D$. The \emph{Uniqueness Criterion} for a family of maps $\mu_U \from W \to \C U$ for $U \in \calS$ from a metric space $W$ asserts:
    \begin{itemize}
     \item there exists a nondecreasing function $\theta\colon[0,\infty)\to[0,\infty)$ such that whenever $w,w' \in W$ satisfy $d_{\C U}(\mu_U(w),\mu_U(w'))\leq R$ for all $U \in \calS$, where $R \geq 0$, then $d_{w}(w,w') \leq \theta(R)$.
    \end{itemize}
    This is equivalent to the product map $\prod_U \mu_U \from W \to \prod_U \C U$ admitting a lower control. To reformulate this in terms of $\prod_U \mu_U \times \prod_{UV} \mu_{UV}$, we use the following observation.

    \begin{lemma}\label{lem:section}
     The 1--Lipschitz projection map $\prod_{U} \C U \times \prod_{UV} (\C U \times \C V)\to \prod_{U} \C U$ admits a 1--Lipschitz section $s$ given by $(y_U)_U \mapsto (y_U)_U \times (y_U, y_V)_{UV}$. In particular, $s$ is an isometric embedding. \qed
    \end{lemma}

    \begin{corollary}\label{cor:prod_ce}
     A family of maps $\mu_U \from W \to \C U$ for $U \in \C U$ satisfies the Uniqueness criterion if and only if $\prod_U \mu_U \times \prod_{UV} \mu_{UV} = s \circ \prod_U \mu_U$ admits a lower control. \qed
    \end{corollary}

    We now give a proof of Theorem \ref{thm:hhs-universal}.

    \begin{theorem}[HHS as universal quasigeodesic \uc cone]\label{thm:hhs-uccone}
    The \uc cone $\lambda \from \mathcal{X} \sqtto D$ as given in Lemma \ref{lem:convert} is a universal quasigeodesic \uc cone over $D$.
    \end{theorem}

    \proof
    By Proposition \ref{prop:equaliser} and Lemma \ref{lem:dense-image}, $\Eq^0(\gamma,\delta)$ serves as an equaliser for the pair $\bar\gamma, \bar\delta$ in $\barBorn$. The Uniqueness criterion holds for the family $\lambda_U$ for $U \in \calS$ as this is one of the HHS axioms. Appealing to Lemma \ref{lem:dense-image} and Corollary \ref{cor:prod_ce}, $\prod_U \lambda_U \times \prod_{UV} \lambda_{UV}$  factors through a coarse equivalence $\X \to \Eq^0(\gamma,\delta)$, and so this product map also realises the equaliser of $\bar\gamma, \bar\delta$ in $\barBorn$. Since $\X$ is quasigeodesic, by Corollary \ref{cor:uc-cat-cone}, it follows that $\lambda \from \mathcal{X} \sqtto D$ is a universal quasigeodesic \uc cone over $D$
    \endproof

    Appealing to Theorem \ref{thm:univ-geod-cone}, there exists a canonical isomorphism $\X \cong \Rips_\sigma\Eq^0(\gamma,\delta)$ in $\barCL$ for $\sigma \geq 0$ large. However, the space $\Eq^0(\gamma,\delta)$ has some redundancy in the secondary co-ordinates: $(x_U)_U \times (x_{UV})_{UV} \in \Eq^0(\gamma, \delta)$ if and only if $x_{UV} = (x_U, x_V)$ for every secondary vertex $UV$. In other words, the projection map $\prod_{U} \C U \times \prod_{UV} \calR_{UV} \to \prod_{U} \C U$ restricts to a bijection from $\Eq^0(\gamma,\delta)$ to its image
     \[\widehat{\Eq^0}(\gamma,\delta) := \left\{(x_U)_U \in \prod_{U} \C U \st (x_U, x_V) \in \calR_{UV} \textrm{ for all secondary } UV\right\}.\]
     By Lemma \ref{lem:section}, this bijection yields an isometry. Consequently, $\X$ may be recovered up to quasi-isometry by taking the stable Rips graph of the space of consistent tuples taken over the primary factors.

    \begin{proposition}
    For $\sigma \geq 0$ sufficiently large, there is a quasi-isometry $\X \to \Rips_\sigma\widehat{\Eq^0}(\gamma,\delta)$ coinciding with the underlying product map $\prod_U \lambda_U$. \qed
    \end{proposition}

    Since the image of $\prod_U \lambda_U$ is coarsely dense in $\widehat{\Eq^0}(\gamma,\delta)$, Theorem \ref{thm:hhs-image} follows.

    \subsection{Pairwise compatible retractions}\label{sec:hier-retract}

    In this section, we prove that pairwise compatible factorwise retractions induce a canonical retraction of an HHS to a subspace. The argument uses only the universal property and therefore applies in greater generality: the same proof works for any family of metric spaces $\{\C U\}_{U \in \calS}$ equipped with pairwise constraint sets $\calR_{UV} \subseteq \C U \times \C V$ for each unordered distinct pair $U,V\in\calS$ (if there is no constraint for a given pair $U,V$, set $\calR_{UV} = \C U \times \C V$). We denote this data as
    \[(\C ; \calR )_\calS := ( \{\C U\}_U ; \{\calR_{UV}\}_{UV})_\calS,\] and refer to it as a \emph{pairwise constrained family} of metric spaces. The notation is intended to evoke a generators-subject-to-relators style intuition.

    A pairwise constrained family $(\C ; \calR )_\calS$ can be encoded by a diagram $D \from \sJ \tto \Born$ in the same manner as for HHS--data, where we take secondary objects $UV$ over all unordered distinct pairs.    A \emph{total space} for $(\C ; \calR )_\calS$ is defined as a universal $\cat$--geodesic \uc cone over $D$.  In particular, any HHS is a total space over its HHS--data (factor spaces and pairwise constraints) working in $\cat = \CLip$.

    Let us consider families of maps respecting the pairwise constraints up to bounded error.

    \begin{definition}[Pairwise compatible]
    Let $(\C ; \calR )_\calS$ and $(\C' ; \calR' )_\calS$ be pairwise constrained families indexed by $\calS$. A uniformly controlled family of maps $\{ \alpha_U \from \C U \to \C' U\}_{U \in \calS}$ is \emph{pairwise compatible} if there exists a \emph{compatibility constant} $r \geq 0$ such that $(\alpha_U \times \alpha_V)(\calR_{UV})$ lies in the (closed) $r$--neighbourhood of $\calR'_{UV}$ in $\C' U \times \C' V$ for all secondary $UV$.
    \end{definition}

    Assume such a family $\{ \alpha_U \from \C U \to \C' U\}_{U \in \calS}$ admits a uniform upper control $\rho$ and compatibility constant $r \geq 0$. We encode this using a \uc natural transformation as follows. Let $D' \from \sJ \tto \Born$ be the \uc diagram encoding the family $(\C' ; \calR' )_\calS$, with the bonding maps (projections to each factor) denoted as $(\beta')^{UV}_U$, $(\beta')^{UV}_V$. The components $\alpha_{UV} \from \calR_{UV} \to \calR'_{UV}$ for each secondary $UV$ are defined by setting $\alpha_{UV}(x_U, x_V)\in \calR'_{UV}$ to be any point within distance $r$ of $(\alpha_U x_U, \alpha_V x_V)$ in $\C' U \times \C' V$. Note that any two such choices lie within distance $2r$ of one another, so the family $\alpha_{UV}$ is well-defined up to uniform closeness.

   \begin{lemma}[Pairwise compatible conversion]\label{lem:pc-convert}
    Let $\{\alpha_U \from \C U \to \C' U\}_{U}$ be a pairwise compatible family of uniformly controlled maps from $(\C ; \calR )_\calS$ to $(\C' ; \calR' )_\calS$.     Then there is a \uc natural transformation $\alpha \from D \sqtto D'$ whose components are the $\alpha_U$ for all primary $U$ and the $\alpha_{UV}$ for all secondary $UV$.     Moreover, all components have upper control $\rho'(t) := \rho(t) + 2r$.
   \end{lemma}

   \proof
   Let $\rho$ be a uniform upper control for all $\alpha_U$ and $r \geq 0$ be a compatibility constant.
   For any $(x_U, x_V) \in \calR_{UV}$, we have that
   \[\alpha_{UV}(x_U, x_V) \approx_r (\alpha_U x_U, \alpha_V x_V) = \left (\alpha_U \beta^{UV}_U (x_U, x_V),  \alpha_V \beta^{UV}_V (x_U, x_V) \right),\]
   hence $(\beta')^{UV}_U \alpha_{UV} \approx_r \alpha_U \beta^{UV}_U$ and $(\beta')^{UV}_V \alpha_{UV} \approx_r \alpha_V \beta^{UV}_V$.
   Furthermore, $\rho'(t) := \rho(t) + 2r$ serves as an upper control for all $\alpha_{UV}$.
      \endproof

   Note that the closeness class of $\alpha \from D \sqtto D'$ is completely determined by its primary components $\{\alpha_U \from \C U \to \C' U\}_{U}$. This allows us to characterise right-invertibility of $\alpha$ (up to closeness) in terms of uniform retractions defined on the primary factors.

   \begin{lemma}[Compatible right inverse]\label{lem:primary-ret}
    Let $\{\alpha_U \from \C U \to \C' U\}_{U}$ and $\{ \omega_U \from \C' U \to \C U\}_{U \in \calS}$ be pairwise compatible uniformly controlled families of maps between $(\C ; \calR )_\calS$ and $(\C' ; \calR' )_\calS$. Assume that $c \geq 0$ is a compatibility constant for the $\omega_U$ family.  If there exists $K \geq 0$ such that $\alpha_U\omega_U \approx_K 1_{\C' U}$ for all $U$ then the induced \uc natural transformations $\alpha \from D \sqtto D'$ and $\omega \from D' \sqtto D$ satisfy $\alpha\omega \approx_{K'} 1_{D'}$, where $K' := K + \rho c + r$. The converse also holds.
    \end{lemma}

   \proof
   Assume $\rho$ serves as a uniform upper control for all $\alpha_U$, and that $r \geq 0$ serves as a compatibility constant for both families of the $\alpha_U$ and $\omega_U$. By assumption, $\alpha_U\omega_U \approx_K 1_{\C' U} $ for all primary $U$. For each secondary $UV$, consider $(x'_U, x'_V) \in \calR'_{UV}$ and
   \[(y_U, y_V) := \omega_{UV}(x'_U, x'_V) \approx_c (\omega_U x'_U, \omega_V x'_V).\]
   Then
   \[\alpha_{UV} \omega_{UV}(x'_U, x'_V) = \alpha_{UV}(y_U, y_V) \approx_r (\alpha_U y_U, \alpha_V y_V) \approx_{\rho c} (\alpha_U\omega_U x'_U, \alpha_V\omega_V x'_V) \approx_K (x'_U, x'_V),\]
   hence $\alpha_{UV} \omega_{UV} \approx_{K + \rho c + r} 1_{\calR'_{UV}}$.
   The converse is immediate.
   \endproof
   
   The following is a direct application of Theorem \ref{thm:uc-retraction}.

   \begin{theorem}[Compatible retractions]\label{thm:hier-retract}
    Let $\{\alpha_U \from \C U \to \C' U\}_{U}$ and $\{ \omega_U \from \C' U \to \C U\}_{U \in \calS}$ be pairwise compatible uniformly controlled families of maps between $(\C ; \calR )_\calS$ and $(\C' ; \calR' )_\calS$ satisfying $\alpha\omega \approx 1_{D'}$. Assume that $(\C ; \calR )_\calS$ admits a total space $\X$. Then $(\C' ; \calR' )_\calS$ admits a total space $\X'$, and the canonical morphisms $\tilde \alpha \from \X \to \X'$ and $\tilde\omega \from \X' \to \X$ in $\barcat$ satisfy $\tilde\alpha \tilde\omega = \bar 1 _{\X'}$. \qed
    \end{theorem}

     The situation can be summarised using the diagram in $\barcat$ presented below.  Here, $\gamma', \delta'$ denote the respective codomain and domain maps for the \uc diagram $D'$. The constants $\sigma,\kappa,\sigma',\kappa' \geq 0$ are chosen so that all horizontal arrows are canonical isomorphisms, as granted by Theorem \ref{thm:univ-geod-cone}.     The vertical arrows in the left- and right-hand columns are obtained by via composing the corresponding vertical arrows in the central column with appropriate horizontal isomorphisms.

    \begin{center}
   \begin{tikzcd}[column sep=huge, row sep=huge] \X \arrow[r, leftrightarrow, "\sim"] \arrow[d, "\tilde{\alpha}", bend left = 15] \arrow[r, leftrightarrow, "\sim"] & \floor{\uclim_{\sJ} D}_{\cat} \arrow[d, "\floor{\uclim_{\sJ} \alpha}_{\cat}", bend left = 15] \arrow[r, leftrightarrow, "\sim"] & \Rips_\sigma{\Eq}^\kappa(\gamma, \delta)  \arrow[d, "\hat{\alpha}", bend left = 15] \\ \X' \arrow[r, "\sim", leftrightarrow] \arrow[u, "\tilde{\omega}", bend left = 15] & \floor{\uclim_{\sJ} D'}_{\cat} \arrow[u, "\floor{\uclim_{\sJ} \omega}_{\cat}", bend left = 15] \arrow[r, leftrightarrow, "\sim"] & \Rips_{\sigma'}{\Eq}^{\kappa'}(\gamma', \delta')  \arrow[u, "\hat{\omega}", bend left = 15] \end{tikzcd}
  \end{center}

  \begin{remark}
    If constants $\sigma, \kappa$ suffice for the upper row then choosing \[\sigma' = \max(\rho'\sigma, K') = \max(\rho\sigma + 2r, K + \rho c + r) \quad\textrm{and}\quad \kappa' = 2K' + \rho'\kappa = 2K + 2\rho c + \rho\kappa + 4r\]
    will work for the lower. This follows from Proposition \ref{prop:rips-tuple-ret} and Lemmas \ref{lem:pc-convert} and \ref{lem:primary-ret}.
  \end{remark}

   By Lemma \ref{lem:left-inverse}, the relation $\tilde \alpha \tilde \omega = \bar 1_{\X'}$ implies that $\tilde \omega \from \X' \to \X$ is represented by a map with lower control in $\F^T$ whose image is a $\cat$--retract. Up to isomorphism in $\barcat$, we may realise $\X'$ as an isometrically embedded subspace of $\X$; in this case, $\tilde \omega$ is represented by the inclusion $\X' \hookrightarrow \X$ and $\tilde \alpha$ by a retraction onto $\X'$.
   The morphism $\hat\alpha$ on underlying tuple-spaces has an explicit factorwise description: map $(x_U)_U\times (x_U, x_V)_{UV} \in {\Eq}^\kappa(\gamma, \delta)$ to a point in $(y_U)_U\times (y_U, y_V)_{UV} \in {\Eq}^{\kappa'}(\gamma', \delta')$ which lies within a uniformly bounded distance of $(\alpha_U x_U)_U \times (\alpha_U x_U, \alpha_V x_V)_{UV}$. Up to closeness, this yields the canonical morphism granted by Lemma \ref{lem:ret-stab}. In particular, $\hat\alpha$ is the unique closeness class whose restriction to the primary factors agrees with the restriction of $\prod_U \alpha_U$. The morphism $\hat \omega$ can be described analogously.

   Finally, we apply Lemma \ref{lem:primary-ret} and Theorem \ref{thm:hier-retract} to the case where the $\alpha_U$ are uniform retractions to subspaces. Theorem \ref{thm:hier-retractions} then follows as a special case.

   \begin{corollary}[Compatible retraction to subspace]\label{cor:hier-ret-subspace}
   Let $\X$ be a total space for $(\C ; \calR )_\calS$. Suppose that $\alpha_U \from \C U \to \C' U$ for $U \in \calS$ is a family of maps to subspaces $\C' U \subseteq \C U$ admitting a common upper control $\rho \in \F$.    Assume the following conditions hold.
      \begin{itemize}
    \item \emph{uniform retractions}: there exists $K\geq 0$ such that for all $U$, we have that $\alpha_Ux_U \approx_K x_U$ for all $x_U \in \C' U$, and
    \item \emph{uniform compatibility}: there exists $c \geq 0$ such that for all $UV$, the set $(\alpha_U \times \alpha_V)(\calR_{UV})$ lies in the $c$--neighbourhood of $\calR_{UV}$ in $\C U \times \C V$.
   \end{itemize}
    Then there is a canonical retraction $\tilde\alpha \from \X \to \X'$ in $\barcat$ to a subspace $\X' \subseteq \X$ induced by the $\alpha_U$. Moreover, $\X'$ realises the total space of $(\C' ; \calR' )_\calS$ where $\calR'_{UV} := (\alpha_U \times \alpha_V)(\calR_{UV})$.
   \end{corollary}

   \proof
   Take each $\omega_U \from \C' U \hookrightarrow \C U$ to be the inclusion.      Then $c$ is a compatibility constant for the inclusions, while the $\alpha_U$ family has compatibility constant zero.
   \endproof

\providecommand{\bysame}{\leavevmode\hbox to3em{\hrulefill}\thinspace}
\providecommand{\MR}{\relax\ifhmode\unskip\space\fi MR }
\providecommand{\MRhref}[2]{  \href{http://www.ams.org/mathscinet-getitem?mr=#1}{#2}
}
\providecommand{\href}[2]{#2}

\end{document}